%% file: draft.tex
\documentclass[12pt]{article}
\usepackage{amsmath, amsthm, amssymb, amsfonts}
\usepackage{graphicx}
\usepackage{mathtools}
\usepackage{blindtext}

\newtheorem{theorem}{Theorem} 

\newtheorem{lemma}{Lemma}
\newtheorem{proposition}{Proposition}
\newtheorem{corollary}{Corollary}

\title{Geodesics in 3-dimensional Euclidean Space with One or Two Analytic Obstacles}
\author{Chengcheng Yang\\
(communited with Prof. Robert Hardt)}

\begin{document}
\maketitle
\include{abstract}
\include{introduction}

\include{part3}

\include{part4}

\include{part7}
\include{conclusion}

\include{cite}
\end{document}

%% file: abstract.tex
\begin{abstract}
In many singular metric spaces, the regularity of a shortest-length curve is unknown.  Algebraic varieties, or more generally sets defined by finitely many polynomial or real analytic equalities or inequalities, all locally partition into finitely many real analytic submanifolds called strata.   So any component of a shortest-length curve which lies completely in one such stratum is a geodesic in the stratum, hence an embedded real analytic curve. The key question thus is whether there are only finitely many components. F. Albrecht and I.D. Berg proved this is true for a geodesic in a closed region  of $n$-dimensional Euclidean space with a smooth real analytic hyper surface as boundary. 
 Here the  curve consists of finitely many interior line segments alternating with boundary hypersurface geodesics. Their bound on the number of these depended on the initial velocity of the geodesic, and they conjectured that is independent. Here we prove this independence in $\mathbb{R}^3$. We also generalize their result to regions whose boundary is locally the boundary of the union of two transversally intersecting analytic hypersurfaces in $\mathbb{R}^3$.



\end{abstract}

%% file: introduction.tex
\begin{section}{introduction}
In many singular metric spaces, the regularity of a shortest-length curve is unknown.  Algebraic varieties, or more generally sets defined by finitely many polynomial or real analytic equalities or inequalities, all locally partition into finitely many real analytic submanifolds called strata.   So any component of a shortest-length curve which lies completely in one such stratum is a geodesic in the stratum, hence an embedded real analytic curve. The key question thus is whether there are only finitely many components. Our Ph.D. thesis verified this finiteness for shortest-length curves in any semi-algebraic subset of $\bf R^2$. Earlier in \cite{AB}, F. Albrecht and I.D. Berg proved this is true for a geodesic in a closed region  of $n$-dimensional Euclidean space with a smooth real analytic hypersurface as boundary. Here the  curve consists of finitely many interior line segments alternating with boundary hypersurface geodesics. Their bound on the number of these depended on the initial velocity of the geodesic, and they conjectured that is independent. Here we prove this independence in $\mathbb{R}^3$. We also generalize their result to regions whose boundary is locally the boundary of the union of two transversally intersecting analytic hypersurfaces in $\mathbb{R}^3$.

Let's start with the properties of real semi-algebraic sets. A semi-algebraic set in $\mathbb{R}^n$ is defined as a finite union of sets of the form $\cap_{j=1}^n \{x \in \mathbb{R}^n | f_i(x) = 0, g_i(x) > 0\}$, where $f_i$, $g_i$ are polynomials. It is easy to notice that a real semi-algebraic set generalizes a real algebraic variety. 
The triangulability question for algebraic and semi-algebraic has been well-studied since van de Waerden in 1929 \cite{BL}. After him, Whitney, Lojasiewicz, Hironaka, Hardt etc. have all proved various stratification theorems for semi-algebraic sets \cite{W}\cite{L}\cite{H}\cite{RH}. Following their footsteps, we wonder where there is a way to stratify any semi-algebraic set such that it also plays nicely with geodesics. More precisely, let $X$ be a semi-algebraic set in $\mathbb{R}^n$. A geodesic on $X$ is defined to be a locally shortest path in $X$. Given any two points $A$, $B$ in $X$, a shortest-length path from $A$ to $B$ is a geodesic whose length is minimal among all such possible geodesics between $A$ and $B$. Then we ask whether there exists some stratification $\mathcal{A}$ such that the intersection of any stratum with any shortest-length path has at most finitely many components. We call this {\it the finiteness property} of $\mathcal{A}$. Intuitively, we don't want a geodesic to oscillate back and forth infinitely often near any point. In the plane, one can show that there exists a cell decomposition satisfying the finiteness property \cite{Y}. However, in a 3-space research is till going on, because the Riemannian obstacle problem comes into the picture. There is a list of papers for references, for example, \cite{A} \cite{B}. In particular, we are interested in obstacles with analytic boundaries, since a semi-algebraic set is definited by polynomial equations which can be described locally by finitely many real-analytic equations \cite{AB} \cite{ABB}. 

In the paper \cite{AB} F. Albrecht and I.D. Berg proved that a geodesic cannot have an accumulation of interior line segments in $\mathbb{R}^n$ or $\mathbb{C}^n$ with an analytic obstacle. More precisely, let $M$ be the closure of the complement of an obstacle in an Euclidean space and let $S$ be its boundary. A geodesic in $M$ (thought of as a string stretching over some obstacle) is a locally shortest path consisting of two types of segments:{\it bounary segment} (touching the boundary with an acceleration outward normal to $S$) and {\it interior segment} or {\it interval} (not touching the boundary with zero acceleration and this is actually a line). The point connecting a boundary segment and an interior line segment is called a {\it switch point}. The regularity of a geodesic has been proved to be necessarily $C^1$ at the switch points and can be parametrized by arc length; furthermore the acceleration exists everywhere except at the switch points \cite{ABB}. When the boundary $S$ is only assumed to be $C^{\infty}$, one can show that the switch points can accumulate to a point, called {\it an intermittent point}. However, when $S$ is analytic, F. Albrecht and I.D. Berg showed that the switch points do not accumulate. That is to say, there exists an $\epsilon > 0$ such that $\gamma(s)$ has no switch point for $0 < s < \epsilon$.

In the first part of our paper we generalize their result by looking at the union of two obstacles in a 3-dimensional Euclidean space. Let $M_1$, $M_2$ be the closure of the complement of two obstacles in an Euclidean space whose surfaces are $S_1$, $S_2$, respectively. Assume that $S_1$ and $S_2$ intersect transversally. Let $M$ be the intersection of $M_1$ and $M_2$ which is the closure of the complement of the union of the two obstacles. Suppose $\gamma$ is a geodesic in $M$, then the same conclusion holds so $\gamma$ does not bounce back and forth between two surfaces infinitely often locally. In other words, $\gamma$ is eventually a boundary segment in one of the surfaces or a line segment. 

The proof consists of two parts. The first part is concerned with the case of $\mathbb{R}^3$ and the angle between $S_1$ and $S_2$ is less than $90^\circ$. The argument uses symmetry. The second part is still dealing with the case of $\mathbb{R}^3$ but the angle can be more than $90^\circ$. Here the symmetry arguemnt in the previous part fails, so we need to come up with an asymmetric argument. Lastly $\mathbb{R}$ can be replaced by $\mathbb{C}$ but for visualization we stay with the real Euclidean space. 

The second half of our paper concerns with the conjecture proposed by  F. Albrecht and I.D. Berg at end of the paper \cite{AB}. Let $M$ be the same as above, and fix one point $p$ on $M$. The conjecture is there exists a uniform bound on the number of intervals (or switch points) within a neighborhood of $p$. We are able to prove this conjecture if the dimension is 3. Namely, there exists an $\epsilon$ such that for any geodesic $\gamma$ initiating from $p$ in $M$, there are at most two intervals within the $\epsilon$-ball of $p$. 

The proof uses the idea that there are finitely many wedges covering the $(x, y)$-plane such that within each wedge such an $\epsilon$ exists. 

The author acknowledge a debt of gratitude to her thesis advisor Professor Robert Hardt, who has given his helpful advices on a regular basis. 
\end{section}

%% file: part3.tex
\begin{section}{Part One}
\begin{theorem}\label{thm1}
Let $M_1$ and $M_2$ be 3-dimensional analytic manifolds with boundary embedded in $\mathbb{R}^3$. Denote the boundary surfaces of $M_1$ and $M_2$ by $S_1$ and $S_2$, respectively. Assume that $S_1$ and $S_2$ intersect transversally whose angle is less than $90^\circ$. Let $M$ be the intersection of $M_1$ and $M_2$. If $\gamma$ is a geodesics stretching over $M$ parametrized by arc length $s$, with $\gamma(0) = p \in S_1 \cap S_2$. Then there exists an $\epsilon > 0$ such that $\gamma$ has no switch point for $0 < s < \epsilon$. 
\end{theorem}

\begin{proof} 
1. Set the coordinate system.

Without loss of generality we may assume that $p$ is the origin, the $x$-axis is tangent to the given geodesic $\gamma$ at $p$, and the outward normal vectors to $S_1$ and $S_2$ at $p$ are $(0, -k, 1)$ and $(0, k, 1)$, respectively. So the normal vectors at $p$ are symmetric with respect to the $z$-axis. Let $ 0 < k < \frac{1}{1 + \beta}$ for some $\beta > 0$ so $k < 1$. The two surfaces $S_1$ and $S_2$ are defined near the origin by analytic equations of the form $z = g(x, y)$ and $z = h(x, y)$, respectively. It follows that the outward normal vector to $S_1$ at the origin is equal to $(-g_x(0, 0), -g_y(0, 0), 1)$, which is also equal to $(0, -k, 1)$ by hypothesis, therefore 
\[g_x (0, 0) = 0, \ \ \  g_y(0, 0) = k.\] 
Similarly, we have 
\[h_x(0, 0) = 0, \ \ \ h_y(0, 0) = -k.\] 
Since
\[ \frac{\partial}{\partial y} \displaystyle \big{|}_{(0, 0)}  [g(x, y) - h(x, y)] = 2k > 0, \]
 the inverse function theorem implies that the intersection of $S_1$ and $S_2$ near $p$ is a real analytic curve defined by the following equations:
 \[ y = \phi(x), \ \ \ z = g(x, \phi(x)) = h(x, \phi(x)),\]
where $\phi$ is real analytic and $\phi(0) = 0$, $\phi'(0) = 0$. In fact the tangent vector to $S_1 \cap S_2$ at the origin is given by 
\begin{eqnarray*}
(1, \phi'(x), g_x (x, \phi(x)) + g_y (x, \phi(x)) \phi'(x)) \big{|}_{x=0} &=& (1, \phi'(0), g_x(0, 0) + g_y(0, 0)\phi'(0)) \\ &=& (1, \phi'(0), k\phi'(0)), 
\end{eqnarray*}
which is normal to $(0, k, 1)$, thus
\[(1, \phi'(0), k\phi'(0)) \cdot (0, k, 1) = 2k \phi'(0) = 0 \Longrightarrow \phi'(0) = 0.\]
So \begin{equation}\label{0}
\phi(x) = a_M x^M + a_{M+1} x^{M+1} + \cdots, 
\end{equation}
where $M \geq 2$ and $a_M \neq 0$.
Notice that we assume $\phi(x)$ is not identically zero and $g(x, 0)$ is not identically zero. Otherwise we will have trivial cases which will be included at the end.
Thus the equation defining $S_1$ near $p$ is of the form
\begin{equation}\label{1}
g(x, y) = ky + x^N a(x, y) + xy b(x, y) + y^2 c(y), 
\end{equation}
where $N \geq 2$, the functions $a, b, c$ are analytic, and $a(0, 0) \neq 0$. Likewise assume that $h(x, 0)$ is not identically zero, the equation defining $S_2$ near $p$ is of the form
\begin{equation}\label{2}
h(x, y) = -ky + x^{\tilde{N}} \tilde{a} (x, y) + xy \tilde{b}(x, y) + y^2 \tilde{c}(y),
\end{equation}
where $\tilde{N} \geq 2$, the functions $\tilde{a}, \tilde{b}, \tilde{c}$ are analytic, and $\tilde{a}(0,0) \neq 0$.
Choose the orientation of the coordinate system so that $\gamma'(0) = (1, 0, 0)$, $M_1 = \{ z \leq g(x, y)\}$, and $M_2 = \{ z \leq h(x, y)\}$. 

2. Project the $S_1 \cap S_2$ onto the $(x, y)$-plane.

The projection of the intersection of $S_1$ and $S_2$ onto the $(x, y)$-plane is a curve given by $(x, \phi(x))$, where $x \in (-\delta, \delta)$ for some $\delta > 0$. Then it divides the vertical strip $(-\delta, \delta) \times (-\infty, \infty)$ into two disconnected regions:
\[\{x < \phi(x)\}, \ \ \ \{x > \phi(x)\}.\]
Suppose $y < 0 = \phi(0)$, then using the linear approximation 
\[g(0, y) - h(0, y) \approx g_y(0, 0)y - h_y(0, 0)y = 2ky  <0.\] Thus connectedness implies that 
\begin{equation}\label{3}
\{g(x, y) < h(x, y)\} = \{y < \phi(x)\}.
\end{equation}
Similarly, 
\begin{equation}\label{4}
\{g(x, y) > h(x, y)\} = \{y > \phi(x)\}.
\end{equation}
Given a point $z = g(x, y)$ in $S_1$. If the point lies in $M_2$, then $z \leq h(x, y)$, thus $g(x, y) \leq h(x, y)$. With (\ref{3}) it follows that the projection of $S_1 \cap M_2$ near $p$ is 
\[\{g(x, y) \leq h(x, y), -\delta < x < \delta\} = \{x \leq \phi(x) , -\delta < x < \delta\}.\]
Likewise with (\ref{4}) the projection of $S_2 \cap M_1$ near $p$ is 
\[\{g(x, y) \geq h(x, y), -\delta < x < \delta\} = \{x \geq \phi(x), -\delta < x < \delta\}.\] 
In conclusion the graph of $\phi(x)$ divides the $(x, y)$-plane into two parts near 0: the part below the graph corresponding to the projection of the surface $S_1$ in $M$ and the part above the graph corresponding to the projection of the surface $S_2$ in $M$. 

3. Concavity of $\phi$, which will become crucial in proving the theorem later.  

Since $\gamma'(0) = (1, 0, 0)$, $x'(s) > 0$ for $0 \leq s \leq \epsilon$ if we choose $\epsilon$ small enough. Therefore $x(s) > 0$ when $0 \leq s \leq \epsilon$.
According to the equation (\ref{0}), 
\[\phi''(x) = M(M-1)a_M x^{M-2} + (M+1)Ma_{M+1} x^{M-1} + \cdots.\]
When $a_M > 0$, $\phi(x)$ is concave upward over the interval $(0, \delta)$; and when $a_M < 0$, $\phi(x)$ is concave downward over the same interval.  Furthermore, we may also assume that $x(s) < \delta$ for $0 \leq s \leq \epsilon$. So there are two cases to consider: $a_M > 0$ and $a_M < 0$. 

4. Approximate $y(s)$ and $y'(s)$ using the normal vectors $N_1(s)$, $N_2(s)$ to $S_1$, $S_2$. 

We denote
\[\gamma(s) = (x(s), y(s), z(s)), \text{ for } 0\leq s \leq \epsilon.\]
If $\gamma(s) \in S_1$, the normal vector to $S_1$ at $\gamma(s)$ is 
\[N_1(s) = (-g_x(x(s), y(s)), -g_y(x(s), y(s)), 1).
\]
From (\ref{1}) it follows that 
\begin{eqnarray*}
g_x (x(s), y(s)) = x(s) m(x(s), y(s)) + y(s) b(x(s), y(s));\\
g_y(x(s), y(s)) = k + x(s)k(x(s), y(s))+ y(s) l(y(s)),
\end{eqnarray*}
where the functions $m, b, k, l$ are bounded near $(0, 0)$. 
Moreover, since $x'(s) > 0$ for $0 \leq s \leq \epsilon$, $x(s)$ has a $C^1$-inverse function $s(x)$ for $s \in [0, \epsilon]$. 
Therefore we can express $y(s)$ as 
\[y(s) = y(s(x)) = \alpha(x),\]
where $\alpha$ is a $C^1$-function and $\alpha(0) = \frac{d \alpha}{dx} (0) = 0$. Then one has $y(s) = o(x(s))$. Hence
\begin{eqnarray*}\label{g_xandg_y}
g_x (x(s), y(s)) = x(s) [m(x(s), y(s)) + \frac{y(s)}{x(s)} b(x(s), y(s))] = x(s) V_1(s);\\
g_y(x(s), y(s)) = k + x(s)[k(x(s), y(s))+ \frac{y(s)}{x(s)} l(y(s))] = k + x(s) V_2(s).
\end{eqnarray*}
Therefore
\begin{equation*}
N_1(s) = (-x(s)V_1(s), -k - x(s) V_2(s), 1),
\end{equation*}
where $V_1(s)$ and $V_2(s)$ are bounded for $s$ near 0.
Let $s$ be such that $\gamma''(s)$ exists, then $\gamma''(s) = z''(s) N_1(s)$ since the acceleration is outward normal to $S_1$. This implies that
\begin{equation}\label{5}
x''(s) = - z''(s)x(s)V_1(s), \ \ \ y''(s) = -z''(s) (k + x(s) V_2(s)).
\end{equation}
For $\epsilon$ sufficiently small, $|x(s) V_2(s)| \leq \beta k$. Therefore $|y''(s)| \leq (1+\beta)k |z''(s)|$ from the second equality in (\ref{5}).

Similarly, if $\gamma(s) \in S_2$ and $\gamma''(s)$ exists, the equality in (\ref{2}) deduces that the normal vector to $S_2$ at $\gamma(s)$ is 
\begin{equation*}
N_2(s) = (-x(s)W_1(s), k - x(s)W_2(s), 1),
\end{equation*}
where $W_1(s)$ and $W_2(s)$ are bounded for $s$ near 0. 
Then it follows from $\gamma''(s) = z''(s) N_2(s)$ that
\begin{equation}\label{6}
x''(s) = - z''(s)x(s)W_1(s), \ \ \ y''(s) = -z''(s) (-k + x(s) W_2(s)).
\end{equation}
Again for $\epsilon$ sufficiently small, one may assume that $|x(s) W_2(s)| \leq \beta k$. Therefore $|y''(s)| \leq (1+\beta)k |z''(s)|$ from the second equality in (\ref{6}). 

When $\gamma(s)$ does not lie on $S_1$ and $S_2$, $\gamma$ is a line segment so $\gamma''(s)$ is equal to 0. Combining with what we've found above, one gets $|y''(s)| \leq (1+\beta)k z''(s)$. Notice that $z(0) = z'(0) = y(0) = y'(0) = 0$. Furthermore $z''(s) \geq 0$ (and hence $z'(s) \geq 0$) on the interval $[0, \epsilon]$, because the outward normal vectors to $S_1$ and $S_2$ have a positive $z$-coordinate of 1 at the origin and $\gamma''(s)$ is directed outward on a boundary segment on $S_1$ or $S_2$. Indeed, $\gamma(s)$ is a locally shortest path and if $\gamma(s)$ lies on the surface of $M_1$ or $M_2$, its acceleration exists everywhere except at the switch points and is outward normal to the surface \cite{ABB}. 
So for $s \in [0, \epsilon]$, one can approximate
\[|y'(s)| = |\displaystyle \int_0^s y''(\sigma) d \sigma| \leq (1+\beta) k \displaystyle \int_0^s z''(\sigma) d \sigma = (1+\beta)k z'(s). \]
Integrating again one obtains
\[|y(s)| \leq (1+\beta)k z(s).\]
If $\gamma(s) \in S_1$, the equality in  (\ref{1}) gives
\begin{eqnarray*}
z(s) & = & g(x(s), y(s)) \\
& = &  ky(s) + x(s)^N a(x(s), y(s)) + x(s) y(s) b(x(s), y(s)) + y^2(s) c(y(s)) \\
& \leq & k|y(s)| + x(s) ^N |a(x(s), y(s))| + |y(s)||x(s) b(x(s), y(s)) + y(s) c(y(s))| \\
& \leq & k|y(s)| + C_1x(s)^N + C_2|y(s)| \\
& \leq & (k+C_2) (1+\beta)k z(s) + C_1x(s)^N,
\end{eqnarray*}
for some constants $C_1$, $C_2$. Since $k < \frac{1}{1+\beta}$, one can choose $\epsilon$ small enough so that $(k+C_2) (1+\beta)k < 1$. Therefore there exists a positive constant $A$ such that 
\begin{equation}\label{zyS1}
z(s) \leq A x(s)^N \Rightarrow |y(s)|  \leq (1+\beta)kA x(s)^N = B x(s)^N.
\end{equation} 

Similarly, if $\gamma(s) \in S_2$, the equality in (\ref{2}) gives us 
\begin{equation}\label{zyS2}
z(s) \leq A x(s)^{\tilde N} \Rightarrow |y(s)| \leq B x(s)^{\tilde N},
\end{equation}
by enlarging $A$ and $B$ if necessary. 

Choosing $\epsilon$ small enough so that $x(s) < 1$ and assuming without loss of generality that $N \leq \tilde{N}$, one has $x(s)^{\tilde{N}} \leq x(s)^N$. Thus with (\ref{zyS1}) and (\ref{zyS2})
\begin{equation}\label{zyS1S2}
z(s) \leq A x(s)^N,  |y(s)| \leq B x(s)^N,
\end{equation}
if $\gamma(s) \in S_1$ or $S_2$. 

Next let's approximate $y'(s)$. 
If $\gamma(s) \in S_1$, differentiating $z(s) = g(x(s), y(s))$ gives
\begin{eqnarray*}
z'(s) &=&  ky'(s) + x(s)^{N-1}[Nx'(s)a(x(s), y(s)) + x(s)(a_x(x(s), y(s))x'(s) + a_y(x(s), y(s))y'(s))] \\
	&  & + y'(s)[x(s)b(x(s), y(s)) + x(s)y(s) b_y(x(s), y(s)) + 2y(s)c(y(s)) + y^2(s)c'(y(s))]\\
	&   & + y(s)[x'(s)b(x(s), y(s)) + x(s)x'(s) b_x(x(s), y(s))] \\
	& \leq& k|y'(s)| + C_1x(s)^{N-1} + C_2|y'(s)| + C_3|y(s)| \\
	& \leq& (k + C_2) (1+\beta)k z'(s) +C_1 x(s)^{N-1} +C_3 B x(s)^N,
\end{eqnarray*}
for some constants $C_1, C_2, C_3$.
Again since $(1+\beta) k< 1$, for $\epsilon$ sufficiently small, one can make $(k + C_2) (1+\beta)k  < 1$, so 
\begin{equation}\label{z'y'S1}
z'(s) \leq C x(s)^{N-1} \Rightarrow |y'(s)| \leq (1+\beta)kG x(s)^{N-1} = D x(s)^{N-1}.
\end{equation}
Similarly, if $\gamma(s) \in S_2$, differentiating $z(s) = h(x(s), y(s))$ gives
\begin{equation}\label{z'y'S2}
z'(s) \leq C x(s)^{\tilde{N}-1} \Rightarrow |y'(s)| \leq (1+\beta)kG x(s)^{\tilde{N}-1} = D x(s)^{\tilde{N}-1}.
\end{equation}
by enlarging $C$ and $D$ if necessary. Combining (\ref{z'y'S1}) and (\ref{z'y'S2}), together with $\tilde{N} \geq N$, one has
\begin{equation}\label{z'y'S1S2}
z'(s) \leq C x(s)^{N-1},  |y'(s)| \leq (1+\beta)kG x(s)^{N-1} = D x(s)^{N-1}.
\end{equation}

Now let's look at the situation when $\gamma(s)$ is in an interior line segment. Considering a line segment in the image of $\gamma$ with two endpoints $\gamma(s_1)$ and $\gamma(s_2)$, we can parametrize $y(s)$ for $s \in [s_1, s_2]$ by 
\[y(s) = y(s_1) + T(x(s) - x(s_1)), \text{ where } T = \frac{d \alpha}{dx} (x(s_1)),\]
where with (\ref{z'y'S1S2})
\begin{equation}\label{T}
\big{|} \frac{d \alpha}{dx} (x(s)) \big{|} = \big{|} \frac{d \alpha}{dx} (x(s_1)) \big{|} = \big{|} \frac{y'(s_1)}{x'(s_2)} \big{|} \leq 2|y'(s_1)| \leq 2D x(s_1)^{N-1} \leq 2Dx(s)^{N-1}, 
\end{equation}
if $x' \geq 1/2$ by choosing $\epsilon$ small enough and the last inequality holds because $x(s)$ is increasing. 
Hence with (\ref{zyS1S2}) and (\ref{T}) one obtains
\begin{eqnarray*}
|y(s)| &\leq& |y(s_1)| + |T| (|x(s)|+|x(s_1)|) \\
		& \leq & B x(s_1)^N + 2Dx(s)^{N-1} (x(s) + x(s)) \\
		& \leq & (B+ 4D) x(s)^N.
\end{eqnarray*}
Replacing $B$ by $B + 4D$, together with (\ref{zyS1S2}),  yields that in general, 
\begin{equation}\label{zyall}
|y(s)| \leq Bx(s)^N \text{ for every } s \in [0, \epsilon]. 
\end{equation}


5. Prove $M \geq N$. 

If the geodesic $\gamma$ moves from $S_1$ to $S_2$ or from $S_2$ to $S_1$, $(x(s), y(s))$ must cross the curve $y = \phi(x)$. 
On the one hand, $|\phi(x(s))| \leq B x(s)^N$ according to (\ref{zyall}); on the other hand, $|\phi(x(s))| \geq \frac{|a_M|}{2} x(s)^M$ by (\ref{0}). Therefore one obtains the following relation:
\[\frac{|a_M|}{2} x(s)^M \leq B x(s)^N \Longrightarrow x(s)^{M-N} \leq \frac{2B}{|a_M|}.\]
Suppose $M < N$ the left-hand side converges to infinity as $s$ approaches 0, a contradiction. This means that if $M < N$ the geodesic $\gamma$ eventually stops bouncing between $S_1$ and $S_2$. Therefore it reduces to the case of one obstacle. 
Hence we proceed with $M \geq N$.


6. Next let's prove that $N = \tilde{N}$.

Case 1: $M > N$. 
On the intersection of $S_1$ and $S_2$ we've shown that $y = \phi(x)$ for $x \in (-\delta, \delta)$ and hence $g(x, \phi(x)) = h(x, \phi(x))$ over the interval $(-\delta, \delta)$. Using the equalities (\ref{0}) and (\ref{1}) one has
\begin{eqnarray*}
g(x, \phi(x)) &=& k\phi(x) + x^N a(x, \phi(x)) + x\phi(x)b(x, \phi(x)) + \phi^2(x) c(\phi(x)) \\
		   &=& kx^M(a_M + a_{M+1} x + \dots) + x^N (a(0,0) + \dots) + \\
		   &&  x^{M+1}(a_M b(0, 0) + \dots) + a_M^2 x^{2M}( c(0) + \dots).
\end{eqnarray*}
Since $M > N$ the first nonzero term in the power serious expansion of $g(x, \phi(x))$ is $a(0, 0) x^N$. Similarly the equalities (\ref{0}) and (\ref{2}) 
gives
\begin{eqnarray*}
h(x, \phi(x)) &=& -k\phi(x) + x^{\tilde{N}} \tilde{a}(x, \phi(x)) + x\phi(x)\tilde{b}(x, \phi(x)) + \phi^2(x) \tilde{c}(\phi(x)) \\
		   &=& -kx^M(a_M + a_{M+1} x + \dots) + x^{\tilde{N}} (\tilde{a}(0,0) + \dots) + \\
		   &&  x^{M+1}(a_M \tilde{b}(0, 0) + \dots) + a_M^2 x^{2M}( \tilde{c}(0) + \dots).
\end{eqnarray*}
By the uniqueness of the power serious expansion one must have $\tilde{N} = N$, otherwise the first nonzero term in the power serious expansion of $h(x, \phi(x))$ has an order of at least $N + 1$ (we assumed $\tilde{N} \geq N$ earlier). 

Case 2: $M = N$. 
If there is an interior line segment in the image of $\gamma$ with two endpoints $\gamma(s_1)$ and $\gamma(s_2)$ such that $\gamma(s_1) \in S_2$ and $\gamma(s_2) \in S_1$. We can parametrize $y(s)$ for $s \in [s_1, s_2]$ by 
\[y(s) = y(s_1) + T(x(s) - x(s_1)), \text{ where } T = \frac{d \alpha}{dx} (x(s_1)). \]
Since $\gamma(s_1) \in S_2$, one can use the second inequality in (\ref{z'y'S2}) to estimate
\[\big{|} \frac{d \alpha}{dx} (x(s_1))\big{|}  = \big{|} \frac{y'(s_1)}{x'(s_1)}\big{|}  \leq 2|y'(s_1)| \leq 2D x(s_1)^{\tilde{N}-1},\]
if $|x'(s)| \geq \frac{1}{2}$ by choosing $\epsilon$ small enough. With (\ref{zyS2}) one gets
\begin{eqnarray*}
|y(s)| &\leq& |y(s_1)| + |T| (|x(s)|+|x(s_1)|) \\
		& \leq & B x(s_1)^{\tilde{N}} + 2D x(s_1)^{\tilde{N}-1} (x(s)+x(s)) \\
		& \leq & (B + 4D) x(s)^{\tilde{N}},
\end{eqnarray*}
where the last inequality holds because $x(s)$ is increasing. For some $s \in (s_1, s_2)$, we have $y(s) = \phi(x(s))$ and so
\[|\phi(x(s))| \leq (B + 4D) x(s)^{\tilde{N}}.\] 
On the other hand, $|\phi(x(s))| \geq \frac{|a_N|}{2} x(s)^N$ if $s$ is sufficiently close to 0 by (\ref{0}). Thus the following relation holds:
\[\frac{|a_N|}{2} x(s)^N \leq (B + 4D) x(s)^{\tilde{N}} \Longrightarrow x(s)^{N - \tilde{N}} \leq \frac{2(B + 4D)}{|a_N|}.\]
Suppose $\tilde{N} > N$ then the left-hand side converges to infinity as $s$ approaches 0, a contradiction. So $\gamma$ eventually stops going from $S_2$ to $S_1$ which reduces to the case of one obstacle. 
Hence we proceed with $\tilde{N} = N$. 


7. Show $a(0, 0) > 0$, $\tilde{a}(0, 0) > 0$. 

Let's prove by contradiction. Suppose that $a(0,0) < 0$. 
Assume $\gamma$ has a switch point inside $S_1$ at $s = s_0$. That is to say $\gamma(s_0) \in S_1$ and for either $s > s_0$ nearby or $s < s_0$ nearly, $\gamma(s)$ is an interior line segment. Denote $(x(s_0), y(s_0))$ by $(x_0, y_0)$ for simplicity. Consider the intersection of the two-dimensional plane $y = y_0 + T(x - x_0)$ with the surface $z = g(x, y)$. Set
\begin{equation*}
f(x) = g(x, y_0 + T(x - x_0)), \text{ where } T = \frac{d \alpha}{dx} (x_0).
\end{equation*}
 It follows that
\begin{equation}\label{f''}
\frac{d^2f}{dx^2} (x_0) = g_{xx}(x_0,y_0) + 2g_{xy}(x_0, y_0)T + g_{yy}(x_0, y_0)T^2. 
\end{equation}
Using (\ref{1}), (\ref{zyall}) and choosing $\epsilon$ sufficiently small one can estimate
\begin{eqnarray}\label{g_xx(x_0, y_0)}
g_{xx}(x_0, y_0) &=& x_0^{N-2}[N(N-1)a(x_0, y_0) + x_0 p(x_0, y_0)] + y_0q(x_0, y_0) \\ \notag
			 &\leq& x_0^{N-2}N(N-1)\frac{1}{2}a(0, 0) + Bx_0^{N}C_1 \\ \notag
			 &\leq& x_0^{N-2}N(N-1)\frac{1}{2}a(0, 0)  - x_0^{N-2}N(N-1)\frac{1}{4}a(0, 0) \\ \notag
			 &=& x_0^{N-2}N(N-1)\frac{1}{4}a(0, 0) \notag
\end{eqnarray}
Furthermore, using the inequality in (\ref{T}) and letting $\epsilon$ be small enough one has
\begin{eqnarray}\label{other}
&& |2g_{xy}(x_0, y_0)T + g_{yy}(x_0, y_0)T^2|\\ \notag
&\leq& C_2 |T| \leq C_2 \cdot 2Dx_0^N \\
&\leq& -x_0^{N-2}N(N-1)\frac{1}{8}a(0, 0). \notag
\end{eqnarray}
So
\[\frac{d^2f}{dx^2} (x_0) \leq x_0^{N-2}N(N-1)\frac{1}{8}a(0, 0) < 0.\]
Therefore $f$ is concave downward at $x_0$ and the tangent line to the curve at $x_0$ is above the graph, a contradiction. 
In other words, $\gamma$ has no switch point on $S_1$ near the origin and so initially stays inside $S_2$ or is a line segment. Similarly $\gamma$ initially stays inside $S_1$ or is a line segment for $\tilde{a}(0, 0) < 0$. Hence this reduces to the case of one obstacle. Thus we proceed with assuming that $a(0, 0) > 0$ and $\tilde{a}(0, 0) > 0$.

8. Show that given $\epsilon$ small enough, if $\gamma$ leaves $S_1$ at a switch point, it will never enter $S_1$ again. Similarly, when $\gamma$ leaves $S_2$, it has to enter $S_1$ at the next switch point. 





Indeed we can prove by contradiction. Suppose $\gamma(s)$ leaves $S_1$ at $s=s_0$ and dives into the interior of $M$ for increasing $s$ until it enters $S_1$ again at $s = s_1$. Again set
\[f(x) = g(x, y_0+ T(x - x_0)),\] 
where $(x_0, y_0) = (x(s_0), y(s_0))$ and $T = \frac{d \alpha}{dx} (x_0)$. 
It follows that for $s\in[s_0, s_1]$ 
\[\frac{d^2f}{dx^2} (x(s)) = g_{xx}(x(s),y(s)) + 2g_{xy}(x(s), y(s))T + g_{yy}(x(s), y(s))T^2. \] Note that when $s = s_0$, this is just the expression in (\ref{f''}). Using an analogous argument as shown in (\ref{g_xx(x_0, y_0)}) for the case $a(0, 0) > 0$ one yields
\[g_{xx}(x(s), y(s)) \geq x(s)^{N-2}N(N-1)\frac{1}{4}a(0, 0).\] Moreover in analogy to (\ref{other}) one has
\begin{eqnarray*}
&&2g_{xy}(x(s), y(s))T + g_{yy}(x(s), y(s))T^2 \\
&\geq& -x_0^{N-2}N(N-1)\frac{1}{8}a(0, 0) \\
&\geq & -x(s)^{N-2}N(N-1)\frac{1}{8}a(0, 0).
\end{eqnarray*}
Hence
\[\frac{d^2f}{dx^2} (x(s)) \geq x(s)^{N-2}N(N-1)\frac{1}{8}a(0, 0) > 0 \text{ for all } s \in [s_0, s_1].\]
Therefore $f'(x(s))$ is increasing as $x(s)$ increases from $x(s_0) = x_0$ to $x(s_1) = x_1$. On the other hand, since the interior line segment is tangent to $S_1$ at the two endpoints, we must have
\[f'(x_0) = f'(x_1),\] a contradiction.
Therefore if $\gamma$ leaves $S_1$ at the switch point $\gamma(s_0)$ for some $s_0 \in [0, \epsilon]$, the geodesic arc beyond this point is a line segment never meeting $S_1$ again. Hence the next switch point (if there is one) lies on the surface $S_2$. The same argument holds for $S_2$ as well. 

9. The global behavior of $\gamma$.

\begin{lemma}\label{lem1}
Near the origin the geodesic is an alternating sequence of a boundary segment on $S_1$, an interval from $S_1$ to $S_2$, a boundary segment on $S_2$, an interval from $S_2$ to $S_1$, and so on.
\end{lemma}
\begin{proof}
Each time the projection of $\gamma$ crosses the graph of $\phi$ at time $s$, there is $l(s) > 0$ such that $\gamma$ is an interior line segment over the interval $[s - l(s), s + l(s)]$. The set $\mathcal{A}$ of such $s$ with $0 \leq s \leq \epsilon$ is therefore countable. Furthermore if $A_1 = \sup \mathcal{A}$ then  $A_1$ is actually the maximum of the set, because there is no $s \in \mathcal{A}$ within the $l(A_1)$-distance of $A_1$. Let $A_2 = \sup (\mathcal{A} - A_1)$, $A_3 = \sup (\mathcal{A} - \{A_1, A_2\})$, and so on. It follows that the set $\mathcal{A}$ can be linearly ordered as 
\[ \mathcal{A}  = \{A_1 > A_2 > A_3 > \dots\},
\]
such that between $A_n$ and $A_{n+1}$ the curve $\gamma$ lies entirely in $S_1$ or $S_2$ for each $n\geq 1$. 
\end{proof}

10. If $a_M > 0$, the curve $y=\phi(x)$ is concave upward for $x> 0$ nearby. We can obtain a contradiction as follows.
\begin{itemize}
\item{If $\gamma$ leaves a point in $S_2$ and enters a point in $S_1$, then $(x(s), y(s))$ crosses $\phi$ from above to below at some $s = s_1$. By concavity one must have $\frac{d \alpha}{dx} (x(s_1))  < \phi(x(s_1))$;}
\item{Later $(x(s), y(s))$ crosses $\phi$ from below to above at some $s = s_2$, then $\frac{d \alpha}{dx} (x(s_2)) > \phi(x(s_2))$;}
\item{In between $\gamma$ stays in $S_1$ all the time.}
\item{Since $\frac{d \alpha}{dx} (x(s_1)) = \frac{y'(s_1)}{x'(s_1)}$ and $\frac{d \alpha}{dx} (x(s_2)) = \frac{y'(s_2)}{x'(s_2)}$, we must have
\begin{eqnarray*}
&&y'(s_1) < x'(s_1)\phi(x(s_1)), y'(s_2) > x'(s_2) \phi(x(s_2)) \\ &\Rightarrow& y'(s_2) - y'(s_1) >  x'(s_2)\phi(x(s_2)) - x'(s_1) \phi(x(s_1))
\end{eqnarray*}}
\item{Therefore it suffices to show that for $0 < s_1 < s_2 < \epsilon$, 
\begin{equation*}
y'(s_2) - y'(s_1) \leq  x'(s_2)\phi(x(s_2)) - x'(s_1) \phi(x(s_1)).
\end{equation*}
}
\end{itemize}

On the one hand, 
\[y'(s_2) - y'(s_1) = \displaystyle\int_{s_1}^{s_2} y''(s) d s,\] 
where $y''(s) = - z''(s)(k + x(s)V_2(s)) \leq -z''(s) (1-\beta)k$ from (\ref{5}). Thus
\[y'(s_2) - y'(s_1) \leq \displaystyle\int_{s_1}^{s_2} -z''(s)(1-\beta)k d s < 0.\]
On the other hand, 
\begin{eqnarray*}
x'(s_2)\phi(x(s_2)) - x'(s_1) \phi(x(s_1)) &=& \displaystyle \int_{s_1}^{s_2} \frac{d}{ds} \big{[} x'(s)\phi'(x(s)) \big{]} \\
&=& \displaystyle \int_{s_1}^{s_2}  x''(s)\phi'(x(s)) + x'(s)^2 \phi''(x(s)) ds.
\end{eqnarray*}
Now let's estimate $x''(s)\phi'(x(s)) + x'(s)^2 \phi''(x(s))$. Since $\gamma(s) \in S_1$ for $s \in (s_1, s_2)$, one has $x''(s) = -z''(s)x(s)V_1(s)$ from (\ref{5}). Therefore $|x''(s)| \leq Ez''(s)x(s)$ for some positive constant $E$. By hypothesis $\gamma$ is parametrized by arc length, so $|x'(s)| \leq 1$.  With (\ref{zyS1}) and (\ref{z'y'S1}) one differentiates $z(s)=g(x(s), y(s))$ twice to obtain
\begin{eqnarray*}
z''(s) &=& g_{xx}(x(s), y(s))x'(s)^2 + 2g_{xy}(x(s), y(s))x'(s)y'(s) + g_{yy}(x(s), y(s))y'(s)^2\\
	&  & + g_x(x(s), y(s))x''(s) + g_y(x(s), y(s)) y''(s) \\
	&=&\left\{x(s)^{N-2}[N(N-1)a(x(s), y(s)) + x(s)p(x(s), y(s))]+ y(s)q(x(s), y(s))\right\}x'(s)^2 \\
	& & + y'(s)[2g_{xy}(x(s), y(s))x'(s) + g_{yy}(x(s), y(s))y'(s)] \\
	& & + [x(s)m(x(s), y(s)) + y(s)b(x(s), y(s))]x''(s)\\
	& &  + [k + x(s)k(x(s), y(s)) + y(s)l(y(s))]y''(s)\\
      &\leq & \left\{ x(s)^{N-2}C_1 + |y(s)|C_2\right\} \cdot 1 + |y'(s)|C_3 + C_4|x''(s)| + [k+\beta k]|y''(s)| \\
	&\leq& x(s)^{N-2}C_1 + Bx(s)^NC_2 + Dx(s)^{N-1}C_3 + C_4Ez''(s)x(s) + (1+\beta)k z''(s)(k+|x(s)V_2(x)|)\\
	&\leq&x(s)^{N-2}(C_1 + Bx(s)^2C_2 + Dx(s)C_3) + z''(s)(C_4Ex(s) + (1+\beta)k \cdot (1+\beta)k)\\
	&\leq&x(s)^{N-2}C_5 + z''(s)C_6 + (1+\beta)^2k^2 z''(s).
\end{eqnarray*}
By hypothesis $(1+\beta)k < 1$, then we can choose $\epsilon$ small enough so that $C_6 < 1-(1+\beta)^2k^2$ implying that
\[z''(s) \leq Fx(s)^{N-2}, \text{ and so } |x''(s)| \leq EF x(s)^{N-1} = Gx(s)^{N-1}.\]
Now let's use (\ref{0}) to approximate $\phi''(x(s))$ and $\phi'(x(s))$:
\[0< \phi'(x) = Ma_Mx^{M-1} + (M+1)a_{M+1}x^M + \dots \leq 2Ma_Mx^{M-1} \text{ for $x$ near 0.}\] So for $\epsilon$ sufficiently small, one has
\begin{eqnarray}\label{9}
&& \phi'(x(s)) x''(s) \geq \phi'(x) \cdot -G x(s)^{N-1} \\ \notag
&\geq& 2Ma_M x(s)^{M-1} \cdot -G x(s)^{N-1} = -2MGa_Mx(s)^{M+N-2}. \notag
\end{eqnarray}
On the other hand, 
\[\phi''(x) = M(M-1)a_Mx^{M-2} + \dots \geq \frac{1}{2}a_MM(M-1)x^{M-2} \text{ for $x$ near 0.} \] So again by choosing $\epsilon$ small enough and assuming $x'(s) \geq \frac{1}{2}$, one obtains
\begin{equation}\label{10}
\phi''(x(s))x'(s)^2 \geq \frac{1}{8}a_MM(M-1)x(s)^{M-2}.
\end{equation}
Combing (\ref{9}) and (\ref{10}) we find that 
\[\phi''(x(s))x'(s)^2 +\phi'(x(s)) x''(s) \geq x(s)^{M-2}a_M(\frac{1}{8}M(M-1) -2MGx(s)^N).\]
Since $N \geq 2$ the above difference can be made positive for every $s \in [0, \epsilon]$ if $\epsilon$ is sufficiently small. Hence 
\[x'(s_2)\phi(x(s_2)) - x'(s_1) \phi(x(s_1)) =  \displaystyle \int_{s_1}^{s_2}  x''(s)\phi'(x(s)) + x'(s)^2 \phi''(x(s)) ds > 0.\] We reach a contradiction. There $\gamma(s)$ eventually stops bouncing between $S_1$ and $S_2$ as $s$ approaches 0, which reduces to the case of one obstacle. 

11. If $a_M < 0$, the curve $y = \phi(x)$ is concave downward for $x > 0$ nearby. Replacing $S_1$ by $S_2$ and $g$ by $h$ in the previous argument gives a contradiction.

12. Trivial Case 1: $\phi(x)$ is identically zero, but $g(x, 0)$ and $h(x, 0)$ are not identically zero. 

Everything is fine until Step 9 by letting $M = \infty$. In Step 10, the proof is as follows:
\begin{itemize}
\item{If $\gamma$ leaves a point in $S_2$ and enters a point in $S_1$, then $(x(s), y(s))$ crosses $\phi$ from above to below at some $s = s_1$. Since the curve $y=\phi(x)$ is the $x$-axis, one must have $\frac{d \alpha}{dx} (x(s_1))  < 0$ and so $y'(s_1) < 0$. }
\item{Later $(x(s), y(s))$ crosses $\phi$ from below to above at some $s = s_2$, thus  $\frac{d \alpha}{dx} (x(s_2)) > 0$ and so $y'(s_2) > 0$. }
\item{In between $\gamma$ stays in $S_1$ all the time where $y''(s)$ is always negative, so $y'(s_2) < y'(s_1)$, a contradiction.}
\end{itemize}

13. Trivial Case 2: one of $g(x, 0)$, $h(x, 0)$ is identically zero, but not both. 

Without loss of generality, let us assume that $g(x, 0)$ is identically zero but $h(x, 0)$ is not. Then we can write $g(x, y)$ as 
\begin{eqnarray*}
g(x, y) = ky + xyb(x, y) + y^2 c(y).
\end{eqnarray*}
If $\phi(x)$ is identically zero, then $g(x, \phi(x))$ is also identically zero. However, with equality (\ref{2})
\[h(x, \phi(x) )= x^{\tilde{N}} \tilde{a}(x, 0) = x^{\tilde{N}} ( \tilde{a}(0, 0) + \dots) \neq 0. \]
So $\phi(x)$ is nonzero whose power serious expansion is still (\ref{0}). 
Everything is fine until step 4. If $\gamma(s) \in S_1$, the equality in (\ref{1}) gives 
\begin{eqnarray*}
z(s) &=& g(x(s), y(s)) \\
&=& ky(s) + x(s)y(s)b(x(s), y(s)) + y(s)^2c(y(s)) \\
&\leq& k|y(s)| + |y(s)||x(s)b(x(s), y(s)) + y(s)c(y(s))| \\
&\leq& k|y(s)| + C_1|y(s)| \\
&\leq& (k + C_1) (1+ \beta)k_1 z(s)
\end{eqnarray*}
By hypothesis $(1+\beta)k < 1$, for $s$ close enough to 0 such that $C_1 < \frac{1 - (1+\beta)k^2}{(1+ \beta)k}$, one obtains $z(s) \leq 0$. On the other hand, $z(s) > 0$ for all $s > 0$. So the geodesic does not touch $S_1$ near the origin, which reduces to the case of one obstacle. 

14. $g(x, 0)$ and $h(x, 0)$ are both identically 0.  One can show that $\phi(x)$ is also identically zero. In this case the geodesic does not touch $S_1$ and $S_2$ near the origin, so it must be a line segment at the beginning. 
\end{proof}

%% file: part4.tex
\begin{theorem}\label{thm2}
Let $M_1$ and $M_2$ be 3-dimensional analytic manifolds with boundary embedded in $\mathbb{R}^3$. Denote the boundary surfaces of $M_1$ and $M_2$ by $S_1$ and $S_2$, respectively. Assume that $S_1$ and $S_2$ intersect transversally whose angle is greater than or equal to $90^\circ$. Let $M$ be the intersection of $M_1$ and $M_2$. Let $\gamma$ be a geodesics  in $M$ parametrized by arc length $s$, with $\gamma(0) = p \in S_1 \cap S_2$. Then there is an $\epsilon > 0$ such that $\gamma$ has no switch point for $0 < s < \epsilon$. 
\end{theorem}

\begin{proof}

1. Set the coordinate system as in Theorem \ref{thm1}, except that 
the outward normal vectors to $S_1$ and $S_2$ at $p$ are no longer symmetrical with respect to the $z$-axis. Tilt the system appropriately so that the normal vectors are $(0, -k_1, 1)$ and $(0, k_2, 1)$, respectively, where $0 < k_1 < 1 < k_2 $. Furthermore, choose $k_1$ and $k_2$ so $k_1k_2 < \frac{1}{1+\beta} < 1$ for some $\beta > 0$. It follows that
\[g_x (0, 0) = 0, \ \ \  g_y(0, 0) = k_1.\] 
And 
\[h_x(0, 0) = 0, \ \ \ h_y(0, 0) = -k_2.\] 
Since
\[ \frac{\partial}{\partial y} \big{|}_{(0, 0)}  [g(x, y) - h(x, y)] = k_1+k_2 > 0,\]
 again the implicit function theorem implies that the intersection of $S_1$ and $S_2$ near $p$ is a real analytic curve defined by the following equations:
 \[ y = \phi(x), \ \ \ z = g(x, \phi(x)) = h(x, \phi(x)),\]
where $\phi$ has the same power expansion as in (\ref{0}) by assuming that $\phi$ is not identically zero. Furthermore let us assume that $g(x, 0)$ and $h(x, 0)$ are not identically zero. Therefore the equation defining $S_1$ near $p$ is of the form
\begin{equation}\label{part4g}
g(x, y) = k_1y + x^N a(x, y) + xy b(x, y) + y^2 c(y), 
\end{equation}
where $N \geq 2$, the functions $a, b, c$ are analytic, and $a(0, 0) \neq 0$.  Moreover the equation defining $S_2$ near $p$ is of the form
\begin{equation}\label{part4h}
h(x, y) = -k_2y + x^{\tilde{N}} \tilde{a} (x, y) + xy \tilde{b}(x, y) + y^2 \tilde{c}(y),
\end{equation}
where $\tilde{N} \geq 2$, the functions $\tilde{a}, \tilde{b}, \tilde{c}$ are analytic, and $\tilde{a}(0,0) \neq 0$.

2. Same as in Theorem \ref{thm1}, the graph of $\phi(x)$ divides the $(x, y)$-plane into two parts near 0: the part below the graph corresponding to the projection of the surface $S_1$ in $M$ and the part above the graph corresponding to the projection of the surface $S_2$ in $M$. 

3. Concavity of $\phi$ is determined by the sign of $a_M$ and there are two cases to consider: $a_M > 0$ and $a_M < 0$.

4. Approximate $y(s)$ and $y'(s)$ using the normal vectors $N_1(s)$, $N_2(s)$ to $S_1$, $S_2$. 

Following the same procedure as in (\ref{5}) and (\ref{6}) one obtains
if $\gamma(s) \in S_1$, 
\[x''(s) = - z''(s)x(s)V_1(s), \ \ \ y''(s) = -z''(s) (k + x(s) V_2(s)),\]
where $|x(s) V_2(s)| \leq \beta k_1$ for $\epsilon$ sufficiently small and therefore $|y''(s)| \leq (1+\beta)k_1 |z''(s)|$; and if $\gamma(s) \in S_2$, 
\[x''(s) = - z''(s)x(s)W_1(s), \ \ \ y''(s) = -z''(s) (-k + x(s) W_2(s)),\]
where $|x(s) W_2(s)| \leq \beta k_2$ for $\epsilon$ small enough and hence $|y''(s)| \leq (1+\beta)k_2 |z''(s)|$.

Case 1, if $\gamma(s) \in S_1$, since $k_1(1+ \beta) < 1$ one can use the same argument as before to obtain
\begin{eqnarray}\label{part4zyS1}
z(s) \leq Ax(s)^N \Rightarrow |y(s)|  \leq (1+\beta)k_1Ax(s)^N = B x(s)^N.
\end{eqnarray} 

Next differentiating $z(s) = g(x(s), y(s))$ once gives
\begin{eqnarray}\label{part4z'y'S1}
z'(s) \leq Cx(s)^{N-1} \Rightarrow |y'(s)| \leq (1+\beta)k_1C x(s)^{N-1} = D x(s)^{N-1}.
\end{eqnarray}

Case 2: if $\gamma(s) \in S_2$, the arguments above fail because $k_2 > 1$. Instead we need to use a different approach. 
Since $y''(s) = -z''(s)(-k_2 + x(s)W_2(s))$, 
\begin{equation}\label{z''1}
z''(s) + k_2 y''(s) = z''(s)[1+ k_2^2 - k_2x(s)W_2(s)] \geq z''(s)[1+(1-\beta)k_2^2].
\end{equation}
On the other hand, if $\gamma(s) \in S_1$, then $y''(s) = -z''(s)(k_1 + x(s)V_2(s))$ implies that
\begin{eqnarray}\label{z''2}
z''(s) + k_2y''(s) &=& z''(s)[1 - k_1k_2 - k_2x(s)V_2(s)]  \\ \notag
&\geq&  z''(s)[1-k_1k_2 - k_2k_1\beta] \\ \notag
&=& z''(s)[1-(1+\beta)k_1k_2], \notag
\end{eqnarray}
Combining (\ref{z''1}) with (\ref{z''2}), together with the hypothesis $(1+\beta)k_1k_2$, there exists a constant $H$ such that if $\gamma(s) \in S_1$ or $S_2$
\begin{equation}\label{part4z''}
z''(s) + k_2y''(s) \geq H z''(s).
\end{equation}
Note that  the inequality (\ref{part4z''}) still holds if $\gamma(s)$ does not touch any surface since the acceleration is 0. It follows that for every $s$ in the interval $[0, \epsilon]$, 
\begin{equation}\label{part4z'}
z'(s) + k_2y'(s) = \displaystyle \int_0^s z''(\sigma) + k_2 y''(\sigma) d \sigma \geq \displaystyle \int_0^s H z''(\sigma) d \sigma = H z'(s).
\end{equation}
Moreover, 
\begin{equation}\label{part4z}
z(s) + k_2y(s) \geq H z(s) \text{ for every } s \in [0, \epsilon]. 
\end{equation}
With (\ref{part4h}) and (\ref{part4z}) one deduces that 
\begin{eqnarray*}
Hz(s) &\leq& x(s)^{\tilde{N}} \tilde{a}(x(s), y(s)) + y(s)[x(s) \tilde{b}(x(s), y(s)) + y(s) \tilde{c}(y(s))] \\
&\leq& C_1x(s)^{\tilde{N}} + |y(s)| C_2 \\
& \leq& C_1 x(s)^{\tilde{N}} + (1+\beta)k_2z(s) C_2.
\end{eqnarray*}
With $\epsilon$ sufficiently small $C_2$ can be made as small as possible so that 
\[H- (1+\beta)k_2C_2 > 0, \]
implying
\begin{eqnarray}\label{part4zyS2}
z(s) \leq A x(s)^{\tilde{N}} \Rightarrow |y(s)| \leq (1+\beta)k_2A x(s)^{\tilde{N}} = B x(s)^{\tilde{N}},
\end{eqnarray}
after enlarging $A$ and $B$ accordingly.

Differentiating $z(s) = h(x(s), y(s))$ once, together with (\ref{part4z'}) and (\ref{part4zyS2}), one obtains
\begin{eqnarray*}
z'(s) &=&  -k_2y'(s)+ x(s)^{\tilde{N}-1}[\tilde{N}x'(s)\tilde{a}(x(s), y(s)) \\
&&+ x(s)(\tilde{a}_x(x(s), y(s))x'(s) + \tilde{a}_y(x(s), y(s))y'(s))] \\
	&  & + y'(s)[x(s)\tilde{b}(x(s), y(s)) + x(s)y(s) \tilde{b}_y(x(s), y(s)) \\
	&&+ 2y(s)\tilde{c}(y(s)) + y^2(s)\tilde{c}'(y(s))]\\
	&   & + y(s)[x'(s)\tilde{b}(x(s), y(s)) + x(s)x'(s) \tilde{b}_x(x(s), y(s))] \\
\text{so } z'(s) + k_2y'(s)	& \leq& C_1x(s)^{N-1} + C_2|y'(s)| + C_3|y(s)| \\
	& \leq& C_1x(s)^{\tilde{N}-1}  + C_2(1+\beta)k_2 z'(s) + C_3 Bx(s)^{\tilde{N}} \\
\text{so }Hz'(s) &\leq&  C_1x(s)^{\tilde{N}-1}  + C_2(1+\beta)k_2 z'(s) + C_3 Bx(s)^{\tilde{N}}
\end{eqnarray*}
With $\epsilon$ sufficiently small $C_2$ can be made as small as possible so that $C_2(1+\beta) < H$ and thus
\begin{eqnarray}\label{part4z'y'S2}
z'(s) \leq C x(s)^{\tilde{N}-1} \Rightarrow |y'(s)| \leq (1+\beta)k_2 z'(s) = D x(s)^{\tilde{N}-1},
\end{eqnarray}
by enlarging $C$ and $D$ accordingly.

Without loss of generality we may assume that $N \leq \tilde{N}$, then $x(s)^{\tilde{N}} \leq x(s)^N$ for $x(s) \leq 1$. Thus (\ref{part4zyS1}) and (\ref{part4zyS2}) imply that  if $\gamma \in S_1$ or $S_2$ 
\begin{equation}\label{part4zyS1S2}
z(s) \leq A x(s)^N, \ \ \ |y(s)| \leq B x(s)^N,
\end{equation}
Furthermore (\ref{part4z'y'S1}) and (\ref{part4z'y'S2}) imply that if $\gamma \in S_1$ or $S_2$
\begin{equation}\label{part4z'y'S1S2}
z'(s) \leq C x(s)^{N-1}, \ \ \ |y'(s)| \leq D x(s)^{N-1}.
\end{equation}
When $\gamma(s)$ does not lie on $S_1$ or $S_2$, following the same argument as in Theorem \ref{thm1}, one can show that
\begin{equation}\label{part4T}
|T(x(s))| \leq 2D x(s)^{N-1}.
\end{equation}
\begin{equation}\label{part4yall}
|y(s)| \leq Bx(s)^N \text{ for every } s \in [0, \epsilon].
\end{equation}

From 5-10, one can copy the proof from Theorem \ref{thm1} word by word to have: $M \geq N$, $N = \tilde{N}$, $a(0, 0)> 0$, $\tilde{a}(0, 0) > 0$, and $\gamma$ is an alternating sequence of  boundaries segments on $S_1$, $S_2$ and line segments between $S_1$ and $S_2$. 

11. If $a_M < 0$, the argument is slightly different.  
\begin{itemize}
\item{If $\gamma$ leaves a point in $S_1$ and enters a point in $S_2$, then $(x(s), y(s))$ crosses $\phi$ from below to above at some $s = s_1$. Since the curve $y = \phi(x)$ is concave downward for $x > 0$, one must have $\frac{d \alpha}{dx} (x(s_1))  > \phi(x(s_1))$;}
\item{Later $(x(s), y(s))$ crosses $\phi$ from above to below at some $s = s_2$, then $\frac{d \alpha}{dx} (x(s_2)) < \phi(x(s_2))$;}
\item{In between $\gamma$ stays in $S_2$ all the time.}
\item{Since $\frac{d \alpha}{dx} (x(s_1)) = \frac{y'(s_1)}{x'(s_1)}$ and $\frac{d \alpha}{dx} (x(s_2)) = \frac{y'(s_2)}{x'(s_2)}$, we must have
\begin{eqnarray*}
y'(s_1) > x'(s_1)\phi(x(s_1)), y'(s_2) < x'(s_2) \phi(x(s_2)) \\
 \Rightarrow y'(s_2) - y'(s_1) <  x'(s_2)\phi(x(s_2)) - x'(s_1) \phi(x(s_1)).
\end{eqnarray*}
}
\item{Therefore it suffices to show that for $0 < s_1 < s_2 < \epsilon$, 
\begin{equation*}
y'(s_2) - y'(s_1) \geq  x'(s_2)\phi(x(s_2)) - x'(s_1) \phi(x(s_1)).
\end{equation*}}
\end{itemize}

On the one hand, 
\[y'(s_2) - y'(s_1) = \displaystyle\int_{s_1}^{s_2} y''(s) d s,\] 
where $y''(s) = - z''(s)(-k_2 + x(s)W_2(x)) \geq z''(s)(1-\beta)k$. Thus
\[y'(s_2) - y'(s_1) \geq \displaystyle\int_{s_1}^{s_2} z''(s)(1-\beta)kd s > 0.\]
On the other hand, 
\begin{eqnarray*}
x'(s_2)\phi(x(s_2)) - x'(s_1) \phi(x(s_1)) &=& \displaystyle \int_{s_1}^{s_2} \frac{d}{ds} \big{[} x'(s)\phi'(x(s)) \big{]} \\
&=& \displaystyle \int_{s_1}^{s_2}  x''(s)\phi'(x(s)) + x'(s)^2 \phi''(x(s)) ds.
\end{eqnarray*}
Now let's estimate $x''(s)\phi'(x(s)) + x'(s)^2 \phi''(x(s))$.  Since $\gamma(s) \in S_2$ for $s \in (s_1, s_2)$, one has $x''(s) = -z''(s)x(s)W_1(s)$ from (\ref{6}). Therefore $|x''(s)| \leq Ez''(s)x(s)$ for some positive constant $E$. By hypothesis $\gamma$ is parametrized by arc length, so $|x'(s)| \leq 1$.  One differentiates $z(s)=h(x(s), y(s))$ twice to get
\begin{eqnarray*}
z''(s) &=& h_{xx}(x(s), y(s))x'(s)^2 + 2h_{xy}(x(s), y(s))x'(s)y'(s) + h_{yy}(x(s), y(s))y'(s)^2\\
	&  & + h_x(x(s), y(s))x''(s) + h_y(x(s), y(s)) y''(s) \\
	&=&\big{\{}x(s)^{N-2}[N(N-1)\tilde{a}(x(s), y(s)) + x(s)\tilde{p}(x(s), y(s))] \\
& & + y(s)\tilde{q}(x(s), y(s))\big{\}}x'(s)^2 \\
	& & + y'(s)[2h_{xy}(x(s), y(s))x'(s) + h_{yy}(x(s), y(s))y'(s)] \\
	& & + [x(s)W_1(s)]x''(s) + [-k_2 + x(s)W_2(s)]y''(s)
\end{eqnarray*}
One moves $-k_2y''(s)$ to the other side of the inequality, together with (\ref{part4zyS2}) and (\ref{part4z'y'S2}),  to obtain
\begin{eqnarray*}
z''(s) + k_2 y''(s)	 &\leq & \left\{ x(s)^{N-2}C_1 + |y(s)|C_2\right\} \cdot 1 + |y'(s)|C_3 + C_4|x''(s)| + |x(s)W_2(s)||y''(s)| \\
	&\leq& x(s)^{N-2}C_1 + Bx(s)^NC_2 + Dx(s)^{N-1}C_3 + C_4Ez''(s)x(s)\\
	&& + |x(s)W_2(s) |z''(s)(k_2+|x(s)W_2(x)|)\\
	&\leq&x(s)^{N-2}(C_1 + Bx(s)^2C_2 + Dx(s)C_3) \\
&& + z''(s)(C_4Ex(s) + \beta k_2(1+\beta) k_2)\\
	&\leq&x(s)^{N-2}C_5 + \beta(1+\beta)k_2^2z''(s)C_6.
\end{eqnarray*}
Using (\ref{part4z''}) and choosing $\epsilon$ small enough so that $\beta(1+\beta)k_2^2C_6 < H$ one gets that
\[z''(s) \leq Fx(s)^{N-2}, \text{ and so } |x''(s)| \leq EF x(s)^{N-1} = Gx(s)^{N-1}.\]
Now let's use (\ref{0}) to approximate $\phi''(x(s))$ and $\phi'(x(s))$: 
\[0> \phi'(x) = Ma_Mx^{M-1} + (M+1)a_{M+1}x^M + \dots \geq 2Ma_Mx^{M-1} \text{ for $x$ near 0.}\] So for $\epsilon$ sufficiently small, one has
\begin{eqnarray}\label{11}
&&\phi'(x(s)) x''(s) \leq \phi'(x(s)) \cdot -Gx(s)^{N-1} \\ \notag
 &\leq&  2Ma_Mx(s)^{M-1} \cdot -G x(s)^{N-1} = -2Ma_MLx(s)^{M+N-2}. \notag
\end{eqnarray}
On the other hand, 
\[\phi''(x) = M(M-1)a_Mx^{M-2} + \dots \leq \frac{1}{2}a_MM(M-1)x^{M-2} \text{ for $x$ near 0.} \] So using $\epsilon$ small enough and assuming $x'(s) \geq \frac{1}{2}$, one obtains
\begin{equation}\label{12}
\phi''(x(s))x'(s)^2 \leq \frac{1}{8}a_MM(M-1)x(s)^{M-2}.
\end{equation}
Combining (\ref{11}) and (\ref{12}) we find that 
\[\phi''(x(s))x'(s)^2 +\phi'(x(s)) x''(s) \leq x(s)^{M-2}a_M(\frac{1}{8}M(M-1) -2MGx(s)^N).\]
Since $N \geq 2$ the above difference can be made negative for every $s \in [0, \epsilon]$ if $\epsilon$ is sufficiently small. Hence 
\[x'(s_2)\phi(x(s_2)) - x'(s_1) \phi(x(s_1)) =  \displaystyle \int_{s_1}^{s_2}  x''(s)\phi'(x(s)) + x'(s)^2 \phi''(x(s)) ds < 0.\] We reach a contradiction. Thus $\gamma$ eventually stops bouncing between $S_1$ and $S_2$ as $s$ approaches 0, which reduces to the case of one obstacle.

The trivial cases from 11-14 follow exactly the same proof in Theorem \ref{thm1}.

\end{proof}

\end{section}

%% file: part7.tex
\begin{section}{Part Two}


\begin{theorem}
Let $M$ be an 3-dimensional manifold with boundary embedded in $\mathbb{R}^{3}$. Denote the boundary surface of $M$ by $S$ and let $\gamma(s)$ be a geodesic on $M$ parametrized by arc length $s$ with $\gamma(0) = p \in S$. Then there exists an $\epsilon > 0$ such that the number of line segments in the image of $\gamma$ within the $\epsilon$-ball of $p$ is uniformly bounded, namely it is independent of the initial velocity $\gamma'(0)$. More precisely, there are at most two complete or partial line segments. 
\end{theorem}

The idea is to show that in each direction there exists a wedge and an $\epsilon$ such that if $\gamma'(0)$ is within this wedge, then $\gamma$ has a uniform bound on the number of switch points within the $\epsilon$-ball of $p$.

1. Set up the coordinate system. 

Choose an orientation of the coordinate system $(x, y, z)$ so that $p$ is the origin, $S$ near $p$ can be parametrized by an analytic function $z = g(x, y)$, and the outward normal vector to $S$ at $p$ is in the positive $z$-axis. Notice that $\gamma'(0)$ is either tangent to $S$ at $p$ in which case $\gamma'(0)$ is in the $(x, y)$-plane. (Or $\gamma'(0)$ has a negative $z$-component pointing towards the interior of $M$. We will look at this case in the end.)

2. Without loss of generalizty, one may assume that the lowest degree in the power series expansion of $g(x, y)$ is 2. The idea for higher degrees is very similar, which will be mentioned at the end. 

Since the $z$-axis is normal to $S$ at $p$ and $p$ is the origin, the Taylor series expansion of $g$ is
\begin{eqnarray*}
g(x, y) = \frac{1}{2}g_{xx}(0, 0) x^2 + g_{xy}(0, 0)xy + \frac{1}{2}g_{yy}(0, 0)y^2 + \text{ higher-order terms, }
\end{eqnarray*}
where 
\[ \frac{1}{2}g_{xx}(0, 0) x^2 + g_{xy}(0, 0)xy + \frac{1}{2}g_{yy}(0, 0)y^2 \] is not the zero polynomial. Let $H$ be the Hessian matrix
\begin{equation*}
H = \left[
\begin{matrix}
g_{xx}(0, 0) & g_{xy}(0, 0) \\
g_{xy}(0, 0) & g_{yy}(0, 0)
\end{matrix} \right].
\end{equation*}
Since $H$ is symmetric, the spectral theorem says that there exists a $2 \times 2$ real orthogonal matrix $P$ such that 
\[PHP^{t} = \left[\begin{matrix} 2a & 0 \\ 0 & 2b \end{matrix} \right].\]
Rotating and/or reflecting the $(x, y)$-plane using $P$, the surface $S$ near $p$ can be parametrized as follows:
\[g(x, y) = ax^2 + by^2 + \text{ higher-order terms, }\]
where $a$ and $b$ are not identically zero. The signs of $a$, $b$ tell us about the shape of $S$ near the origin, namely
\begin{itemize}
\item{When $a > 0$, $b \geq 0$ or $a \geq 0$, $b > 0$, $S$ is concave upward near the origin and $\gamma$ has at most one switch point near $p$.}
\item{When $a < 0$, $b \leq 0$ or $a \leq 0$, $b < 0$, $S$ is concave downward near the origin and $\gamma$ has no switch point near $p$.}
\item{When one of $a, b$ is positive and the other is negarive, $S$ has a saddle point at $p$ and we are going to investigate this case in details. Without loss of generality, we may assume that $a > 0$, $b <0$. Moreover replacing $b$ by $-b$ yields
\begin{equation*}
g(x, y) = ax^2 - by^2 + \text{ higher-order terms, where } a> 0, b> 0.
\end{equation*}
Notice that when $ax^2 - by^2 = 0$, $y = \pm \displaystyle \sqrt{\frac{a}{b}} x,$ which gives rise to two lines dividing the plane into four different regions. There exists an angle $0 < \theta_0 < \pi/2$ such that
\[\tan(\theta_0) = \displaystyle \sqrt{\frac{a}{b}}.\]
Then we are going to prove the following:
\begin{enumerate}
\item{In the positive direction of the $x$-axis: for any small positive $\delta$, there exists an $\epsilon > 0$ such that if $\gamma'(0)$ lies inside the wedge $[-\theta_0 + \delta, \theta_0 - \delta]$, then $\gamma$ has at most one switch point within the $\epsilon$-ball of $p$. Similar statements hold in the negative direction of the $x$-axis and positive and negative directions of the $y$-axis.}
\item{In the positive direction of $y = \displaystyle \sqrt{\frac{a}{b}} x:$
there exist  an $\eta > 0$ and an $\epsilon > 0$ such that if $\gamma'(0)$ is in the wedge $[\theta_0 - \eta, \theta_0+ \eta]$, then $\gamma$ has at most two switch points within the $\epsilon$-ball of $p$. Similar statements also hold in the negative direction of  $y = \displaystyle \sqrt{\frac{a}{b}} x$ and the positive and negative directions of $y = - \displaystyle \sqrt{\frac{a}{b}} x$.}
\item{Combining 1 and 2, the theorem follows.}
\end{enumerate}
}
\end{itemize}

\begin{proposition}\label{proposition1}
In the positive direction of the $x$-axis: for any $0 < \delta < \theta_0$, there exists an $\epsilon > 0$ such that if $\gamma'(0)$ lies inside the wedge $[-\theta_0 + \delta, \theta_0 - \delta]$, then $\gamma$ has at most one switch point within the $\epsilon$-ball of $p$. 
\end{proposition}

\begin{proof}
(1) Set up the frame. The first estimate of $\epsilon$ comes from  that $S$ is parametrized by $g(x, y)$ within the $\epsilon$-ball of $p$. 
Since $\gamma'(0)$ is a unit vector tangent to $S$, one has $x'(0) = \cos\theta$ and $y'(0) = \sin \theta$, where $\theta \in [-\theta_0 + \delta, \theta_0 - \delta]$ by hypothesis. 

(2)  Show that $x'(s) > 0$ if $\epsilon$ is chosen small enough. This is the second estimate of $\epsilon$. 
Here we are assuming that if $\gamma(s)$ is within the $\epsilon$-ball of $p$, then for every $0 \leq \sigma \leq s$, $\gamma(\sigma)$ also lies within the $\epsilon$-neighbhorhood of $p$. 

Let $\gamma(s) = (x(s), y(s), z(s))$ with $|x(s)|$,$|y(s)|$, $|z(s)|$ less than or equal to $\epsilon$, so that $\gamma(s)$ is within the $\epsilon$-ball of $p$. If $\gamma(s) \in S$, then the normal vector at $\gamma(s)$ is 
\[N(s) = (-g_x(x(s), y(s)), -g_y(x(s), y(s)), 1).\]
Let 
\begin{equation}\label{gabcdef}
g(x, y) = ax^2 - by^2 + x^3c(x, y) + x^2yd(x, y) + xy^2e(x, y) + y^3f(x, y),
\end{equation}
then 
\begin{eqnarray*}
g_x(x, y) &=& 2ax + (3x^2c + x^3c_x) + (2xyd+ x^2yd_x) + (y^2e + xy^2e_x) + y^3f_x.\\
g_y(x, y) &=& -2by + x^3c_y + (x^2d + x^2yd_y) + (2xye + xy^2e_y) + (3y^2 f + y^3f_y).
\end{eqnarray*}
There exists a positive constant $A$ such that
\begin{equation}\label{gxgy}
g_x(x, y) \leq A \text{ and } g_y(x, y) \leq A, \text{ if } |x|, |y| \leq \epsilon,
\end{equation}
where $A \to 0$ as $\epsilon \to 0$. 
Let $s$ be such that $\gamma''(s)$ exists, then $\gamma''(s) = z''(s)N(s)$. This implies that 
\[x''(s)= -z''(s)g_x(x(s), y(s)) \Rightarrow |x''(s)| \leq Az''(s). \]
Here $z''(s) \geq 0$, because within the $\epsilon$-ball of $p$ the outward normal vector to $S$ has a positive $z$-coordinate of 1 and $\gamma''(s)$ directs outward on a boundary segment on $S$. Indeed $\gamma(s)$ is a locally shortest path. If $\gamma(s)$ lies on the surface of $M$, its acceleration exists everywhere except at the switch points and it outward normal to the surface. On the other hand, if $\gamma(s)$ lies on a line segment in the interior of $M$, then the acceleration is zero. So we obtain
\begin{equation}\label{x'positive}
x'(s) = x'(0) + \displaystyle \int_0^s x''(\sigma) d\sigma \geq \cos \theta - A \displaystyle \int_0^s z''(\sigma) d\sigma = \cos \theta - Az'(s). 
\end{equation}

Next let's approximate $z'(s)$. If $\gamma(s) \in S$, then 
\[z(s) = g(x(s), y(s)) = ax(s)^2 - by(s)^2 + x(s)^3c + x(s)^2y(s)d + x(s)y(s)^2e + y(s)^3f,\] so
\begin{eqnarray*}
z'(s) &=& 2ax(s)x'(s) - 2by(s)y'(s) + 3x(s)^2x'(s) c(x(s), y(s)) + x(s)^3 c_x(x(s), y(s))x'(s) \\
&&+  x(s)^3 c_y(x(s), y(s))y'(s) + 2x(s)x'(s)y(s)d(x(s), y(s)) + x(s)^2y'(s)d(x(s), y(s)) \\
&&+ x(s)^2y(s)d_x(x(s), y(s))x'(s) + x(s)^2y(s)d_y(x(s), y(s))y'(s) x'(s)y(s)^2 e(x(s), y(s)) \\
&& + 2x(s)y(s)y'(s)e(x(s), y(s)) + x(s)y(s)^2e_x(x(s), y(s))x'(s) + x(s)y(s)^2e_y(x(s), y(s))y'(s) \\
&& 3y(s)^2y'(s)f(x(s), y(s)) + y(s)^3f_x(x(s), y(s))x'(s) +  y(s)^3f_y(x(s), y(s))y'(s).
\end{eqnarray*}
Since $\gamma(s)$ is parametrized by arc length, $|x'(s)|$ and $|y'(s)|$ are no more than 1. Since each term in the above expression has either an $x(s)$ or $y(s)$, there exists a positive constant $B$ such that 
\[|z'(s)| \leq B, \text{ if } |x(s)|, |y(s)| \leq \epsilon,\]
where $B \to 0$ as $\epsilon \to 0$. On the other hand, if $\gamma(s)$ is within an interior line segment, then $\gamma'(s)$ is constant and equal to the value at the endpoints. Therefore $|z'(s)|$ is still bounded by $B$. Thus we can choose $\epsilon$ small enough so that $B < \frac{\cos |\theta|}{A}$, where $|\theta| \leq \theta_0 -\delta$. It follows from (\ref{x'positive}) that 
\[x'(s) > 0,\]
if $|x(s)|, |y(s)| \leq \epsilon$ for $\epsilon$ sufficiently small. Therefore $x(s) > 0$. 

(3) Approximate $y'(s)$ and $y(s)$. This is the third estimate of $\epsilon$. 

If $\gamma(s) \in S$, then the normal vector to $S$ at $\gamma(s)$ is 
\begin{eqnarray*}
N(s) = (-g_x(x(s), y(s)), -g_{y}(x(s), y(s)) ,1).
\end{eqnarray*}
Let $s$ be such that $\gamma''(s)$ exists, then $\gamma''(s) = z''(s)N(s)$. From (\ref{gxgy}) one obtains
\begin{eqnarray*}
y''(s) = -z''(s)g_y(x(s), y(s)) \Rightarrow |y''(s)| \leq A z''(s), 
\end{eqnarray*}
where $A \to 0$ as $\epsilon \to 0$. Thus
\begin{eqnarray*}
y'(s) = y'(0) + \displaystyle \int_0^s y''(\sigma) d\sigma \leq \sin \theta + A\displaystyle \int_0^s z''(\sigma) d\sigma = \sin \theta + Az'(s).
\end{eqnarray*}
Let $\tan |\theta| < c < \tan \theta_0$, where $|\theta| \leq \theta_0 - \delta$. 
Choose $\epsilon$ sufficiently small so that 
\[\sin |\theta| + Az'(s) \leq c\left[\cos \theta - A z'(s)\right], \text{ i.e. } A \leq \frac{c\cos \theta -\sin |\theta|}{1+c}.\]
Therefore combining with (\ref{x'positive})
\begin{equation}\label{y'x'}
|y'(s)| \leq cx'(s).
\end{equation}
With $y(0) = x(0) = 0$, one has 
\begin{eqnarray}\label{yx}
|y(s)| \leq cx(s).
\end{eqnarray}

(4) Concavity. This is the last estimate of $\epsilon$. 

Suppose $\gamma(s)$ leaves $S$ at a switch point when $s = s_0$ and dives into the interior of $M$ for increasing $s$ until it enters $S$ again at $s = s_1$. 

Since $x'(s) > 0$ within the $\epsilon$-ball of $p$, then $x(s)$ has a $C^1$-inverse function $s(x)$ for $s \in [0, s_1]$. Therefore we can express $y(s)$ as 
\[y(s) = y(s(x)) = \alpha(x),\]
where $\alpha$ is a $C^1$-function and $\alpha(0) = \frac{d\alpha}{dx} = 0$. 

Consider the intersection of the two-dimensional plane $y = y_0 + T(x - x_0)$ with the surface $z = g(x, y)$, where $(x_0, y_0) = (x(s_0), y(s_0))$, and $T = \frac{d\alpha}{dx}(x_0)$. With (\ref{y'x'})
\[|T| = \left|\frac{y'(s_0)}{x'(s_0)} \right| \leq c. \]
Set
\[f(x) = g(x, y_0 + T(x - x_0)),\]then
\[\frac{d^2 f}{d x^2}(x) = g_{xx} + 2g_{xy}T + g_{yy}T^2,\] where
\begin{eqnarray*}
g_{xx} &=& 2a + (6xc + 6x^2c_x + x^3c_{xx}) + (2yd+ 4xyd_x+ x^2yd_{xx}) +(2y^2e_x +xy^2e_{xx}) + y^3f_{xx}\\
&=& 2a + x(\ldots) + y(\ldots) \\
g_{xy} &=& (3x^2c_y + x^3c_{xy}) + (2xd + 2xyd_y + x^2d_x + x^2yd_{xy}) + (2ye + y^2e_y + 2xye_x + xy^2e_{xy}) + x^3f_{xy} \\
&=& x(\ldots) + y(\ldots) \\
g_{yy} &=& -2b + x^3c_{yy} + (2x^2d_y + x^2yd_{yy}) + (2xe + 4xye_y + xy^2e_{yy}) + (6yf + 6y^2f_y + y^3f_{yy}) \\
&=&-2b + x(\ldots) + y(\ldots)
\end{eqnarray*}
So with $y(s) = y_0 + T(x - x_0)$ for $s \in [s_0, s_1]$ and (\ref{y'x'}), (\ref{yx})
\begin{eqnarray*}
\frac{d^2 f}{d x^2}(x) &=& 2a + x(\ldots) + y(\ldots) + 2Tx(\ldots) + 2Ty(\ldots) -2bT^2 + T^2x(\dots) + T^2y(\ldots), \text{ so }\\
\frac{d^2 f}{d x^2}(x(s))&\geq& 2a - x(s)C_1 - cx(s)C_2 - 2cx(s)C_3 - 2c^2x(s)C_4 - 2bc^2 - c^2x(s)C_5 - c^3x(s)C_6, 
\end{eqnarray*}
where $C_1, \ldots, C_6$ are constants bounding the terms inside the corresponding parentheses. 
By assumption $c < \tan \theta_0 = \displaystyle \sqrt{\frac{a}{b}}$, we can choose $\epsilon$ sufficiently small such that the right side of the above inequality is positive, i.e., 
\[x(s)C_1 + cx(s)C_2 + 2cx(s)C_3 + 2c^2x(s)C_4 + c^2x(s)C_5 + c^3x(s)C_6 < 2a - 2bc^2. \]
Therefore $\frac{d^2 f}{d x^2}(x(s)) > 0$ for $s \in [s_0, s_1]$. It implies that $f'(x(s))$ is increasing as $x(s)$ increases from $x(s_0)= x_0$ to $x(s_1) = x_1$. On the other hand, since the interior line segment is tangent to $S$ at the two endpoints, one must have 
\[f'(x_0) = f'(x_1),\]
a contradiction. Therefore if $\gamma$ leaves $S$ at the switch point $\gamma(s_0)$, the geodesic arc beyond this point is a line segment never returning to $S$ again. So $\gamma$ has at most one switch point within the $\epsilon$-ball of $p$, as desired. 
\end{proof}

\begin{proposition}\label{propositon2}
In the positive direction of the line $y = \displaystyle \sqrt{\frac{a}{b}}$, there exist an $\eta > 0$ and an $\epsilon > 0$ such that if $\gamma'(0)$ lies in the wedge $[\theta_0 - \eta, \theta_0 + \eta]$, then $\gamma$ has at most two switch points within the $\epsilon$-ball of $p$. 
\end{proposition}

\begin{proof}

(1) Set up the frame. 

Let's rotate the $(x, y)$-plane so that the $x$-axis points in the positive direction of the line $y = \displaystyle \sqrt{\frac{a}{b}}$ by the matrix
\begin{equation*}
\left[ \begin{matrix}
\cos \theta_0 & -\sin \theta_0 \\
\sin \theta_0 & \cos \theta_0 
\end{matrix}\right].
\end{equation*}
Thus with respect to this new coordinate system and in connection with (\ref{gabcdef}) the surface $S$ near $p$ can be parametrized by 
\begin{eqnarray*}
&&g(\cos \theta_0 x -\sin \theta_0 y, \sin \theta_0 x  + \cos \theta_0 y) \\
&=& a(\cos \theta_0 x -\sin \theta_0 y)^2 - b(\sin \theta_0 x  + \cos \theta_0 y)^2 + \\
&& (\cos \theta_0 x -\sin \theta_0 y)^3c + (\cos \theta_0 x -\sin \theta_0 y)^2(\sin \theta_0 x + \cos \theta_0 y)d\\
&&+ (\cos \theta_0 x -\sin \theta_0 y)(\sin \theta_0 x  + \cos \theta_0 y)^2 e + (\sin \theta_0 x  + \cos \theta_0 y)^3 f\\
&=& -2\sqrt{ab}xy + (a-b)y^2 \\
&&+ x^3\left[(\cos \theta_0)^3 c+ (\cos \theta_0)^2\sin \theta_0d + \cos \theta_0(\sin \theta_0)^2 e+ (\sin \theta_0)^3 f \right] \\
&& + x^2y \left[3(\cos \theta_0)^2\sin \theta_0c +\cos \theta_0(1- 3(\sin \theta_0)^2)d + \sin \theta_0 (3 (\cos \theta_0)^2-1) e + 3 (\sin \theta_0)^2 \cos \theta_0 f \right]\\
&& + xy^2 \left[ -3 \cos \theta_0(\sin \theta_0)^2 c + ((\cos \theta_0)^3 - 2(\cos \theta_0)^2\sin \theta_0(d+e) + 3\sin \theta_0(\cos\theta_0)^2 f\right] \\
&& + y^3 \left[ -(\sin \theta_0)^3c + (\sin \theta_0)^2 \cos \theta_0 d - \sin \theta_0 (\cos \theta_0)^2e + (\cos \theta_0)^3 f\right]
\end{eqnarray*}
where $c, d, e, f$ are evaluated at $ (\cos \theta_0 x -\sin \theta_0 y, \sin \theta_0 x  + \cos \theta_0 y)$. 
For convenience, let's still use $g(x, y)$ to denote the new parametrization
\begin{eqnarray}\label{gtheta0}
g(x, y) &=& -2\sqrt{ab}xy + (a-b)y^2 \\ \notag
&&+x^3 c(x, y) + x^2y d(x, y) + xy^2e(x, y) + y^3f(x, y), \notag
\end{eqnarray}
where $c, d, e, f$ are the (new) functions inside the corresponding brackets. 

(2) Notice that within the $\epsilon$-ball of $p$, we have $|x|, |y| \leq \epsilon$ no matter how we rotate the $(x, y)$-plane, because the distance to the origin is fixed. So if $\gamma(s)$ is within the $\epsilon$-ball of $p$, then we have
\[ |\gamma(s)| \leq \epsilon \Rightarrow |x(s)| \leq \epsilon, |y(s)| \leq \epsilon, |z(s)| \leq \epsilon.\]

(3) Concavity with respect to the new frame. If one moves slightly above the $x$-axis, namely in the first quadrant, then the surface is concave downward; if one moves slightly below the $x$-axis, namely in the fourth quadrant, then the surface is concave upward. 

(4) If $\gamma'(0)$ points in the positive $x$-direction.

One can rewrite the $g(x, y)$ again as follows:
\begin{equation}\label{gdelta0}
g(x, y) = x^N h(x, y) + xyi(x, y) + y^2j(x, y),
\end{equation}
where $N \geq 3$, $h(0, 0) \neq 0$, $i(0, 0)= -2\sqrt{ab}$, and $j(0, 0) = a-b$. According to the theorem for a fixed direction (cite here), if $\gamma'(0) = \frac{\partial}{\partial x}$ then there exists an $\epsilon > 0$ such that $\gamma$ has at most one switch point before leaving the $\epsilon$-ball of $p$. 

(5) $\gamma'(0)$ is in the wedge $[-\eta, \eta]$ where $0 < \eta < \min(2\theta_0, \pi - 2\theta_0)$. This is the first estimate of $\eta$. 

Let $\delta \in [-\eta, \eta]$ and $\delta \neq 0$. So $\gamma'(0)$ does NOT point in the positive $x$-direction.
If we rotate the $(x,y)$-plane according to the angle $\delta$, the Taylor expansion of $g(\cos \delta x-\sin \delta y, \sin \delta x + \cos \delta y)$, denoted as $g_{\delta}(x, y)$, has a nonzero $x^2$ term, because with (\ref{gtheta0})
\begin{eqnarray}\label{gdelta}
g_{\delta}(x, y) &=& -2\sqrt{ab}(\cos \delta x-\sin \delta y)(\sin \delta x + \cos \delta y) + (a - b)(\sin \delta x + \cos \delta y)^2 \\ \notag
&&+ (\cos \delta x-\sin \delta y)^3c + (\cos \delta x-\sin \delta y)^2(\sin \delta x + \cos \delta y)d \\ \notag
&&+ (\cos \delta x-\sin \delta y)(\sin \delta x + \cos \delta y)^2e + (\sin \delta x + \cos \delta y)^3 f\\ \notag
&=&x^2 \left[b\frac{\cos^2 (\theta_0 + \delta)}{\cos^2 \theta_0} - b \right] + xy \left[ -2\sqrt{ab}\cos 2\delta + (a-b)\sin 2\delta \right]\\ \notag
&& + y^2 \left[ \sqrt{ab} \sin 2\delta + (a-b) \cos^2 \delta \right] \\ \notag
&&+ x^3\left[(\cos \delta)^3 c+ (\cos \delta)^2\sin \delta d + \cos \delta(\sin \delta)^2 e+ (\sin \delta)^3 f \right] \\ \notag
&& + x^2y \left[3(\cos \delta)^2\sin \delta c +\cos \delta(1- 3(\sin \delta)^2)d + \sin \delta (3 (\cos \delta)^2-1) e + 3 (\sin \delta)^2 \cos \delta f \right]\\ \notag
&& + xy^2 \left[ -3 \cos \delta(\sin \delta)^2 c + ((\cos \delta)^3 - 2(\cos \delta)^2\sin \delta(d+e) + 3\sin \delta(\cos\delta)^2 f\right] \\ \notag
&& + y^3 \left[ -(\sin \delta)^3c + (\sin \delta)^2 \cos \delta d - \sin \delta (\cos \delta)^2e + (\cos \delta)^3 f\right] \notag
\end{eqnarray}
where $c, d, e, f$ are evaluated at $(\cos \delta x-\sin \delta y, \sin \delta x + \cos \delta y)$.
Notice that if the coefficient of $x^2$ is zero, then 
\[b\frac{\cos^2 (\theta_0 + \delta)}{\cos^2 \theta_0} - b = 0 \Rightarrow \cos(\theta_0 + \delta) = \pm \cos \theta_0 \Rightarrow \delta = 0, -2\theta_0, \text{ or } \pi - 2\theta_0,    \]
a contradiction. Therefore we could rewrite $g_{\delta}(x, y)$ to include the angle $\delta = 0$. 
\begin{eqnarray}\label{gdeltanew}
g_{\delta}(x, y) &=& a_2(\delta)x^2 + a_3(\delta)x^3 + \cdots + a_{N-1}(\delta)x^{N-1} \\ \notag
&&+ x^Nh_{\delta}(x, y) + xyi_{\delta}(x, y) + y^2j_{\delta}(x, y), \notag
\end{eqnarray}
where $a_2(\delta), \ldots, a_{N-1}(\delta)$ are constant coefficients of $x^2, \ldots, x^{N-1}$, respectively. Moreover at $\delta = 0$, these constants vanish, and $h_0(x, y) = h(x, y)$, $i_0(x, y) = i(x, y)$, $j_0(x, y) = j(x, y)$. 

(6) There exists an $\epsilon$ that works for all $\delta \in [-\eta, \eta]$. This is the second estimate of $\eta$ (in step 4). There are four different cases depending on the sign of $a_2(\delta)$ and $h(0, 0)$. 

{\bf Case 1: $a_2(\delta) > 0$, $h(0, 0) > 0$.}

When $a_2(\delta) >0 $, the angle $\delta < 0.$
There is actually a relationship between $a_2(\delta), \ldots, a_{N-1}(\delta)$.
\begin{eqnarray*}
a_2(\delta) &=& \frac{b}{\cos^2 \theta_0}[\sin^2 (\theta_0) - \sin^2(\theta_0 + \delta)] \\
&=&  \frac{b}{\cos^2 \theta_0}[\sin (\theta_0) + \sin(\theta_0 + \delta)][\sin (\theta_0) - \sin(\theta_0 + \delta)]\\
&=& \frac{b}{\cos^2 \theta_0}[\sin (\theta_0) + \sin(\theta_0 + \delta)][\sin (\theta_0) -\sin(\theta_0)\cos \delta - \cos(\theta_0) \sin \delta] \\
&=& \frac{b}{\cos^2 \theta_0}[\sin (\theta_0) + \sin(\theta_0 + \delta)][- \cos(\theta_0) \sin \delta + \sin (\theta_0)(1 -\cos \delta)]\\
&\geq& \frac{b}{\cos^2 \theta_0}\sin (\theta_0)[- \cos(\theta_0) \sin \delta] = b\tan(\theta_0) |\sin \delta|.
\end{eqnarray*}
Furthermore, there are constants $c_i$ and $M$ sufficiently large such that
\begin{eqnarray*}
a_3(\delta)&=& c_1\cos^2 \delta \sin \delta + c_2 \cos \delta \sin^2 \delta + c_3 \sin^3 \delta \\
a_4(\delta)&=& c_4\cos^3 \delta \sin \delta + c_5 \cos^2 \delta \sin^2 \delta + c_6 \cos \delta \sin^3 \delta + c_7 \sin^4 \delta \\
&& \vdots\\
a_{N-1}(\delta) &=& c_8 \cos^{N-2} \delta \sin \delta + \cdots + c_9 \sin^{N-1} \delta \\
\Rightarrow && |a_3(\delta)|, \ldots, |a_{N-1}(\delta)| \leq M |\sin \delta| \\
\end{eqnarray*}

(i) Let's continue with the $\epsilon$ in (4). Then $x'(s) \geq \frac{1}{2}$ if $\epsilon$ is chosen small enough. Here we are assuming that if $\gamma(s)$ is within the $\epsilon$-ball of $p$, then for every $0 \leq \sigma \leq s$, $\gamma(\sigma)$ also lies within the $\epsilon$-ball of $p$. 

(ii) Let $\gamma(s) = (x(s), y(s), z(s))$. If $\gamma(s) \in S$, then the normal vector to $S$ at $\gamma(s)$ is 
\[N(s) = (-(g_{\delta})_x(x(s), y(s)), -(g_{\delta})_y(x(s), y(s)), 1).\]
 From (\ref{gdeltanew}) it follows that 
\begin{eqnarray*}
(g_{\delta})_x (x, y) &=& 2a_2(\delta)x + \cdots + (N-1)a_{N-1}(\delta)x^{N-2} + Nx^{N-1}h_{\delta} + x^N(h_{\delta})_x \\
&& + yi_{\delta} + xy(i_{\delta})_x + y^2(j_{\delta})_x.
\end{eqnarray*}
There exists a positive constant $A$ such that 
\[|(g_{\delta})_x (x, y)| \leq A \text{ if } |x|, |y| \leq \epsilon \text{ and } \delta \in [- \eta,  \eta].\]
Moreover $A \to 0$ as $\epsilon \to 0$ for a fixed $\eta$. 
Let $s$ be such that $\gamma''(s)$ exists, then $\gamma''(s) = z''(s)N(s)$. This implies that 
\[x''(s) = -z''(s)(g_{\delta})_x(x(s), y(s)) \Rightarrow |x''(s)| \leq Az''(s). \]
Here $z''(s) \geq 0$, because within the $\epsilon$-ball of $p$ the outward normal vector to $S$ has a positive $z$-coordinate of 1 and $\gamma''(s)$ directs outward on a boundary segment on $S$. Indeed $\gamma(s)$ is a locally shortest path. If $\gamma(s)$ lies on the surface of $M$, its acceleration exists everywhere except at the switch points and is outward normal to the surface. 

On the other hand, if $\gamma(s)$ lies on a line segment in the interior of $M$, then the acceleration $\gamma''(s)$ is zero. So the previous inequality still holds. Thus
\begin{eqnarray}\label{x'positive2}
x'(s) = x'(0) + \displaystyle \int_0^s x''(\sigma) d\sigma \geq 1 - A\displaystyle \int_0^s z''(\sigma) d\sigma = 1 - Az'(s). 
\end{eqnarray}

Next let's approximate $z'(s)$. If $\gamma(s) \in S$, using (\ref{gdeltanew}) one has 
\begin{eqnarray*}
z'(s) &=& 2a_2(\delta)x(s)x'(s)  + \cdots + (N-1)a_{N-1}(\delta)x(s)^{N-2}x'(s) \\
&&+ Nx^{N-1}x'(s) h_{\delta} + x^N[(h_{\delta})_xx'(s) + (h_{\delta})_yy'(s)] \\
&& + x'(s)y(s)i_{\delta} +x(s)y'(s)i_{\delta}+ x(s)y(s)[(i_{\delta})_xx'(s) + (i_{\delta})_yy'(s)] \\
&&+ 2y(s)y'(s)j_{\delta} + y(s)^2[(j_{\delta})_xx'(s) + (j_{\delta})_yy'(s)].
\end{eqnarray*}
Since $\gamma(s)$ is parametrized by arc length, $|x'(s)|$ and $|y'(s)|$ are no more than 1. So there exists a positive constant $B$ such that
\[|z'(s)| \leq B \text{ if } |x(s)|, |y(s)| \leq \epsilon \text{ and } \delta \in [- \eta, \eta], \]
where $B \to 0$ as $\epsilon \to 0$ for a fixed $\eta$. On the other hand, if $\gamma(s)$ is within an interior line segment, then $\gamma'(s)$ is constant and equal to the value at the endpoints. Therefore $|z'(s)|$ is still bounded by $B$. Thus we can choose $\epsilon$ small enough so that $B \leq \frac{1}{2A}$. It follows from (\ref{x'positive2}) that 
\begin{equation}\label{x'1/2}
x'(s) \geq \frac{1}{2},
\end{equation}
if $\epsilon$ is chosen sufficiently small.
Notice that this $\epsilon$ works for all $\delta$ because we can bound $h_{\delta}$, $i_{\delta}$, $j_{\delta}$ and their partial derivatives uniformly. 

(iii) Approximate $z(s), z'(s), y(s), y'(s)$. If $\gamma(s) \in S$, then 
\[y''(s) = -z''(s) (g_{\delta})_y(x(s), y(s)) \Rightarrow |y''(s)| \leq Az''(s),\]
where $A \to 0$ if $\epsilon \to 0$ and $A$ does not depend on $\delta$. If $\gamma(s) \not \in S$, then $\gamma(s)$ is within an interior line segment, so the inequality still holds. With respect to the new frame after rotating the $(x, y)$-plane by $\delta$, $\gamma'(0) = \frac{\partial}{\partial x}$ and so with $z(0) = z'(0) = y(0) = y'(0) = 0$ one deduces that
\begin{equation}\label{y'y}
|y'(s)| \leq Az'(s), |y(s)| \leq Az(s). 
\end{equation}
Now let's estimate $z(s)$. If $\gamma(s) \in S$, with (\ref{gdeltanew}) and (\ref{y'y})
\begin{eqnarray*}
z(s) &\leq &  a_2(\delta)x(s)^2 + |a_3(\delta)|x(s)^3 + \cdots + |a_{N-1}(\delta)|x(s)^{N-1} \\
&&+ x(s)^N|h_{\delta}(x(s), y(s))| + x(s)Az(s)|i_{\delta}(x(s), y(s))| + A^2z(s)^2|j_{\delta}(x(s), y(s))| 
\end{eqnarray*}
First, 
\begin{eqnarray*}
&& |a_3(\delta)|x(s) + \cdots + |a_{N-1}(\delta)|x(s)^{N-3} \\
&\leq& M|\sin \delta|[x(s) + x(s)^2 + \cdots + x(s)^{N-3}] \\
&\leq& b\tan(\theta_0) |\sin \delta| \leq a_2(\delta),
\end{eqnarray*}
if we choose $\epsilon$ sufficiently small. 
Second, 
\[|h_{\delta}(x(s), y(s))| \leq 2h(0, 0),\]
if $\eta$ and $\epsilon$ are sufficiently small. 
Third, 
\[ x(s)Az(s)|i_{\delta}(x(s), y(s))| + A^2z(s)^2|j_{\delta}(x(s), y(s))| \leq x(s)Az(s)C_1 + A^2z(s)^2C_2,\]
where $C_1, C_2$ are some constants and we can choose $\epsilon$ small enough so that 
\[Ax(s)C_1 + A^2z(s)C_2 \leq C_3 < 1.\]
Therefore there exists a positive constant $C$ such that 
\begin{eqnarray}\label{zy}
&& z(s) \leq (2a_2(\delta)x(s)^2 + 2h(0, 0)x(s)^N)C \\ \notag
&\Rightarrow& |y(s)| \leq (2a_2(\delta)x(s)^2 + 2h(0, 0)x(s)^N)AC \notag
\end{eqnarray}
On the other hand, 
\begin{eqnarray*}
z'(s) &\leq& 2a_2(\delta)x(s) + \cdots + (N-1)|a_{N-1}(\delta)|x(s)^{N-2} \\
&&+ Nx(s)^{N-1}|h_{\delta}| + x(s)^N[|(h_{\delta})_x| + |(h_{\delta})_y|] \\
&&+ |y(s)||i_{\delta}| + x(s)Az'(s)|i_{\delta}| + x(s)|y(s)|[|(i_{\delta})_x| + |(i_{\delta})_y|] \\
&&+ 2y(s)Az'(s)|j_{\delta}| + |y(s)|^2[|(j_{\delta})_x| + |(j_{\delta})_y|]\\
&\leq& 3a_2(\delta)x(s)+ 2Nx(s)^{N-1}h(0, 0) +  (2a_2(\delta)x(s)^2 + 2h(0, 0)x(s)^N)ACC_4 + C_3z'(s)\\
&\leq& 4a_2(\delta)x(s) + 4Nx(s)^{N-1}h(0, 0) + C_3z'(s),
\end{eqnarray*}
if $\eta$, $\epsilon$ are sufficiently small and $C_3 < 1$, $C_4$ are some constants. So there is a constant $C$ such that 
\begin{eqnarray}\label{z'y'}
&&z'(s) \leq (4a_2(\delta)x(s) + 4Nx(s)^{N-1}h(0, 0))C  \\ \notag
&\Rightarrow& y'(s) \leq (4a_2(\delta)x(s) + 4Nx(s)^{N-1}h(0, 0))AC. \notag
\end{eqnarray}
It follows that if $y = \alpha(x)$ and $T(s) = \alpha'(x(s))$, then
\[|T(s)| = \left|\frac{y'(s)}{x'(s)} \right| \leq 2|y'(s)| \leq (8a_2(\delta)x(s) + 8Nx(s)^{N-1}h(0, 0))AC.\]

If $\gamma(s) \not \in S$, then $\gamma(s)$ is in some line segment where $y(s) = y(s_0) + T(s_0)(x(s) - x(s_0))$ for some switch point at $s_0 < s$. Since $x(s)$ is increasing, it follows that
\begin{eqnarray*}
|y(s)| &\leq& |y(s_0)| + |T(s_0)| [|x(s)|+ |x(s_0)|] \\
&\leq&  (2a_2(\delta)x(s_0)^2 + 2h(0, 0)x(s_0)^N)AC \\
&&+  2(4a_2(\delta)x(s_0) + 4Nx(s_0)^{N-1}h(0, 0))AC [|x(s)|+ |x(s_0)|] \\
&\leq& (2a_2(\delta)x(s)^2 + 2h(0, 0)x(s)^N)AC + 2(4a_2(\delta)x(s) + 4Nx(s)^{N-1}h(0, 0))AC \cdot 2x(s)\\
&\leq& (18a_2(\delta)x(s)^2 + 18Nh(0,0)x(s)^N)AC.
\end{eqnarray*}
The same inequality for $T(s)$ still holds as before because $T(s) = T(s_0)$. 

(iv) Now we are ready to show that there is at most one switch point. Set $f(x) = g(x, y_0 + T(x-x_0))$ where $(x_0, y_0) = (x(s_0), y(s_0))$, $T = T(s_0)$, and $s_0$ is arbitrary. Then 
\begin{eqnarray*}
f''(x_0) &=& g_{xx}(x_0, y_0) + 2g_{xy}(x_0, y_0)T + g_{yy}(x_0, y_0)T^2,
\end{eqnarray*}
where
\begin{eqnarray*}
g_{xx}(x_0, y_0) &=& 2a_2(\delta) + 6a_3(\delta)x_0 + \cdots + (N-1)(N-2)a_{N-1}(\delta)x_0^{N-3} \\
&& + x_0^{N-2}[N(N-1)h_{\delta} + 2Nx_0(h_{\delta})_x + x_0^2(h_{\delta})_{xx}] \\
&& + y_0[2(i_{\delta})_x + x_0(i_{\delta})_{xx} + y_0(j_{\delta})_{xx}].
\end{eqnarray*}
First, for $\epsilon$ sufficiently small,
\begin{eqnarray*}
&&|6a_3(\delta)x_0 + \cdots + (N-1)(N-2)a_{N-1}(\delta)x_0^{N-3}| \\
&\leq& M|\sin \delta| (6x_0 + \cdots + (N-1)(N-2)x_0^{N-3}) \\
&\leq& b\tan(\theta_0)|\sin \delta| \leq a_2(\delta).
\end{eqnarray*}
Second, for $\eta$ and $\epsilon$ sufficiently small, 
\[N(N-1)h_{\delta} + 2Nx_0(h_{\delta})_x + x_0^2(h_{\delta})_{xx} \geq \frac{1}{2}N(N-1)h(0, 0).\]
Third, there are constants $C_5, C_6, C_7$ such that 
\[|2(i_{\delta})_x + x_0(i_{\delta})_{xx} + y_0(j_{\delta})_{xx}| \leq C_5, |2g_{xy}(x_0, y_0)| \leq C_6, |g_{yy}(x_0, y_0)T| \leq C_7.\]
So 
\begin{eqnarray*}
f''(x_0) &\geq& 2a_2(\delta) - a_2(\delta) + x_0^{N-2}\frac{1}{2}N(N-1)h(0, 0) \\
&&- (18a_2(\delta)x_0^2 + 18Nh(0,0)x_0^N)ACC_5 - (8a_2(\delta)x_0+ 8Nx_0^{N-1}h(0, 0))ACC_6 \\
&&- (8a_2(\delta)x_0+ 8Nx_0^{N-1}h(0, 0))ACC_7\\
&=&a_2(\delta)[1-18ACC_5x_0^2 - 8ACC_6x_0 - 8ACC_7x_0]\\
&& + x_0^{N-2}h(0, 0)[\frac{1}{2}N(N-1) - 18ACC_5Nx_0^2 - 8ACC_6Nx_0 - 8ACC_7Nx_0].
\end{eqnarray*}
Therefore if $\epsilon$ is sufficiently small, for all $|x_0| \leq \epsilon$, one has
\begin{eqnarray*}
1-18ACC_5x_0^2 - 8ACC_6x_0 - 8ACC_7x_0 > 0, \\
\frac{1}{2}N(N-1) - 18ACC_5Nx_0^2 - 8ACC_6Nx_0 - 8ACC_7Nx_0 > 0, 
\end{eqnarray*}
implying that $\gamma$ can't have a line segment within the $\epsilon$-ball. Thus in the case when $a_2(\delta) > 0$ and $h(0, 0) > 0$, $\gamma$ has at most one switch point within the $\epsilon$-ball for any $\delta \in [-\eta, \eta]$. 

{\bf Case 2: $a_2(\delta) < 0$, $h(0, 0) < 0$. }

Suppose for the sake of contradiction, $\gamma'(0)$ is in the direction of $\frac{\partial}{\partial x}$. Then we can approximate $y(s)$, $T(s)$ as in the paper (cite here) to show that $\gamma$ has no switch point close to the origin. It follows that $\gamma$ has to be on the surface initially. 
Given any two points  $\gamma(s_1)$, $\gamma(s_2)$ on the geodesic close to the origin, one can show that the line segment connecting them is actually in the interior of $M$ contradicting that $\gamma$ is locally shortest. 

Let $(x_1, y_1) = (x(s_1), y(s_1))$ and $(x_2, y_2) = (x(s_2), y(s_2))$. Then for $t \in [0, 1]$, set
\[f(t) = g_{\delta}(x_1 + t(x_2 - x_1), y_1 + t(y_2- y_1)).\] So
\[f''(t) = (g_{\delta})_{xx}(x_2-x_1)^2 + 2(g_{\delta})_{xy}(x_2-x_1)(y_2-y_1) + (g_{\delta})_yy(y_2-y_1)^2. \] By the mean value theorem, 
\[x_2 - x_1 = (s_2 - s_1)x'(\tilde{s}), y_2 - y_1 = (s_2 - s_1)y'(\hat{s}),\] for some $\tilde{s}, \hat{s}$ in $(s_1, s_2)$. Let
\[T = \frac{y_2 - y_1}{x_2 - x_1} = \frac{y'(\hat{s})}{x'(\tilde{s})} \to 0 \text{ as } s_1, s_2 \to 0, \] 
because $\gamma'(0) = (1, 0, 0)$. 
Therefore
\[f''(t) = (x_2 - x_1)^2 [(g_{\delta})_{xx} + 2(g_{\delta})_{xy}T + (g_{\delta})_{yy} T^2]. \]
Since $g_{xx}(x, y) \to 2a_2(\delta)$ as $x, y \to 0$, then 
\[f''(t) \to (x_2-x_1)^2[2a_2(\delta)]\]
as $s_1, s_2 \to 0$ and for every $t \in [0, 1]$. Thus the shortest path two points on $\gamma$ close to the orgin is the line segment in between which lies stricly below the surface, a contradiction. So $\gamma$ is a straight line initially. 

The surface in the $(x, z)$-plane is the curve with equation
\begin{eqnarray*}
z(x) &=& g_{\delta}(x, 0) \\
&=& a_2(\delta)x^2 + \cdots + a_{N-1}(\delta)x^{N-1} + x^N h_{\delta}(x, 0),
\end{eqnarray*}
which implies that
\[z'(0) = 0, z''(0) = 2a_2(\delta) < 0.\]
So the slope of the line segment must be negative. If the line segment re-enters the surface at some switch point, then the surface can't be concave downward there. Otherwise the line lies above the surface. 
\begin{eqnarray*}
z'(x) &=& 2a_2(\delta)x + \cdots + (N-1)a_{N-1}(\delta)x^{N-2} + Nx^{N-1}h_{\delta}(x, 0) + x^N(h_{\delta})_x(x, 0) \\
z''(x) &=& 2a_2(\delta) + 6a_3(\delta)x + \cdots + (N-1)(N-2)a_{N-1}(\delta)x^{N-3} \\
&&+ N(N-1)x^{N-2}h_{\delta}(x, 0) +2 Nx^{N-1}(h_{\delta})_x(x, 0) + x^N(h_{\delta})_{xx}(x, 0)
\end{eqnarray*}

When $a_2(\delta) < 0$, the angle $\delta > 0$. As before, one has 
\begin{eqnarray*}
- a_2(\delta) &=&\frac{b}{\cos^2 \theta_0}[\sin(\theta_0) + \sin (\theta_0 + \delta)][\sin (\theta_0 + \delta) - \sin (\theta_0)]\\
&=& \frac{b}{\cos^2 \theta_0}[\sin(\theta_0) + \sin (\theta_0 + \delta)][\delta \cos (\theta)],
\end{eqnarray*}
where  $\theta \in (\theta_0, \theta_0 + \delta)$ by the mean value theorem. For $\eta$ sufficiently small, we have
\[\delta = \frac{\delta}{\sin \delta} \sin \delta \geq \frac{1}{2} \sin \delta, \]
since $\frac{\delta}{\sin \delta} \to 1$ as $\delta \to 0$. Thus 
\[- a_2(\delta) \geq \frac{b}{\cos^2 \theta_0} \sin (\theta_0) \frac{1}{2} \sin \delta \cos(\theta_0 + \eta) \geq \frac{1}{4}b\tan (\theta_0) \sin \delta,\] if $\eta$ is sufficiently close to 0. 
Furthermore, there is $M$ sufficiently large such that 
\[|a_3(\delta)|, \ldots, |a_{N-1}(\delta)| \leq M \sin \delta.\]
It follows that 
\begin{eqnarray*}
&&|6a_3(\delta)x + \cdots + (N-1)(N-2)a_{N-1}(\delta)x^{N-3}| \\
&\leq& M \sin \delta (6x + \cdots + (N-1)(N-2)x^{N-3}) \\
&\leq& \frac{1}{4}b\tan (\theta_0) \sin \delta \leq -a_2(\delta),
\end{eqnarray*}
for all $|x| \leq \epsilon$ if $\epsilon$ is sufficiently small. 
Moreover, if $\eta$ and $\epsilon$ are sufficiently small, then
\[N(N-1)h_{\delta}(x, 0) +2 Nx(h_{\delta})_x(x, 0) + x^2(h_{\delta})_{xx}(x, 0) \leq \frac{1}{2}N(N-1)h(0, 0),\]
Therefore 
\[z''(x) \leq 2a_2(\delta) -a_2(\delta) + \frac{1}{2}N(N-1)h(0, 0) x^{N-2} < 0. \]
So $\gamma$ has no switch point unless it terminates at a point on the surface. 
As a summary in the case when $a_2(\delta) < 0$ and $h(0, 0) < 0$, $\gamma$ is either a straight line exiting the $\epsilon$-ball or a line segment terminating at some point on the surface within the $\epsilon$-ball. 

{\bf Case 3: $a_2(\delta) < 0$, $h(0, 0) > 0$.}

Since $a_2(\delta) < 0$, $\gamma$ is initially a straight line just as shown in Case 2. Suppose the angle between $\gamma'(0)$ and the positive $x$-axis is $-\beta$. The line either terminates at some point on the surface, or exits the $\epsilon$-ball, or enters the surface at some switch point at time $s_0$. Denote $x(s_0)$ as $x_0$. First, the intersection of the line with the surface at $\gamma(s_0)$ satisfies $-\tan(\beta) x_0 = g_{\delta}(x_0, 0)$ and so
\[-\tan(\beta) x_0 = a_2(\delta)x_0^2 + a_3(\delta)x_0^3 + \cdots + a_{N-1}(\delta)x_0^{N-1} + x_0^Nh_{\delta}(x_0, 0).\]
Next, the line is tangent to the surface at $\gamma(s_0)$, so $-\tan(\beta) = (g_{\delta})_x(x_0, 0)$ and
\[-\tan(\beta) = 2a_2(\delta)x_0 + 3a_3(\delta)x_0^2 + \cdots + (N-1)a_{N-1}(\delta)x_0^{N-2} + Nx_0^{N-1} h_{\delta}(x_0, 0) + x_0^N (h_{\delta})_x(x_0, 0).\]
Since $x_0 > 0$, the above two equalities imply the following:
\begin{eqnarray*}
&& a_2(\delta)x_0+ a_3(\delta)x_0^2 + \cdots + a_{N-1}(\delta)x_0^{N-2} + x_0^{N-1}h_{\delta}(x_0, 0) \\
&=& 2a_2(\delta)x_0 + 3a_3(\delta)x_0^2 + \cdots + (N-1)a_{N-1}(\delta)x_0^{N-2} + Nx_0^{N-1} h_{\delta}(x_0, 0) + x_0^N (h_{\delta})_x(x_0, 0)\\
\Rightarrow&& a_2(\delta)+ a_3(\delta)x_0 + \cdots + a_{N-1}(\delta)x_0^{N-3} + x_0^{N-2}h_{\delta}(x_0, 0) \\
&=& 2a_2(\delta) + 3a_3(\delta)x_0 + \cdots + (N-1)a_{N-1}(\delta)x_0^{N-3} + Nx_0^{N-2} h_{\delta}(x_0, 0) + x_0^{N-1} (h_{\delta})_x(x_0, 0)\\
\Rightarrow&& a_2(\delta) + 2a_3(\delta)x_0 + \cdots + (N-2)a_{N-1}(\delta)x_0^{N-3} + (N-1)x_0^{N-2}h_{\delta}(x_0, 0) + x_0^{N-1}(h_{\delta})_x(x_0, 0) \\
&=&0.
\end{eqnarray*}
Let $0 < c < 1$ be a constant to be determined later. Then for $\epsilon$ sufficiently small, 
\begin{eqnarray*}
&&|2a_3(\delta)x_0 + \cdots + (N-2)a_{N-1}(\delta)x_0^{N-3}| \\
&\leq& M \sin \delta [2x_0 + \cdots + (N-2)x_0^{N-3}] \\
&\leq& c \frac{1}{4}b\tan (\theta_0) \sin \delta \leq -ca_2(\delta).
\end{eqnarray*}
Therefore
\[-a_2(\delta) - 2a_3(\delta)x_0 - \cdots - (N-2)a_{N-1}(\delta)x_0^{N-3} \geq -a_2(\delta) + ca_2(\delta).\]
On the other hand, for $\eta$ and $\epsilon$ sufficiently small, we have
\[|(N-1)x_0^{N-2}h_{\delta}(x_0, 0) + x_0^{N-1}(h_{\delta})_x(x_0, 0)| \leq (N-1)x_0^{N-2} (1+c) h(0, 0) \]
Thus 
\begin{eqnarray}\label{x_0N-2}
&&(1-c)|a_2(\delta)| \leq (N-1)x_0^{N-2} (1+c) h(0, 0)\\ \notag
&\Rightarrow& x_0^{N-2} \geq \frac{(1-c)|a_2(\delta)|}{(N-1)(1+c)h(0, 0)}.
\end{eqnarray}

Now we are going to first shift our coordinates to have the origin at $\gamma(s_0) = (x_0, 0, z_0= -\tan(\beta)x_0)$ and then rotate the $(x, z)$-plane so that the $\gamma'(s_0)$ points in the positive $x$-axis. Let's use $(u, v, w)$ for the new coordinates, then with respect to the new frame 
\[x = \cos \beta u + \sin \beta w + x_0 , y = v, z = - \sin \beta u + \cos \beta w + z_0,\]
so the surface $z = g_{\delta}(x, y)$ satisfies the equation
\[ - \sin \beta u + \cos \beta w + z_0 = g_{\delta}(\cos \beta u + \sin \beta w + x_0, v). \]

Check that we can still solve for $w$ analytically in terms of $u, v$ within the $\epsilon$-ball. Taking the partial derivative of 
\[\sin \beta u - \cos \beta w - z_0 + g_{\delta}(\cos \beta u + \sin \beta w + x_0, v)\]
with respect to $w$ yields
\[-\cos \beta + (g_{\delta})_x  \sin \beta = \cos \beta [(g_{\delta})_x \tan \beta -1] = \cos \beta [-(g_{\delta})_x(g_{\delta})_x(x_0, 0) -1],\]
where $-\tan \beta = (g_{\delta})_x(x_0, 0)$ from before.
Since $(g_{\delta})_x(x, y) \to 0$ as $x, y \to 0$, for $\epsilon$ sufficiently small, 
\[-(g_{\delta})_x(g_{\delta})_x(x_0, 0) < 1.\]
Therefore there exists a real analytic function $k_\delta$ such that $w = k_\delta(u, v)$ with $k_\delta(0, 0) = 0$, $(k_\delta)_u(0, 0) = 0$, $(k_\delta)_v(0, 0) = 0$. 

Estimate $\gamma(s)$ in the new frame starting from the point $(x_0, 0, z_0)$. After replacing $s$ by $s-s_0$, $\gamma'(0)$ is equal to $\frac{\partial}{\partial u}$. 

1. $u'(s) \geq \frac{1}{2}$ if $\epsilon$ is chosen small enough.

By triangular inequality, $|u(s)|, |v(s)|, |w(s)|$ are less than or equal to 2$\epsilon$. If $\gamma(s) \in S$, then the normal vector to $S$ at $\gamma(s)$ is 
\[N(s) = (-(k_\delta)_u(u(s), v(s)), -(k_\delta)_v(u(s), v(s)), 1).\]
Since the lowest degree in $k_\delta$ is at least two, there exists a positive constant $A$ such that 
\[|(k_\delta)_u(u, v)| \leq A, |(k_\delta)_v(u, v)| \leq A, \text{ if } |u|, |v| \leq 2\epsilon,\]
where $A \to 0$ as $\epsilon \to 0$. Let $s$ be such that $\gamma''(s)$ exists, then $\gamma''(s) = w''(s)N(s)$. This implies that 
\[u''(s) = -w''(s)(k_\delta)_u(u(s), v(s)) \Rightarrow |u''(s)| \leq A w''(s).\]
Here $w''(s) \geq 0$ because within the $\epsilon$-ball of $p$ the surface $S$ has the parametrization $w = k_\delta(u, v)$ and thus the outward normal vector to $S$ has a positive $w$-coordinate of 1 and $\gamma''(s)$ is outward normal on a boundary segment in $S$. 
If $\gamma(s) \not \in S$, $\gamma''(s) = 0$ except at the switch points. 
Thus 
\[u'(s) = u'(0) + \int_0^s u''(\sigma) d\sigma \geq 1 - Aw'(s). \]

Next approximate $w'(s)$. If $\gamma(s) \in S$, then 
\[w(s) = k_\delta(u(s), v(s)) = u(s)^2a(u(s), v(s)) + u(s)v(s)b(u(s), v(s)) + v(s)^2c(u(s), v(s)),\] for some analytic functions $a, b, c$. 
Since $\gamma$ is parametrized by arclength, $|u'(s)| \leq 1$ and $|v'(s)| \leq 1$. Since each term in $w'(s)$ has either $u(s)$ or $v(s)$ and 
$u'(s)$ or $v'(s)$, there exists a positive constant $B$ such that 
\[|w'(s)| \leq B, \text{ if } |u(s)|, |v(s)| \leq 2\epsilon,\]
where $B \to 0$ as $\epsilon \to 0$. On the other hand, if $\gamma(s) \not \in S$, $\gamma'(s)$ is contant and equal to the value at the endpoints. Therefore $|w'(s)|$ is still bounded by $B$. Thus one can choose $\epsilon$ small enough so that $B < \frac{1}{2A}$. It follows that 
\[u'(s) \geq \frac{1}{2}.\]

2. Approximate $v'(s)$ and $v(s)$. 

If $\gamma(s) \in S$, then the normal vector to $S$ at $\gamma(s)$ is 
\[N(s) = (-(k_\delta)_u(u(s), v(s)), -(k_\delta)_v(u(s), v(s)), 1).\]
Let $s$ be such that $\gamma''(s)$ exists, then $\gamma''(s) = w''(s)N(s)$. This implies that
\[v''(s) = -w''(s)(k_\delta)_v(u(s), v(s)) \Rightarrow |v''(s)| \leq A w''(s).\]
Thus with $v'(0)=w'(0)=0$,
\[|v'(s)| \leq \int_0^s |v''(\sigma)| d\sigma \leq A\int_0^s w''(\sigma) d\sigma = Aw'(s).\]
With $v(0)=w(0)=0$, 
\[|v(s)| \leq Aw(s).\]

3. Coefficients of $k_\delta(u, v)$. Denote $k_\delta(u, v)$ as 
\[k_\delta(u, v) = b_2(\delta)u^2 + \cdots + b_{N-1}(\delta)u^{N-1} + u^N l_{\delta}(u, v) + uv m_{\delta}(u, v) + v^2 n_{\delta}(u, v),\]
where $b_2(\delta), \ldots, b_{N-1}(\delta)$ are constants and $l_{\delta}, m_{\delta}, n_{\delta}$ are analytic functions of $u, v$. Observe that for $n$ between 2 and $N-1$,
\[n!b_n(\delta) = \frac{\partial^n k_\delta}{\partial u^n}(0, 0).\]
The following lemma finds $\frac{\partial^n k_\delta}{\partial u^n}(u, v)$ for $n \geq 2$ by induction. 
\begin{lemma}\label{lem3}
Let $A$ be $\cos \beta + \sin \beta (k_\delta)_u(u, v)$. Then for each $n \geq 2$, 
\begin{eqnarray}\label{coefficients of k}
\cos \beta \frac{\partial^n k_\delta}{\partial u^n}(u, v) = \sum_{p = 0}^{n-1} \frac{\partial^{n-p} g_{\delta}}{\partial x^{n-p}} \sum_{I} c_I A^{n-p-|I|} (\frac{\partial A}{\partial u})^{i_1} (\frac{\partial^2 A}{\partial u^2})^{i_2}\cdots (\frac{\partial^p A}{\partial u^p})^{i_p},
\end{eqnarray}
where $I = (i_1, i_2, \ldots, i_p)$, $i_1 + 2i_2 + \cdots + pi_p = p$, $|I| = i_1 + i_2 + \cdots + i_p \leq n-p$, $c_I \geq 0$, and the partial derivatives of $g_{\delta}$ are evaluated at $(\cos \beta u + \sin \beta k_\delta(u, v) + x_0, v)$. 
\end{lemma}

\begin{proof}
When $n=2$, differentiating the equation 
\[-\sin \beta u + \cos \beta k_\delta(u, v) + z_0 = g_{\delta}(\cos \beta u + \sin \beta k_\delta(u, v) + x_0, v)\]
once with respect to $u$ gives
\[-\sin \beta + \cos \beta (k_\delta)_u = (g_{\delta})_x[\cos \beta + \sin \beta (k_\delta)_u] = (g_{\delta})_x A. \] Then taking the partial derivative with respect to $u$ once more gives 
\[\cos \beta (k_\delta)_{uu} = (g_{\delta})_{xx}A^2 + (g_{\delta})_x \partial_u A.\]
In (\ref{coefficients of k}) when $p = 0$, there is no $I$ so we have $c_0A^{2-0-0} = c_0 A^2$ where $c_0 = 1$;  when $p = 1$, there is only one $I = (1)$ so we have $c_1A^{2-1-1}(\frac{\partial A}{\partial u})^1 = c_1 \partial_u A$ where $c_1 = 1$. This coincides with the expression above. 

When $n \geq 2$, by inductive hypothesis we take the partial derivative of (\ref{coefficients of k}) with respect to $u$. The left-hand side is 
$\cos \beta \partial_u^{n+1} k_\delta$. The right-hand side consists of three parts due to the product rule.
\begin{eqnarray*}
&(1)& \sum_{p = 0}^{n-1} \frac{\partial^{n+1-p} g_{\delta}}{\partial x^{n+1-p}} \sum_{I} c_I A^{n+1-p-|I|} (\frac{\partial A}{\partial u})^{i_1} (\frac{\partial^2 A}{\partial u^2})^{i_2}\cdots (\frac{\partial^p A}{\partial u^p})^{i_p}, 
\end{eqnarray*}
where $n$ becomes $n+1$ and $p$ stays the same.
\begin{eqnarray*}
&(2)& \sum_{p = 0}^{n-1} \frac{\partial^{n-p} g_{\delta}}{\partial x^{n-p}} \sum_{I} c_I (n-p-|I|) A^{n-p-|I|-1} (\frac{\partial A}{\partial u})^{i_1+1} (\frac{\partial^2 A}{\partial u^2})^{i_2}\cdots (\frac{\partial^p A}{\partial u^p})^{i_p},
\end{eqnarray*}
where $n$, $p$, $i_1$ become $n+1$, $p+1$, $i_1 +1$, respectively. By letting $i_{p+1}$ be zero, one again has
\[ i_1 + 1 + 2i_2 + \cdots + pi_p + (p+1)i_{p+1} = p+ 1, n-p-|I|-1 = (n+1)-(p+1)-(|I| + 1).\]
\begin{eqnarray*}
&(3)& \sum_{p = 0}^{n-1} \frac{\partial^{n-p} g_{\delta}}{\partial x^{n-p}} \sum_{I} c_I A^{n-p-|I|} \sum_{i_j \neq 0} (\frac{\partial A}{\partial u})^{i_1} \cdots i_j (\frac{\partial^j A}{\partial u^j})^{i_j -1}(\frac{\partial^{j+1} A}{\partial u^{j+1}})^{i_j+1} \cdots (\frac{\partial^p A}{\partial u^p})^{i_p},
\end{eqnarray*}
where $n$, $p$ become $n+1$, $p+1$, respectively. When $j < p$, $i_j$ and $i_{j+1}$ are replaced by $i_j+1$ and $i_{j+1}+1$. By letting $i_{p+1}=0$ one has
\[\cdots + j(i_j-1) + (j+1)(i_{j+1}+1) + \cdots + (p+1)i_{p+1} = p+1,  \cdots + (i_j-1) + (i_{j+1} +1) + \cdots + i_{p+1} = |I|.\] 
On the other hand, when $j=p$, $i_p= i_{p+1} = 1$ and so 
\[p(i_p-1) + (p+1)i_{p+1} = p+1, (i_p-1)+ i_{p+1} = 1 = |I|.\]
It follows that $c_I$ are nonnegative integers and (\ref{coefficients of k}) is true. 
\end{proof}

\begin{corollary}\label{cor1}
The coefficient of $(g_{\delta})_x\partial_u^{n-1}A$  in (\ref{coefficients of k}) is always 1.
\end{corollary}

\begin{proof}
When $n = 2$, we've shown in the above lemma that the coefficient of $(g_{\delta})_x \partial_u A$ is 1. When $n \geq 2$, if $p = n-1$ then 
\[i_1 + 2i_2 + \cdots + (n-1)i_{n-1} = n-1 \text{ and } i_1 + i_2 + \cdots + i_{n-1} \leq 1\] imply that 
\[i_1 = \ldots = i_{n-2}=0 \text{ and } i_{n-1} =1.\]
There is only one term of $(g_{\delta})_x\partial_u^{n-1}A$ whose coefficient is 1 by induction. Taking its derivative with respect to $u$ yields
\[(g_{\delta})_{xx}A\partial_u^{n-1}A + (g_{\delta})_x\partial_u^{n}A,\]
so the coefficient of $(g_{\delta})_x\partial_u^{n}A$ is still 1 completing the induction. 
\end{proof}

\begin{corollary}\label{cor2}
The coefficient of $\frac{\partial^n g_{\delta}}{\partial x^n}A^n$ in (\ref{coefficients of k}) is always 1.
\end{corollary}

\begin{proof}
When $n = 2$, we've shown in the above lemma that the coefficient of $(g_{\delta})_{xx}A^2$ is 1. When $n \geq 2$, if $p=0$ there is no $I$ since $|I| = 0$ and in (\ref{coefficients of k}) we have only one term $\frac{\partial^n g_{\delta}}{\partial x^n}A^n$ whose coefficient is 1 by induction. Taking its derivative with respect to $u$ yields
\[\frac{\partial^{n+1} g_{\delta}}{\partial x^{n+1}}A^{n+1} + \frac{\partial^n g_{\delta}}{\partial x^n}nA^{n-1}\frac{\partial A}{\partial u},\]
so the coefficient of $\frac{\partial^{n+1} g_{\delta}}{\partial x^{n+1}}A^{n+1}$ is still 1 completing the induction.
\end{proof}

Since $A = \cos \beta + \sin \beta (k_\delta)_u(u, v)$ and $(k_\delta)_u(0, 0) = 0$, $A(0, 0) = \cos \beta$ and for $1 \leq p \leq N-2$
\[\frac{\partial^p A}{\partial u^p}(0, 0) = \sin \beta \frac{\partial^{p+1} k_\delta}{\partial u^{p+1}}(0, 0) = \sin \beta (p+1)!b_{p+1}(\delta).\] It follows that
\begin{eqnarray*}
\cos \beta n!b_n(\delta) &=& \sum_{p = 0}^{n-1} \frac{\partial^{n-p} g_{\delta}}{\partial x^{n-p}}(x_0, 0) \sum_{I} c_I (\cos \beta)^{n-p-|I|} (\sin \beta)^{|I|} b_2(\delta)^{i_1} b_3(\delta)^{i_2}\cdots b_{p+1}(\delta)^{i_p}\\
&&2!^{i_1}3!^{i_2}\cdots (p+1)!^{i_p},
\end{eqnarray*}
for $2 \leq n \leq N-1$. Furthermore the Corollary (\ref{cor1}) says that the term corresponding to $p=n-1$ in the above expression is
\[(g_{\delta})_x(x_0, 0)\sin \beta n!b_n(\delta) = -\tan \beta \sin \beta  n!b_n(\delta).\] So moving it to the other side yields
\begin{eqnarray*}
&&(\cos \beta + \tan \beta \sin \beta) n!b_n(\delta) = \\
&&\sum_{p = 0}^{n-2} \frac{\partial^{n-p} g_{\delta}}{\partial x^{n-p}}(x_0, 0) \sum_{I} c_I (\cos \beta)^{n-p-|I|} (\sin \beta)^{|I|} b_2(\delta)^{i_1} b_3(\delta)^{i_2}\cdots b_{p+1}(\delta)^{i_p}\\
&&2!^{i_1}3!^{i_2}\cdots (p+1)!^{i_p}\\
&\Rightarrow& n!b_n(\delta) = \sum_{p = 0}^{n-2} \frac{\partial^{n-p} g_{\delta}}{\partial x^{n-p}}(x_0, 0) \sum_{I} c_I (\cos \beta)^{n+1-p-|I|} (\sin \beta)^{|I|} b_2(\delta)^{i_1} b_3(\delta)^{i_2}\cdots b_{p+1}(\delta)^{i_p}\\
&&2!^{i_1}3!^{i_2}\cdots (p+1)!^{i_p},
\end{eqnarray*}
using $\cos \beta + \tan \beta \sin \beta = \sec \beta$. So $b_n(\delta)$ depends on the previous constants for $3 \leq n \leq N-1$. 

\begin{lemma}\label{lem3.1}
$b_n(\delta) > 0$ for $n$ between 2 and $N-1$.
\end{lemma}

\begin{proof}
Before proceeding with the proof, we need to first estimate $\frac{\partial^{p} g_{\delta}}{\partial x^{p}}(x_0, 0)$ for $2 \leq p \leq N-1$. By induction one can show that
\begin{eqnarray*}
\frac{\partial^{p} g_{\delta}}{\partial x^{p}}(x_0, 0) &=& p!a_p(\delta) + \cdots + (N-1)(N-2)\cdots(N-p)a_{N-1}(\delta)x_0^{N-1-p} \\
&&+ \sum_{q=0}^p \left(\begin{matrix}p \\ q\end{matrix}\right)N(N-1)\cdots(N-q+1)x_0^{N-q} \ \partial_x^{p-q}h_{\delta}(x_0, 0).
\end{eqnarray*}
On the one hand, if $\epsilon$ is sufficiently small 
\begin{eqnarray*}
&&|p!a_p(\delta) + \cdots + (N-1)(N-2)\cdots(N-p)a_{N-1}(\delta)x_0^{N-1-p}| \\
&\leq& M\sin \delta[p! + \cdots + (N-1)(N-2)\cdots(N-p)x_0^{N-1-p}]\\
&\leq&c\frac{1}{4}b\tan(\theta_0) \sin \delta \leq c|a_2(\delta)|. 
\end{eqnarray*}
On the other hand, if $\eta$ and $\epsilon$ are sufficiently small
\begin{eqnarray*}
&&\sum_{q=0}^p \left(\begin{matrix}p \\ q\end{matrix}\right)N(N-1)\cdots(N-q+1)x_0^{N-q} \ \partial_x^{p-q}h_{\delta}(x_0, 0)\\
&\geq& N(N-1)\cdots(N-p+1)x_0^{N-p}h(0, 0)(1-c).
\end{eqnarray*}
Combining the two inequalities, together with (\ref{x_0N-2}), yields
\begin{eqnarray*}
\frac{\partial^{p} g_{\delta}}{\partial x^{p}}(x_0, 0) &\geq& N(N-1)\cdots(N-p+1)x_0^{N-p}h(0, 0)(1-c) - c|a_2(\delta)| \\
&\geq& \frac{N(N-1)(1-c)^2h(0, 0)|a_2(\delta)|}{(N-1)(1+c)h(0, 0)} - c|a_2(\delta)|= \left[\frac{N(1-c)^2}{1+c}-c\right]|a_2(\delta)|,
\end{eqnarray*}
which is positive if we choose $c$ as follows. 
\[N > \frac{c(1+c)}{(1-c)^2} \Rightarrow 2 > \frac{c(1+c)}{(1-c)^2} \Rightarrow 0 < c < \frac{5-\sqrt{17}}{2} < 1.\]
Denote the constant in the brackets as $L = L(c, N)$,  then for $2 \leq p \leq N-1$
\[\frac{\partial^{p} g_{\delta}}{\partial x^{p}}(x_0, 0) \geq L|a_2(\delta)|.\]

Let's determine the signs of $b_n(\delta)$ for $2 \leq n \leq N-1$. When $n=2$,
\[2!b_2(\delta) = \cos^3 \beta (g_{\delta})_{xx}(x_0, 0) \geq \cos^3 \beta L|a_2(\delta)| > 0 \Rightarrow b_2(\delta) > 0.\] 
When $n\geq 3$, by induction
\begin{eqnarray*}
n!b_n(\delta) &\geq& \sum_{p=0}^{n-2} L|a_2(\delta)| \sum_{I}c_I(\cos \beta)^{n+1-p-|I|} (\sin \beta)^{|I|} b_2(\delta)^{i_1} b_3(\delta)^{i_2}\cdots b_{p+1}(\delta)^{i_p}\\
&&2!^{i_1}3!^{i_2}\cdots (p+1)!^{i_p} > 0. 
\end{eqnarray*}
Indeed, when $p=0$, the corresponding term in the above sum, together with Corollary \ref{cor2}, is
\[\frac{\partial^n g_\delta}{\partial x^n}(x_0, 0)\cos^{n+1} \beta \geq L|a_2(\delta)| \cos^{n+1} \beta > 0.\]
So $b_n(\delta) > 0$, as desired. 
\end{proof}

The following lemma shows that the sign of $l_{\delta}(0, 0)$ is also positive. Furthermore, it gives a lower bound of $l_{\delta}(0, 0)$.
\begin{lemma}\label{lem3.2}
\[l_{\delta}(0, 0) \geq \frac{1-c}{2^{N+1}}h(0, 0) > 0.\]
\end{lemma}

\begin{proof}
By Lemma \ref{lem3},
\[\cos \beta \frac{\partial^N k_\delta}{\partial u^N}(u, v) = \sum_{p = 0}^{N-1} \frac{\partial^{N-p} g_{\delta}}{\partial x^{N-p}} \sum_{I} c_I A^{N-p-|I|} (\frac{\partial A}{\partial u})^{i_1} (\frac{\partial^2 A}{\partial u^2})^{i_2}\cdots (\frac{\partial^p A}{\partial u^p})^{i_p}.\]
When $p = N-1$, Corollary \ref{cor1} suggests that we have
\[(g_\delta)_x \frac{\partial^{N-1} A}{\partial x^{N-1}}(u, v).\]
Evaluating at $(u, v) = (0, 0)$ gives us
\begin{eqnarray*}
&&\cos \beta N!l_{\delta}(0, 0) - (g_\delta)_x(x_0, y_0) \sin \beta N! l_{\delta}(0, 0) \\
&=& \sum_{p = 0}^{N-2} \frac{\partial^{N-p} g_{\delta}}{\partial x^{N-p}}(x_0, 0) \sum_{I} c_I (\cos \beta)^{N-p-|I|} (\sin \beta)^{|I|} b_2(\delta)^{i_1} b_3(\delta)^{i_2}\cdots b_{p+1}(\delta)^{i_p}\\
&&2!^{i_1}3!^{i_2}\cdots (p+1)!^{i_p}.
\end{eqnarray*}
Since $(g_\delta)_x(x_0, y_0) = -\tan \beta$, 
\begin{eqnarray*}
N!l_{\delta}(0, 0) &=& \sum_{p = 0}^{N-2} \frac{\partial^{N-p} g_{\delta}}{\partial x^{N-p}}(x_0, 0) \sum_{I} c_I (\cos \beta)^{N+1-p-|I|} (\sin \beta)^{|I|} b_2(\delta)^{i_1} b_3(\delta)^{i_2}\cdots b_{p+1}(\delta)^{i_p}\\
&&2!^{i_1}3!^{i_2}\cdots (p+1)!^{i_p}.
\end{eqnarray*}
When $p=0$, the corresponding term in the above summation by Corollary \ref{cor2} is
\[\frac{\partial^{N} g_{\delta}}{\partial x^{N}}(x_0, 0) \cos^{N+1} \beta.\]
Thus
\begin{eqnarray*}
N!l_{\delta}(0, 0) &=& \frac{\partial^{N} g_{\delta}}{\partial x^{N}}(x_0, 0) \cos^{N+1} \beta \\
&&+ \sum_{p = 1}^{N-2} \frac{\partial^{N-p} g_{\delta}}{\partial x^{N-p}}(x_0, 0) \sum_{I} c_I (\cos \beta)^{N+1-p-|I|} (\sin \beta)^{|I|} b_2(\delta)^{i_1} b_3(\delta)^{i_2}\cdots b_{p+1}(\delta)^{i_p}\\
&&2!^{i_1}3!^{i_2}\cdots (p+1)!^{i_p}.
\end{eqnarray*}
where the second term is positive by Lemma \ref{lem3.1}. Moreover, if $\eta$ and $\epsilon$ are sufficiently small
\begin{eqnarray*}
\frac{\partial^{N} g_{\delta}}{\partial x^{N}}(x_0, 0) &=& \sum_{q=0}^N \left(\begin{matrix}N \\ q\end{matrix}\right)N(N-1)\cdots(N-q+1)x_0^{N-q} \ \partial_x^{q}h_{\delta}(x_0, 0)\\
&\geq& N!h(0, 0)(1-c) > 0.
\end{eqnarray*}
Therefore 
\[N!l_{\delta}(0, 0) \geq N!h(0, 0)(1-c)\cos^{N+1} \beta \Rightarrow l_\delta(0, 0) \geq h(0, 0)(1-c) \cos^{N+1} \beta.\]
Since $-\tan \beta = (g_\delta)_x(x_0, 0)$ and $g_x(0, 0) = 0$, it follows that
\[\beta \to 0, \text{ as } \delta, x_0 \to 0.\]
Therefore for $\eta$ and $\epsilon$ sufficiently small, one can have
\[\cos \beta \geq \frac{1}{2}.\]
So
\[l_\delta(0, 0) \geq \frac{1-c}{2^{N+1}}h(0, 0) > 0.\]
\end{proof}

4. Approximate $v(s)$ and $v'(s)$ using the normal vector $N(s)$ to $S$. We denote $\gamma(s) = (u(s), v(s), w(s))$. If $\gamma(s) \in S$, the normal vector to $S$ at $\gamma(s)$ is 
\[N(s) = (-(k_\delta)_u(v(s), v(s), -(k_\delta)_v(u(s), v(s), 1).\]
Since $u'(s) \geq \frac{1}{2}$, $u(s)$ has a $C^1$-inverse function $s(u)$. Therefore we can express $v(s)$ as
\[v(s) = v(s(u)) = \alpha(u),\]
where $\alpha$ is a $C^1$-function and $\alpha(0) = \frac{d\alpha}{du}(0) = 0$. Then one has $v(s)=o(u(s))$. Hence
\[(k_\delta)_u(u(s), v(s)) = u(s)V_1(s); (k_\delta)_v(u(s), v(s)) = u(s)V_2(s), \] where $V_1(s), V_2(s)$ are bounded by some constant $A$. Let $\gamma''(s)$ exist, then $\gamma''(s) = w''(s)N(s)$, so 
\[u''(s) = -w''(s)u(s)V_1(s), v''(s) = -w''(s)u(s)V_2(s).\]
When $\gamma(s)$ does not touch the surface, the equalities still hold since $\gamma''(s) = 0$. With $w(0) = w'(0) = v(0) = v'(0) = 0$ one can approximate
\[|v'(s)| \leq \int_0^s |v''(\sigma)| d\sigma \leq Au(s)w'(s) \Rightarrow |v(s)| \leq Au(s)w(s).\]
If $\gamma(s) \in S$, then 
\begin{eqnarray*}
w(s) &=& k_\delta(u(s), v(s)) \\
&=&b_2(\delta)u(s)^2 + \cdots + b_{N-1}(\delta)u(s)^{N-1} + u(s)^Nl_{\delta} + v(s)[u(s)m_{\delta} + v(s)n_{\delta}]\\
&\leq& b_2(\delta)u(s)^2 + \cdots + b_{N-1}(\delta)u(s)^{N-1} + u(s)^N|l_{\delta}| + Au(s)w(s)C_1,
\end{eqnarray*}
where $C_1 \to 0$ as $\epsilon \to 0$. Furthermore one can choose $\epsilon$ so small that 
\[|l_{\delta}(u(s), v(s)) - l_{\delta}(0, 0)| \leq ch(0, 0)(1-c) \leq cl_{\delta}(0, 0) \Rightarrow |l_{\delta}(u(s), v(s))| \leq (1+c)l_{\delta}(0, 0), \] by uniform continuity. Therefore there is a constant $B$ such that
\[w(s) \leq B[b_2(\delta)u(s)^2 + \cdots + b_{N-1}(\delta)u(s)^{N-1} + u(s)^N(1+c)l_{\delta}(0, 0)].\] So
\[v(s) \leq Au(s)B[b_2(\delta)u(s)^2 + \cdots + b_{N-1}(\delta)u(s)^{N-1} + u(s)^N(1+c)l_{\delta}(0, 0)].\]

Next let's pproximate $v'(s)$. If $\gamma(s) \in S$, differentiating $w(s)= k_\delta(u(s), v(s))$ gives
\begin{eqnarray*}
w'(s) &=& 2b_2(\delta)u(s)u'(s) + \cdots + (N-1)b_{N-1}(\delta)u(s)^{N-2}u'(s) \\
&& + Nu(s)^{N-1}u'(s)l_{\delta} + u(s)^N[(l_{\delta})_uu'(s) + (l_{\delta})_vv'(s)] \\
&&+ u'(s)v(s)m_{\delta} + u(s)v'(s)m_{\delta} + u(s)v(s)[(m_{\delta})_uu'(s) + (m_{\delta})_vv'(s)]\\
&& + 2v(s)v'(s)n_{\delta} + v(s)^2[(n_{\delta})_uu'(s) + (n_{\delta})_vv'(s)]\\
&\leq&2b_2(\delta)u(s) + \cdots + (N-1)b_{N-1}(\delta)u(s)^{N-2} + u(s)^{N-1}|Nl_{\delta} +u(s)[(l_{\delta})_u + (l_{\delta})_v]| \\
&& + |v(s)|C_1 + u(s)|v'(s)|C_2 + u(s)|v(s)|C_3 + 2|v(s)||v'(s)|C_4 + v(s)^2C_5\\
&\leq&2b_2(\delta)u(s) + \cdots + (N-1)b_{N-1}(\delta)u(s)^{N-2}+u(s)^{N-1}N(1+2c)l_{\delta}(0, 0) \\
&&+(u(s)C_2 + 2|v(s)|C_4)Au(s)w'(s) + (C_1 + u(s)C_3 + |v(s)|C_5)Au(s)B \cdot \\
&&[b_2(\delta)u(s)^2 + \cdots + b_{N-1}(\delta)u(s)^{N-1} + u(s)^N(1+c)l_{\delta}(0, 0)],\\
&\leq& 3b_2(\delta)u(s) + \cdots + Nb_{N-1}(\delta)u(s)^{N-2} + u(s)^{N-1}N(2+2c)l_{\delta}(0, 0)\\
&&+(u(s)C_2 + 2|v(s)|C_4)Au(s)w'(s),
\end{eqnarray*}
where one can choose $\epsilon$ so small that 
\begin{eqnarray*}
u(s)|(l_{\delta})_u + (l_{\delta})_v| \leq Nc\frac{1-c}{2^{N+1}}h(0, 0) \leq Ncl_{\delta}(0, 0). \\
(C_1 + u(s)C_3 + |v(s)|C_5)Au(s)B \leq 1, u(s)\leq 1, 1+c < N
\end{eqnarray*}
By making $(u(s)C_2 + 2|v(s)|C_4)Au(s) < 1$, there exists a constant $C$ such that 
\[w'(s) \leq C[3b_2(\delta)u(s) + \cdots + Nb_{N-1}(\delta)u(s)^{N-2} + u(s)^{N-1}N(2+2c)l_{\delta}(0, 0)].\]
So
\[|v'(s)|\leq Au(s)C[3b_2(\delta)u(s) + \cdots + Nb_{N-1}(\delta)u(s)^{N-2} + u(s)^{N-1}N(2+2c)l_{\delta}(0, 0)].\]

Now let's look at the situation when $\gamma(s)$ is on an interior line segment. Consider a line segment in the image of $\gamma$ with two endpoints $\gamma(s_1)$ and $\gamma(s_2)$, we can parametrize $v(s)$ for $s \in [s_1, s_2]$ by
\[v(s) = v(s_1) + T(u(s) - u(s_1)), \text{ where } T = \frac{d\alpha}{du}(u(s_1)) \text{ and } |T| = \left| \frac{v'(s_1)}{u'(s_1)} \right| \leq 2|v'(s_1)|.\]
\begin{eqnarray*}
|T| &\leq&2Au(s_1)C[3b_2(\delta)u(s_1) + \cdots + Nb_{N-1}(\delta)u(s_1)^{N-2} + u(s_1)^{N-1}N(2+2c)l_{\delta}(0, 0)] \\
&\leq&2Au(s)C[3b_2(\delta)u(s) + \cdots + Nb_{N-1}(\delta)u(s)^{N-2} + u(s)^{N-1}N(2+2c)l_{\delta}(0, 0)],
\end{eqnarray*}
where the last inequality holds because $u(s)$ is increasing. Hence
\begin{eqnarray*}
|v(s)| &\leq& |v(s_1)| + |T|(u(s) + u(s)) \\
&\leq& Au(s)B[b_2(\delta)u(s)^2 + \cdots + b_{N-1}(\delta)u(s)^{N-1} + u(s)^N(1+c)l_{\delta}(0, 0)] \\
&&+4Au(s)^2C[3b_2(\delta)u(s) + \cdots + Nb_{N-1}(\delta)u(s)^{N-2} + u(s)^{N-1}N(2+2c)l_{\delta}(0, 0)]\\
&\leq&(AB + 4AC\cdot2N)u(s)[b_2(\delta)u(s)^2 + \cdots + b_{N-1}(\delta)u(s)^{N-1} + u(s)^N(1+c)l_{\delta}(0, 0)]\\
&=&Du(s)[b_2(\delta)u(s)^2 + \cdots + b_{N-1}(\delta)u(s)^{N-1} + u(s)^N(1+c)l_{\delta}(0, 0)], 
\end{eqnarray*}
where $D = (AB + 4AC\cdot2N)$. Furthermore, since $v'(s)= v'(s_1)$, 
\begin{eqnarray*}
|v'(s)|&\leq& Au(s)C[3b_2(\delta)u(s) + \cdots + Nb_{N-1}(\delta)u(s)^{N-2} + u(s)^{N-1}N(2+2c)l_{\delta}(0, 0)]\\
&\leq&AC \cdot 2Nu(s)[b_2(\delta)u(s) + \cdots + b_{N-1}(\delta)u(s)^{N-2} + u(s)^{N-1}(1+c)l_{\delta}(0, 0)]\\
&\leq&Du(s)[b_2(\delta)u(s) + \cdots + b_{N-1}(\delta)u(s)^{N-2} + u(s)^{N-1}(1+c)l_{\delta}(0, 0)].
\end{eqnarray*}

5. Concavity. Suppose for the sake of contradiction that $\gamma(s)$ leaves $S$ at a switch point when $s=s_0$ and dives into the interior of $M$ for increasing $s$ until it re-enters $S$ again at $s = s_1$.

Consider the intersection of the two-dimensional plane $v = v_0 + T(u - u_0)$ with the surface $w = k_\delta(u, v)$, where $(u_0, v_0) = (u(s_0), v(s_0))$ and $T = \frac{d\alpha}{du}(u_0)$. 
Set $f(u) = k_\delta(v, v_0 + T(u-u_0))$, then with $v(s) = v_0 + T(u(s)-u_0)$,
\[\frac{d^2f}{du^2}(u(s)) = (k_\delta)_{uu}(u(s), v(s)) + 2(k_\delta)_{uv}(u(s), v(s))T + (k_\delta)_{vv}(u(s), v(s))T^2.\]
On the one hand, 
\begin{eqnarray*}
(k_\delta)_{uu} &=& 2b_2(\delta) + 6b_3(\delta)u(s) + \cdots + (N-1)(N-2)b_{N-1}(\delta)u(s)^{N-3} \\
&& + N(N-1)u(s)^{N-2}l_{\delta} + 2Nu(s)^{N-1}(l_{\delta})_u + u(s)^N(l_{\delta})_{uu} \\
&& + v(s)[2(m_{\delta})_u + u(s)(m_{\delta})_{uu}] + v(s)^2(n_{\delta})_{uu},
\end{eqnarray*}
where one can choose $\epsilon$ small enough so that
\[|l_{\delta}| \geq (1-c)l_{\delta}(0, 0), \ \frac{|2Nu(s)(l_{\delta})_u + u(s)^2(l_{\delta})_{uu}|}{N(N-1)} \leq cl_{\delta}(0, 0), 
|2(m_{\delta})_u + u(s)(m_{\delta})_{uu} + v(s)(n_{\delta})_{uu}| \leq C_1.\]
Thus
\begin{eqnarray*}
(k_\delta)_{uu} &\geq& 2b_2(\delta) + 6b_3(\delta)u(s) + \cdots + (N-1)(N-2)b_{N-1}(\delta)u(s)^{N-3} \\
&&+ N(N-1)u(s)^{N-2}(1-2c) l_{\delta}(0, 0) \\
&&- C_1Du(s)[b_2(\delta)u(s)^2 + \cdots + b_{N-1}(\delta)u(s)^{N-1} + u(s)^N(1+c)l_{\delta}(0, 0)]\\
&\geq& 2b_2(\delta) + 6b_3(\delta)u(s) + \cdots + (N-1)(N-2)b_{N-1}(\delta)u(s)^{N-3} \\
&&+ N(N-1)u(s)^{N-2}(1-2c) l_{\delta}(0, 0) \\
&&- \frac{1}{2}[b_2(\delta) + \cdots + b_{N-1}(\delta)u(s)^{N-3} + u(s)^{N-2}(1+c)l_{\delta}(0, 0)],\\
\end{eqnarray*}
where we can choose $\epsilon$ small enough so that $C_1Du(s) \leq \frac{1}{2}$ and $u(s) \leq 1$. Note here that $1-2c > 0$ because
\[0 < c < \frac{5-\sqrt{17}}{2} < \frac{1}{2}.\]
On the other hand, 
\begin{eqnarray*}
|T| &\leq&2Au(s_0)C[3b_2(\delta)u(s_0) + \cdots + Nb_{N-1}(\delta)u(s_0)^{N-2} + u(s_0)^{N-1}N(2+2c)l_{\delta}(0, 0)] \\
&\leq&4ACNu(s_0)[b_2(\delta)u(s_0)^2 + \cdots + b_{N-1}(\delta)u(s_0)^{N-1} + u(s_0)^N(1+c)l_{\delta}(0, 0)]\\
&\leq&Du(s)[b_2(\delta)u(s)^2 + \cdots + b_{N-1}(\delta)u(s)^{N-1} + u(s)^N(1+c)l_{\delta}(0, 0)].
\end{eqnarray*}
So for $\epsilon$ sufficiently small
\begin{eqnarray*}
&&|2(k_\delta)_{uv}T+(k_\delta)_{vv}T^2| \leq C_1|T| \\
&\leq& C_1Du(s)[b_2(\delta)u(s)^2 + \cdots + b_{N-1}(\delta)u(s)^{N-1} + u(s)^N(1+c)l_{\delta}(0, 0)]\\
&\leq& \frac{1}{2}[b_2(\delta)+ \cdots + b_{N-1}(\delta)u(s)^{N-3} + u(s)^{N-2}(1+c)l_{\delta}(0, 0)],
\end{eqnarray*}
if $C_1Du(s) \leq \frac{1}{2}$ and $u(s)\leq 1$. 
It follows that
\begin{eqnarray*}
&&\frac{d^2f}{du^2}(u(s)) \\
&\geq& 2b_2(\delta) + 6b_3(\delta)u(s) + \cdots + (N-1)(N-2)b_{N-1}(\delta)u(s)^{N-3} \\
&&+ N(N-1)(1-2c)u(s)^{N-2}l_{\delta}(0, 0) \\
&& - b_2(\delta) -b_3(\delta)u(s) - \cdots - b_{N-1}(\delta)u(s)^{N-3} - u(s)^{N-2}(1+c)l_{\delta}(0, 0)\\
&=&b_2(\delta) + 5b_3(\delta)u(s) +\cdots+[(N-1)(N-2)-1]b_{N-1}(\delta)u(s)^{N-3} \\
&&+ [N(N-1)(1-2c)-(1+c)]u(s)^{N-2}l_{\delta}(0, 0),
\end{eqnarray*}
where we want
\[N(N-1)(1-2c) \geq 2(1-2c) > 1 + c \Rightarrow 0 < c < \frac{1}{5} < \frac{5-\sqrt{17}}{2}.\]
Therefore 
\[\frac{d^2f}{du^2}(u(s)) > 0,\]
so $f'(u(s))$ is increasing as $u(s)$ increases from $u(s_0)=u_0$ to $u(s_1)=u_1$. On the other hand, since the interior line segment is tangent to $S$ at the two endpoints, we must have
\[f'(u_0) = f'(u_1), \]
a contradiction. Therefore if $\gamma$ leaves $S$ at a switch point $\gamma(s_0)$, the geodesic arc beyond this point is a line segment either exiting the $\epsilon$-ball or terminating at a point on $S$. So $\gamma$ has at most two switch points in this case. 

{\bf Case 4: $a_2(\delta) > 0$, $h(0, 0) < 0$.}

This is a combination of Case 1, Case 2, and Case 3. 
When $a_2(\delta) > 0$, the angle $\delta < 0$. From Case 1 the constants $a_2(\delta)$, $a_3(\delta)$, $\ldots$, $a_{N-1}(\delta)$ satisfy the following: 
\[a_2(\delta) \geq b\tan(\theta_0)|\sin \theta|; |a_3(\delta)|, \ldots, |a_{N-1}(\delta)| \leq M|\sin \delta|. \]
Assume $\gamma'(0) = \frac{\partial }{\partial x}$. 

1. $x'(s) \geq \frac{1}{2}$. The argument is exactly the same as in Case 1 when we proved for (\ref{x'1/2}). 

2. Approximate $z(s), z'(s), y(s), y'(s)$. If $\gamma(s) \in S$, imitating the proof for (\ref{zy}) in Case 1 one can show that there exists a positive constant $C$ such that 
\begin{eqnarray*}
&& z(s) \leq (2a_2(\delta)x(s)^2 + 2|h(0, 0)|x(s)^N)C \\ \notag
&\Rightarrow& |y(s)| \leq (2a_2(\delta)x(s)^2 + 2|h(0, 0)|x(s)^N)AC \notag
\end{eqnarray*}
Simiarly, imitating the proof for (\ref{z'y'}) in Case 1 one can also show that 
\begin{eqnarray*}
&&z'(s) \leq (4a_2(\delta)x(s) + 4Nx(s)^{N-1}|h(0, 0)|)C  \\ \notag
&\Rightarrow& y'(s) \leq (4a_2(\delta)x(s) + 4Nx(s)^{N-1}|h(0, 0)|)AC. \notag
\end{eqnarray*}
It follows that if $y = \alpha(x)$ and $T(s) = \alpha'(x(s))$, then
\[|T(s)| = \left|\frac{y'(s)}{x'(s)} \right| \leq 2|y'(s)| \leq (8a_2(\delta)x(s) + 8Nx(s)^{N-1}|h(0, 0)|)AC.\]
If $\gamma(s) \not \in S$, then $\gamma(s)$ is in some line segment where $y(s) = y(s_0) + T(s_0)(x(s) - x(s_0))$ for some switch point at $s_0 < s$. Since $x(s)$ is increasing, it follows that
\begin{eqnarray*}
|y(s)| &\leq& |y(s_0)| + |T(s_0)| [|x(s)|+ |x(s_0)|] \\
&\leq&  (2a_2(\delta)x(s_0)^2 + 2|h(0, 0)|x(s_0)^N)AC \\
&&+  (8a_2(\delta)x(s_0) + 8Nx(s_0)^{N-1}|h(0, 0)|)AC [|x(s)|+ |x(s_0)|] \\
&\leq& (2a_2(\delta)x(s)^2 + 2|h(0, 0)|x(s)^N)AC + (8a_2(\delta)x(s) + 8Nx(s)^{N-1}|h(0, 0)|)AC \cdot 2x(s)\\
&\leq& (18a_2(\delta)x(s)^2 + 18N|h(0,0)|x(s)^N)AC.
\end{eqnarray*}
Furthermore, 
\begin{eqnarray*}
|z(s)| &\leq& |z(s_0)| + 2|z'(s_0)|[|x(s)|+ |x(s_0)|] \\
&\leq&  (2a_2(\delta)x(s_0)^2 + 2|h(0, 0)|x(s_0)^N)C + 2(4a_2(\delta)x(s_0) + 4Nx(s_0)^{N-1}|h(0, 0)|)C \cdot 2x(s)\\
&\leq& (2a_2(\delta)x(s)^2 + 2|h(0, 0)|x(s)^N)C +4x(s)(4a_2(\delta)x(s) + 4Nx(s)^{N-1}|h(0, 0)|)C\\
&=&(18a_2(\delta)x(s)^2 + 18N|h(0, 0)|x(s)^N)C.
\end{eqnarray*}
The same inequalities for $z'(s)$, $y'(s)$, and $T(s)$ still hold as before because 
\[z'(s) = z'(s_0), y'(s) = y'(s_0), T(s) = T(s_0).\]

3. Now we are ready to estimate the location of the first switch point. Set $f(x) = g_{\delta}(x, y_0 + T(x-x_0))$ where $(x_0, y_0) = (x(s_0), y(s_0))$ and $T = T(s_0)$. Then 
\begin{eqnarray*}
f''(x_0) &=& (g_{\delta})_{xx}(x_0, y_0) + 2(g_{\delta})_{xy}(x_0, y_0)T + (g_{\delta})_{yy}(x_0, y_0)T^2,
\end{eqnarray*}
where
\begin{eqnarray*}
(g_{\delta})_{xx}(x_0, y_0) &=& 2a_2(\delta) + 6a_3(\delta)x_0 + \cdots + (N-1)(N-2)a_{N-1}(\delta)x_0^{N-3} \\
&& + x_0^{N-2}[N(N-1)h_{\delta} + 2Nx_0(h_{\delta})_x + x_0^2(h_{\delta})_{xx}] \\
&& + y_0[2(i_{\delta})_x + x_0(i_{\delta})_{xx} + y_0(j_{\delta})_{xx}].
\end{eqnarray*}
As $x_0 \to 0$, $(g_{\delta})_{xx}(x_0, y_0) \to 2a_2(\delta)$ and $T(s_0) \to 0$. Therefore $f''(x_0) > 0$ at the beginning and the geodesic initially stays on the surface if it is not a straight line for which we will discuss later. 

Suppose $f''(x_0) = 0$. Let's estimate $x_0$. 
First, for $\epsilon$ sufficiently small,
\begin{eqnarray*}
&&|6a_3(\delta)x_0 + \cdots + (N-1)(N-2)a_{N-1}(\delta)x_0^{N-3}| \\
&\leq& M|\sin \delta| (6x_0 + \cdots + (N-1)(N-2)x_0^{N-3}) \\
&\leq& b\tan(\theta_0)|\sin \delta| \leq a_2(\delta).
\end{eqnarray*}
Second, for $\eta$ and $\epsilon$ sufficiently small, 
\[N(N-1)h_{\delta} + 2Nx_0(h_{\delta})_x + x_0^2(h_{\delta})_{xx} \geq (1+\frac{c}{2})N(N-1)h(0, 0),\]
for some $0 < c < 1$ to be determined later. 
Third, there are constants $C_5, C_6$ such that 
\[|2(i_{\delta})_x + x_0(i_{\delta})_{xx} + y_0(j_{\delta})_{xx}| \leq C_5, |2(g_{\delta})_{xy}(x_0, y_0)+(g_{\delta})_{yy}(x_0, y_0)T| \leq C_6.\]
If $f''(x_0) = 0$, then 
\begin{eqnarray*}
&&-x_0^{N-2}[N(N-1)h_{\delta} + 2Nx_0(h_{\delta})_x + x_0^2(h_{\delta})_{xx}] \\
&=& 2a_2(\delta) + 6a_3(\delta)x_0 + \cdots + (N-1)(N-2)a_{N-1}(\delta)x_0^{N-3} \\
&& + y_0[2(i_{\delta})_x + x_0(i_{\delta})_{xx} + y_0(j_{\delta})_{xx}] + 2(g_{\delta})_{xy}(x_0, y_0)T + (g_{\delta})_{yy}(x_0, y_0)T^2. 
\end{eqnarray*}
On the one hand, 
\[\text{LHS} \leq -(1+\frac{c}{2})N(N-1)h(0, 0)x_0^{N-2} =(1+\frac{c}{2})N(N-1)|h(0, 0)|x_0^{N-2} .\]
On the other hand, 
\begin{eqnarray*}
\text{RHS} &\geq& 2a_2(\delta) - a_2(\delta) - C_5(18a_2(\delta)x_0^2 + 18N|h(0, 0)|x_0^N)AC \\
&& - C_6(8a_2(\delta)x_0 + 8Nx_0^{N-1}|h(0, 0)|)AC.
\end{eqnarray*}
Combining the two inequalities we get
\begin{eqnarray*}
&&N(N-1)|h(0, 0)|x_0^{N-2} [(1+\frac{c}{2}) + \frac{18ACC_5x_0^2 + 8ACC_6x_0}{N-1}] \\
&\geq&a_2(\delta)[1-18ACC_5x_0^2 - 8ACC_6x_0]
\end{eqnarray*}
Choose $\epsilon$ small enough so that
\[\frac{18ACC_5x_0^2 + 8ACC_6x_0}{N-1} \leq \frac{c}{2} \text{ and } 18ACC_5x_0^2 + 8ACC_6x_0 \leq c.\]
Therefore we obtain a lower bound for $x_0$ 
\begin{eqnarray}\label{x_0N-22}
&&N(N-1)|h(0, 0)|x_0^{N-2} (1+c) \geq a_2(\delta)(1-c) \\ \notag
&\Rightarrow& x_0^{N-2} \geq \frac{a_2(\delta)(1-c)}{N(N-1)(1+c)|h(0, 0)|}. 
\end{eqnarray}

Similarly, we can also get an upper bound for $x_0$. Again suppose $f''(x_0)=0$. First, for $\epsilon$ sufficiently small, we still have 
\[|6a_3(\delta)x_0 + \cdots + (N-1)(N-2)a_{N-1}(\delta)x_0^{N-3}| \leq a_2(\delta).\]
Second, for $\eta$ and $\epsilon$ sufficiently small, 
\[N(N-1)h_{\delta} + 2Nx_0(h_{\delta})_x + x_0^2(h_{\delta})_{xx} \leq (1-\frac{c}{2})N(N-1)h(0, 0).\]
Third, there are still constants $C_5$, $C_6$ such that
\[|2(i_{\delta})_x + x_0(i_{\delta})_{xx} + y_0(j_{\delta})_{xx}| \leq C_5, |2(g_{\delta})_{xy}(x_0, y_0)+(g_{\delta})_{yy}(x_0, y_0)T| \leq C_6.\]
If $f''(x_0) = 0$, then 
\begin{eqnarray*}
&&-x_0^{N-2}[N(N-1)h_{\delta} + 2Nx_0(h_{\delta})_x + x_0^2(h_{\delta})_{xx}] \\
&=& 2a_2(\delta) + 6a_3(\delta)x_0 + \cdots + (N-1)(N-2)a_{N-1}(\delta)x_0^{N-3} \\
&& + y_0[2(i_{\delta})_x + x_0(i_{\delta})_{xx} + y_0(j_{\delta})_{xx}] + 2(g_{\delta})_{xy}(x_0, y_0)T + (g_{\delta})_{yy}(x_0, y_0)T^2. 
\end{eqnarray*}
On the one hand, 
\[\text{ LHS } \geq -(1-\frac{c}{2})N(N-1)h(0, 0)x_0^{N-2} = (1-\frac{c}{2})N(N-1)|h(0, 0)|x_0^{N-2}.\]
On the other hand, 
\begin{eqnarray*}
\text{RHS} &\leq& 2a_2(\delta) + a_2(\delta) + C_5(18a_2(\delta)x_0^2 + 18N|h(0, 0)|x_0^N)AC \\
&& + C_6(8a_2(\delta)x_0 + 8Nx_0^{N-1}|h(0, 0)|)AC.
\end{eqnarray*}
Combining the two yields
\begin{eqnarray*}
&&(1-\frac{c}{2})N(N-1)|h(0, 0)|x_0^{N-2} \leq 3a_2(\delta) + (18ACC_5x_0^2 + 8ACC_6x_0)a_2(\delta) \\
&& + x_0^{N-2}N|h(0, 0)|(18ACC_5x_0^2 + 8ACC_6x_0)\\
&\Rightarrow&N(N-1)|h(0, 0)|x_0^{N-2}[(1-\frac{c}{2}) - \frac{18ACC_5x_0^2 + 8ACC_6x_0}{N-1}] \\
&\leq& 3a_2(\delta) + (18ACC_5x_0^2 + 8ACC_6x_0)a_2(\delta).
\end{eqnarray*}
So 
\[N(N-1)|h(0, 0)|x_0^{N-2}(1-c) \leq (3+c) a_2(\delta) \Rightarrow x_0^{N-2} \leq \frac{(3+c) a_2(\delta)}{N(N-1)|h(0, 0)|(1-c)}.\]

4. Now like in Case 3 we are going to shift our coordinates to have the origin at $\gamma(s_0) = (x_0, y_0, z_0)$ and then rotate the $(x, y, z)$-space so that $\gamma'(s_0) = (x'(s_0), y'(s_0), z'(s_0))$ points in the positive $x$-axis. Let's use $(u, v, w)$ for the new coordinates, then with respect to the new frame there is a rotation matrix $P \in SO_3(\mathbb{R})$ such that
\begin{eqnarray*}
&&\left(\begin{matrix} x-x_0 \\ y-y_0 \\ z-z_0\end{matrix}\right) = P\left( \begin{matrix}u \\ v\\ w \end{matrix}\right), \text{ where } 
P = \left( \begin{matrix} x'(s_0) & p_{12} & p_{13} \\ y'(s_0) & p_{22} & p_{23} \\ z'(s_0) & p_{32} & p_{33} \end{matrix}\right).
\end{eqnarray*}
So
\begin{eqnarray*}
x &=& x_0 + x'(s_0)u + p_{12}v + p_{13}w, \\
y &=& y_0 + y'(s_0)u + p_{22}v + p_{23}w, \\
z &=& z_0 + z'(s_0)u + p_{32}v + p_{33}w.
\end{eqnarray*}
The second and third columns of $P$ can be further specificed as below:
\begin{equation}\label{P2P3}
\left(\begin{matrix} p_{12} \\ p_{22} \\ p_{32}\end{matrix}\right) = \left( \begin{matrix} -y'(s_0) \\ x'(s_0) \\ 0\end{matrix}\right)Y; \left( \begin{matrix} p_{13} \\ p_{23} \\ p_{33}\end{matrix}\right) = \left(\begin{matrix} -x'(s_0)^2z'(s_0) \\-x'(s_0)y'(s_0)z'(s_0)  \\ x'(s_0)[x'(s_0)^2 + y'(s_0)^2] \end{matrix} \right)Z,
\end{equation}
where $Y = 1/\sqrt{x'(s_0)^2 + y'(s_0)^2}$ and $Z = 1/x'(s_0)\sqrt{x'(s_0)^2 + y'(s_0)^2}$. 

Thus the surface $z=g_{\delta}(x, y)$ satisfies the equation
\begin{eqnarray*}
 z_0 + z'(s_0)u + p_{32}v + p_{33}w = g_{\delta}(x_0 + x'(s_0)u + p_{12}v + p_{13}w, y_0 + y'(s_0)u + p_{22}v + p_{23}w).
\end{eqnarray*}
Check that we can still solve for $w$ analytically in terms of $u, v$ within the $\epsilon$-ball. Let's take the partial derivate of 
\[- z_0 - z'(s_0)u - p_{32}v - p_{33}w + g_{\delta}(x_0 + x'(s_0)u + p_{12}v + p_{13}w, y_0 + y'(s_0)u + p_{22}v + p_{23}w)\]
with respect to $w$:
\[-p_{33} + (g_{\delta})_xp_{13} + (g_{\delta})_yp_{23}.\] Since $Z > 0$, it is equivalent to show that 
\[\frac{-p_{33} + (g_{\delta})_xp_{13} + (g_{\delta})_yp_{23}}{Z} \neq 0.\]
On the one hand, 
\begin{eqnarray*}
\frac{(g_{\delta})_xp_{13} + (g_{\delta})_yp_{23}}{Z} =-x'(s_0)^2z'(s_0)(g_{\delta})_x -x'(s_0)y'(s_0)z'(s_0)(g_{\delta})_y\leq |(g_{\delta})_x|+|(g_{\delta})_y|. 
\end{eqnarray*}
Since $(g_{\delta})_x(x, y)$, $(g_{\delta})_y(x, y)$ $\to 0$ as $x, y \to 0$, for $\epsilon$ sufficiently small
\[|(g_{\delta})_x|+|(g_{\delta})_y| \leq \frac{1}{16}.\] 
On the other hand, 
\[\frac{p_{33}}{Z} = x'(s_0)[x'(s_0)^2 + y'(s_0)^2] \geq x'(s_0)^3 \geq \frac{1}{8}.\]
so
\[\frac{-p_{33} + (g_{\delta})_xp_{13} + (g_{\delta})_yp_{23}}{Z} \leq -\frac{1}{16} < 0.\]
Therefore there exists a real analytic function $k_\delta$ such that $w = k_\delta(u, v)$ such that $k_\delta(0, 0) = 0$, $(k_\delta)_u(0, 0) = 0$, and $(k_\delta)_v(0, 0)=0$. 

5. Estimate $\gamma(s)$ in the new frame starting from $(x_0, y_0, z_0)$. After replacing $s$ by $s-s_0$, we denote $\gamma(s)$ as $(u(s), v(s), w(s))$. 



(1). Coefficients of $k_{\delta}(u, v)$. Denote $k_{\delta}(u, v)$ as 
\[k_\delta(u, v) = b_2(\delta)u^2 + \cdots + b_{N-1}(\delta)u^{N-1} + u^N l_{\delta}(u, v) + uvm_{\delta}(u, v) + v^2n_{\delta}(u, v),\]
where $b_2(\delta), \ldots, b_{N-1}(\delta)$ are constants and $l_{\delta}$, $m_{\delta}$, $n_{\delta}$ are analytic functions. 
Observe that for $n$ between 2 and $N-1$,
\[b_n(\delta) = \frac{1}{n!}\frac{\partial^n k_\delta}{\partial u^n}(0, 0).\]

Next let's look for $\frac{\partial^n k_\delta}{\partial u^n}(u, v)$ for $n \geq 2$ by induction. 
\begin{lemma}\label{lem4}
Let $A$ be $x'(s_0) + p_{13}(k_\delta)_u$ and $B$ be $y'(s_0) + p_{23}(k_\delta)_u$, then for $n \geq 2$, 
\begin{eqnarray*}
&&p_{33}\frac{\partial^n k_\delta}{\partial u^n}(u, v) = \sum_{a+b = 1}^{n} \frac{\partial^{a+b} g_{\delta}}{\partial x^{a}\partial y^b} \sum_{I, J} c_{I, J} A^{a-|I|}B^{b-|J|} (\frac{\partial A}{\partial u})^{i_1}\cdots (\frac{\partial^p A}{\partial u^p})^{i_p}(\frac{\partial B}{\partial u})^{j_1}\cdots (\frac{\partial^p B}{\partial u^p})^{j_p},
\end{eqnarray*}
where $p = n-(a+b)$, $I = (i_1, i_2, \ldots, i_p)$, $J = (j_1, j_2, \ldots, j_p)$, $(i_1 + 2i_2 + \cdots + pi_p) + (j_1 + 2j_2 + \cdots + pj_p) = p $, $|I| = i_1 + i_2 + \cdots + i_p \leq a$, $|J| = j_1 + j_2 + \cdots + j_p \leq b$, and the partial derivatives of $g_{\delta}$ are evaluated at $(x_0 + x'(s_0)u + p_{12}v + p_{13}k_\delta(u, v), y_0 + y'(s_0)u + p_{22}v + p_{23}k_\delta(u, v))$. 
\end{lemma}

\begin{proof}
When $n=2$, differentiating the following equation 
\[ z_0 + z'(s_0)u + p_{32}v + p_{33}k_\delta(u, v) = g_{\delta}(x_0 + x'(s_0)u + p_{12}v + p_{13}k_\delta(u, v), y_0 + y'(s_0)u + p_{22}v + p_{23}k_\delta(u, v))\]
with respect to $u$ once gives us
\begin{equation}\label{firstderivative}
z'(s_0) + p_{33}(k_\delta)_u= (g_{\delta})_x[x'(s_0) + p_{13}(k_\delta)_u] + (g_{\delta})_y[y'(s_0) + p_{23}(k_\delta)_u].
\end{equation}
Let $A = x'(s_0) + p_{13}(k_\delta)_u$ and $B = y'(s_0) + p_{23}(k_\delta)_u$, then
\[z'(s_0) + p_{33}(k_\delta)_u=  (g_{\delta})_xA +  (g_{\delta})_yB.\]
Taking the partial derivative with respect to $u$ once more gives
\[p_{33}(k_\delta)_{uu} = (g_{\delta})_{xx}A^2 + (g_{\delta})_x\partial_u A + (g_{\delta})_{yy}B^2 + (g_{\delta})_y\partial_uB + 2(g_{\delta})_{xy}AB.\]
When $a+b=2$, $p=0$ and so there are no $I$ and $J$. There are three terms corresponding to: $a=2, b=0$; $a=0, b=2$; $a=1, b=1$, respectively: 
\[(g_{\delta})_{xx}A^2,  (g_{\delta})_{yy}B^2,  (g_{\delta})_{xy}AB.\]
When $a+b=1$, either $a=1, b=0$ with $I=(1), J=0$ or $a=0, b=1$ with $I = 0, J= (1)$. It follows that there are two terms
\[(g_{\delta})_x\partial_u A, (g_{\delta})_y\partial_u B.\]

When $n \geq 2$, by inductive hypothesis we can take the partial derivative of the expression in the lemma with respect to $u$. The left-hand side is $p_{33}\frac{\partial^{n+1} k_\delta}{\partial u^{n+1}}(u, v)$. The right-hand side consists of three parts due to the product rule. 
\begin{eqnarray*}
&(1)&\sum_{a+b = 1}^{n} \frac{\partial^{a+b+1} g_{\delta}}{\partial x^{a+1}\partial y^b} \sum_{I, J} c_{I, J} A^{a+1-|I|}B^{b-|J|} (\frac{\partial A}{\partial u})^{i_1}\cdots (\frac{\partial^p A}{\partial u^p})^{i_p}(\frac{\partial B}{\partial u})^{j_1}\cdots (\frac{\partial^p B}{\partial u^p})^{j_p} \\
&& + \sum_{a+b = 1}^{n} \frac{\partial^{a+b+1} g_{\delta}}{\partial x^{a}\partial y^{b+1}} \sum_{I, J} c_{I, J} A^{a-|I|}B^{b+1-|J|} (\frac{\partial A}{\partial u})^{i_1}\cdots (\frac{\partial^p A}{\partial u^p})^{i_p}(\frac{\partial B}{\partial u})^{j_1}\cdots (\frac{\partial^p B}{\partial u^p})^{j_p},
\end{eqnarray*}
where $a$ becomes $a+1$ in the first term, $b$ becomes $b+1$ in the second term, and $p$ stays the same since $(n+1)-(a+1+b) = p$. 
\begin{eqnarray*}
&(2)&\sum_{a+b = 1}^{n} \frac{\partial^{a+b} g_{\delta}}{\partial x^{a}\partial y^b} \sum_{I, J} c_{I, J}(a-|I|) A^{a-|I|-1}B^{b-|J|} (\frac{\partial A}{\partial u})^{i_1+1}\cdots (\frac{\partial^p A}{\partial u^p})^{i_p}(\frac{\partial B}{\partial u})^{j_1}\cdots (\frac{\partial^p B}{\partial u^p})^{j_p}\\
&& + \sum_{a+b = 1}^{n} \frac{\partial^{a+b} g_{\delta}}{\partial x^{a}\partial y^b} \sum_{I, J} c_{I, J}(b-|J|) A^{a-|I|}B^{b-|J|-1} (\frac{\partial A}{\partial u})^{i_1+1}\cdots (\frac{\partial^p A}{\partial u^p})^{i_p}(\frac{\partial B}{\partial u})^{j_1+1}\cdots (\frac{\partial^p B}{\partial u^p})^{j_p},
\end{eqnarray*}
where $a$, $b$ stay the same, so $p$ becomes $p+1 = (n+1)-(a+b)$. Moreover, $i_1$ becomes $i_1 + 1$ in the first term and $j_1$ becomes $j_1+1$ in the second term. So 
\[(i_1 +1 + 2i_2 + \cdots + pi_p) + (j_1 + 2j_2 + \cdots + pj_p) = (i_1 + 2i_2 + \cdots + pi_p) + (j_1 +1 + 2j_2 + \cdots + pj_p) = p+1.\]
\begin{eqnarray*}
&(3)& \sum_{a+b = 1}^{n} \frac{\partial^{a+b} g_{\delta}}{\partial x^{a}\partial y^b} \sum_{I, J} c_{I, J} A^{a-|I|}B^{b-|J|}\sum_{i_k \neq 0}\cdots i_k(\frac{\partial A}{\partial u})^{i_k -1}(\frac{\partial A}{\partial u})^{i_{k+1} +1} \cdots\\
&& + \sum_{a+b = 1}^{n} \frac{\partial^{a+b} g_{\delta}}{\partial x^{a}\partial y^b} \sum_{I, J} c_{I, J} A^{a-|I|}B^{b-|J|}\sum_{j_k \neq 0}\cdots j_k(\frac{\partial B}{\partial u})^{j_k -1}(\frac{\partial B}{\partial u})^{j_{k+1} +1}\cdots,
\end{eqnarray*}
where $a$, $b$ stay the same, so $p$ becomes $p+1=(n+1)-(a+b)$. Moreover, $i_k$, $i_{k+1}$ become $i_k-1$, $i_{k+1}+1$ in the first term and $j_k$, $j_{k+1}$ become $j_k-1$, $j_{k+1}+1$ in the second term. So 
\[\cdots + k(i_k-1)+(k+1)(i_{k+1}+1)+\cdots = \cdots + k(j_k-1)+(k+1)(j_{k+1}+1)+\cdots = p+1.\]
It follows that  $c_{I, J}$ are nonnegative integers and the lemma is true. 
\end{proof}

The following two corollaries are analogous to Corollaries \ref{cor1} and \ref{cor2}.
\begin{corollary}\label{cor3}
The coefficients of $(g_{\delta})_x\partial_u^{n-1}A$ and $(g_{\delta})_y\partial_u^{n-1}B$ are always 1.
\end{corollary}

\begin{corollary}\label{cor4}
The coefficients of $\frac{\partial^{n} g_{\delta}}{\partial x^{n}}A^{n}$ and $\frac{\partial^{n} g_{\delta}}{\partial y^{n}}B^{n}$ are always 1. 
\end{corollary}

Now we let $(u, v) = (0, 0)$, then $(k_\delta)_u(0, 0)=0$ implies that
\[A(0, 0) = x'(s_0), B(0, 0) = y'(s_0).\]
Furthermore for $p \geq 1$,
\begin{eqnarray*}
\frac{\partial^p A}{\partial u^p}(0, 0) = p_{13} \frac{\partial^{p+1} k_\delta}{\partial u^{p+1}}(0, 0) = p_{13}(p+1)!b_{p+1}(\delta)\\
\frac{\partial^p B}{\partial u^p}(0, 0) = p_{23} \frac{\partial^{p+1} k_\delta}{\partial u^{p+1}}(0, 0) = p_{23}(p+1)!b_{p+1}(\delta).
\end{eqnarray*}
It follows that for $2 \leq n \leq N-1$
\begin{eqnarray*}
p_{33}n!b_n(\delta) &=& \sum_{a+b = 1}^{n} \frac{\partial^{a+b} g_{\delta}}{\partial x^{a}\partial y^b}(x_0, y_0) \sum_{I, J} c_{I, J} x'(s_0)^{a-|I|}y'(s_0)^{b-|J|} (p_{13})^{|I|}(p_{23})^{|J|}\\
&&  (2!)^{i_1+j_1}(3!)^{i_2+j_2} \cdots (p+1)!^{i_p + j_p}b_2(\delta)^{i_1+j_1}\cdots b_{p+1}(\delta)^{i_p+j_p}.
\end{eqnarray*}
Furthermore Corollary (\ref{cor3}) suggests that the terms corresponding to $a+b = 1$ or $p=n-1$ in the expression are
\[(g_{\delta})_x(x_0, y_0)p_{13}n!b_n(\delta) + (g_{\delta})_y(x_0, y_0)p_{23}n!b_n(\delta).\]
Moving them to the other side of the expression yields
\begin{eqnarray*}
&&(p_{33}-(g_{\delta})_x(x_0, y_0)p_{13} -(g_{\delta})_y(x_0, y_0)p_{23})n!b_n(\delta) \\
&=&\sum_{a+b = 2}^{n} \frac{\partial^{a+b} g_{\delta}}{\partial x^{a}\partial y^b}(x_0, y_0) \sum_{I, J} c_{I, J} x'(s_0)^{a-|I|}y'(s_0)^{b-|J|} (p_{13})^{|I|}(p_{23})^{|J|}\\
&&  (2!)^{i_1+j_1}(3!)^{i_2+j_2} \cdots (p+1)!^{i_p + j_p}b_2(\delta)^{i_1+j_1}\cdots b_{p+1}(\delta)^{i_p+j_p},
\end{eqnarray*}
where $p+1 = n-(a+b)+1 \leq n-1$.
So $b_n(\delta)$ depends on the previous coefficients $b_2(\delta)$, $\ldots$, $b_{n-1}(\delta)$. 

Let's calculate the coefficient $b_n(\delta)$. Evaluating (\ref{firstderivative}) at $(u, v)=(0, 0)$ arrives
\[z'(s_0) = (g_{\delta})_x(x_0, y_0)x'(s_0) + (g_{\delta})_y(x_0, y_0)y'(s_0),\]
together with (\ref{P2P3}), then
\begin{eqnarray*}
&&\frac{1}{Z}\left[p_{33}-(g_{\delta})_x(x_0, y_0)p_{13} -(g_{\delta})_y(x_0, y_0)p_{23}\right]\\
&=&x'(s_0)[x'(s_0)^2 + y'(s_0)^2] + (g_{\delta})_x(x_0, y_0)x'(s_0)^2z'(s_0) + (g_{\delta})_y(x_0, y_0)x'(s_0)y'(s_0)z'(s_0) \\
&=&x'(s_0)[x'(s_0)^2 + y'(s_0)^2] + x'(s_0)z'(s_0)^2 \\
&=&x'(s_0)[x'(s_0)^2 + y'(s_0)^2 + z'(s_0)^2] = x'(s_0)\cdot 1 = x'(s_0).
\end{eqnarray*}
So
\begin{eqnarray*}
&&\frac{n!b_n(\delta)}{\sqrt{x'(s_0)^2 + y'(s_0)^2}}\\
&=&\sum_{a+b = 2}^{n} \frac{\partial^{a+b} g_{\delta}}{\partial x^{a}\partial y^b}(x_0, y_0) \sum_{I, J} c_{I, J} x'(s_0)^{a-|I|}y'(s_0)^{b-|J|} (p_{13})^{|I|}(p_{23})^{|J|}\\
&&(2!)^{i_1+j_1}(3!)^{i_2+j_2} \cdots (p+1)!^{i_p + j_p}b_2(\delta)^{i_1+j_1}\cdots b_{p+1}(\delta)^{i_p+j_p}\\
&=& \sum_{a+b = 2}^{n} \frac{\partial^{a+b} g_{\delta}}{\partial x^{a}\partial y^b}(x_0, y_0) \sum_{I, J} c_{I, J} x'(s_0)^{a-|I|}y'(s_0)^{b-|J|}\\
&&\left(\frac{-x'(s_0)^2z'(s_0)}{x'(s_0)\sqrt{x'(s_0)^2 + y'(s_0)^2}}\right)^{|I|}\left(\frac{-x'(s_0)y'(s_0)z'(s_0)}{x'(s_0)\sqrt{x'(s_0)^2 + y'(s_0)^2}}\right)^{|J|}\\
&&(2!)^{i_1+j_1}(3!)^{i_2+j_2} \cdots (p+1)!^{i_p + j_p}b_2(\delta)^{i_1+j_1}\cdots b_{p+1}(\delta)^{i_p+j_p}\\
&=&\sum_{a+b = 2}^{n} \frac{\partial^{a+b} g_{\delta}}{\partial x^{a}\partial y^b}(x_0, y_0) \sum_{I, J} c_{I, J} (-1)^{|I| + |J|}x'(s_0)^ay'(s_0)^b\left(\frac{z'(s_0)}{\sqrt{x'(s_0)^2 + y'(s_0)^2}}\right)^{|I|+|J|}\\
&&(2!)^{i_1+j_1}(3!)^{i_2+j_2} \cdots (p+1)!^{i_p + j_p}b_2(\delta)^{i_1+j_1}\cdots b_{p+1}(\delta)^{i_p+j_p}.
\end{eqnarray*}

\begin{lemma}\label{lem5}
The signs of $b_n(\delta)$ for $n$ between 2 and $N-1$ are all negative. 
\end{lemma}

\begin{proof}
Before proving the lemma, one needs to estimate $\frac{\partial^{p} g_{\delta}}{\partial x^{p}}(x_0, y_0)$ for $2 \leq p \leq N-1$. By induction one can show that 
\begin{eqnarray*}
\frac{\partial^{p} g_{\delta}}{\partial x^{p}}(x_0, y_0) &=& p! a_p(\delta) + \cdots + (N-1)(N-2)\cdots(N-p)a_{N-1}(\delta)x_0^{N-1-p} \\
&& + \sum_{q=0}^p \left(\begin{matrix}p \\ q \end{matrix}\right)N(N-1)\cdots(N-q+1)x_0^{N-q}\partial_x^{p-q}h_{\delta}(x_0, y_0) \\
&&+ y_0 \left[\sum_{q=0}^p \left(\begin{matrix}p \\ q \end{matrix}\right) \frac{d^q x}{d x^q}(x_0) \partial_x^{p-q}i_{\delta}(x_0, y_0)\right] + y_0^2 \partial_x^p j_{\delta}(x_0, y_0).
\end{eqnarray*}
First if $\epsilon$ is sufficiently small, 
\begin{eqnarray*}
&&|p! a_p(\delta) + \cdots + (N-1)(N-2)\cdots(N-p)a_{N-1}(\delta)x_0^{N-1-p}| \\
&\leq& M|\sin \delta|[p! + \cdots + (N-1)(N-2)\cdots(N-p)x_0^{N-1-p}] \\
&\leq& \frac{c}{2} b\tan(\theta_0)|\sin \delta| \leq \frac{c}{2}a_2(\delta).
\end{eqnarray*}
Second if $\epsilon$ and $\eta$ are sufficiently small, 
\begin{eqnarray*}
&&\sum_{q=0}^p \left(\begin{matrix}p \\ q \end{matrix}\right)N(N-1)\cdots(N-q+1)x_0^{N-q}\partial_x^{p-q}h_{\delta}(x_0, y_0)\\
&\leq& N(N-1)\cdots(N-p+1)x_0^{N-p}h(0, 0)(1-\frac{c}{2}).
\end{eqnarray*}
Third if $\epsilon$ is sufficiently small,
\begin{eqnarray*}
&&\left| y_0 \left[\sum_{q=0}^p \left(\begin{matrix}p \\ q \end{matrix}\right) \frac{d^q x}{d x^q}(x_0) \partial_x^{p-q}i_{\delta}(x_0, y_0)\right] + y_0^2 \partial_x^p j_{\delta}(x_0, y_0)\right| \\
&\leq& |y_0|C_1 \leq (18a_2(\delta)x_0^2 + 18N|h(0, 0)|x_0^N)ACC_1\\
&\leq&\frac{c}{2}a_2(\delta) + N(N-1)\cdots(N-p+1)x_0^{N-p}|h(0, 0)|\frac{c}{2}.
\end{eqnarray*}
Combining the above three inequalities, together with (\ref{x_0N-22}), yields
\begin{eqnarray*}
\frac{\partial^{p} g_{\delta}}{\partial x^{p}}(x_0, y_0) &\leq& ca_2(\delta) + N(N-1)\cdots(N-p+1)x_0^{N-p}h(0, 0)(1-c) \\
&\leq& ca_2(\delta) -\frac{N(N-1)\cdots(N-p+1)|h(0, 0)|(1-c)^2a_2(\delta)}{N(N-1)(1+c)|h(0, 0)|}\\
&=& -\left[\frac{(N-2) \cdots (N-p+1)(1-c)^2}{1+c} -c\right]a_2(\delta)\\
&\leq& -\left[\frac{(1-c)^2}{1+c}-c\right]a_2(\delta),
\end{eqnarray*}
which is negative if we choose $c$ as follows:
\[\frac{(1-c)^2}{1+c}-c > 0 \Rightarrow (1-c)^2 > c(1+c) \Rightarrow0 < c < \frac{1}{3}.\]
Denote the constant in the brackets as $L=L(c)$, then for $2 \leq p \leq N-1$
\[\frac{\partial^{p} g_{\delta}}{\partial x^{p}}(x_0, y_0) \leq -La_2(\delta).\]


Now we are ready to prove the lemma by induction. When $n=2$,
\[\frac{2b_2(\delta)}{\sqrt{x'(s_0)^2 + y'(s_0)^2}} = (g_{\delta})_{xx}(x_0, y_0)x'(s_0)^2 + 2(g_{\delta})_{xy}(x_0, y_0)x'(s_0)y'(s_0) +(g_{\delta})_{yy}(x_0, y_0)y'(s_0)^2.\]
On the one hand, since $x'(s_0) \geq \frac{1}{2}$,
\[(g_{\delta})_{xx}(x_0, y_0)x'(s_0)^2 \leq -\frac{L}{4}a_2(\delta).\]
On the other hand, 
\begin{eqnarray*}
&&y'(s_0) \left[2(g_{\delta})_{xy}(x_0, y_0)x'(s_0) +(g_{\delta})_{yy}(x_0, y_0)y'(s_0)\right] \leq C_1|y'(s_0)| \\
&\leq& C_1(4a_2(\delta)x_0 + 4Nx_0^{N-1}|h(0, 0)|)AC \\
&\leq& 4ACC_1x_0a_2(\delta) + 4ACC_1x_0 \frac{(3+c)a_2(\delta)N|h(0, 0)|}{N(N-1)|h(0, 0)|(1-c)}\\
&=&4ACC_1x_0a_2(\delta)\left[1 + \frac{3+c}{(N-1)(1-c)}\right].
\end{eqnarray*}
If we choose $\epsilon$ small enough so that
\[4ACC_1x_0\left[1 + \frac{3+c}{(N-1)(1-c)}\right] \leq \frac{L}{8},\] then
\[2(g_{\delta})_{xy}(x_0, y_0)x'(s_0)y'(s_0) +(g_{\delta})_{yy}(x_0, y_0)y'(s_0)^2 \leq  \frac{L}{8}a_2(\delta).\]
It follows that
\[\frac{2b_2(\delta)}{\sqrt{x'(s_0)^2 + y'(s_0)^2}} \leq -\frac{L}{8}a_2(\delta) \Rightarrow b_2(\delta) < 0.\]
By inductive hypothesis, suppose $b_2(\delta), \ldots, b_{n-1}(\delta)$ are all negative, then
it suffices to show that 
\[\frac{n!b_n(\delta)}{\sqrt{x'(s_0)^2 + y'(s_0)^2}} < 0.\]
There are two cases. 
Case 1: when $b=0$, $|J| = 0$ since $|J| \leq b$, then we have
\begin{eqnarray*}
&&\sum_{a=2}^n \frac{\partial^a g_{\delta}}{\partial x^a}(x_0, y_0) \sum_{I, J}c_{I, J}(-1)^{|I|}x'(s_0)^a\left(\frac{z'(s_0)}{\sqrt{x'(s_0)^2 + y'(s_0)^2}} \right)^{|I|} \\
&&2!^{i_1}3!^{i_2}\cdots(p+1)!^{i_p}b_2(\delta)^{i_1}b_3(\delta)^{i_2}\cdots b_{p+1}(\delta)^{i_p}.
\end{eqnarray*}
Since $z'(s_0) > 0$ and the sign of $b_2(\delta)^{i_1}b_3(\delta)^{i_2}\cdots b_{p+1}(\delta)^{i_p}$ is $(-1)^{|I|}$, it follows that for each $I$, $J$,
\[c_{I, J}(-1)^{|I|}x'(s_0)^a\left(\frac{z'(s_0)}{\sqrt{x'(s_0)^2 + y'(s_0)^2}} \right)^{|I|} 2!^{i_1}3!^{i_2}\cdots(p+1)!^{i_p}b_2(\delta)^{i_1}b_3(\delta)^{i_2}\cdots b_{p+1}(\delta)^{i_p} \geq 0.\]
Since $\partial_x^a g_{\delta} (x_0, y_0) < 0$, the above sum is negative. Especially when $a=n$, $p=n-a-b=0$ and so there is only one term
\[\frac{\partial^n g_{\delta}}{\partial x^n}(x_0, y_0)x'(s_0)^n \leq -\frac{1}{2^n}La_2(\delta).\]
whose coefficient is 1 by Corollary \ref{cor4}.

Case 2: when $b\neq 0$, there is at least one copy of $y'(s_0)$ in the summation, so we can write the rest of the terms as 
\begin{eqnarray*}
&&\left| y'(s_0) \left[\sum_{ a+b = 2, b \geq 1}^{n} \frac{\partial^{a+b} g_{\delta}}{\partial x^{a}\partial y^b}(x_0, y_0) \sum_{I, J} c_{I, J} (-1)^{|I| + |J|}x'(s_0)^ay'(s_0)^{b-1}\left(\frac{z'(s_0)}{\sqrt{x'(s_0)^2 + y'(s_0)^2}}\right)^{|I|+|J|}\right. \right.\\
&&\left. \left. (2!)^{i_1+j_1}(3!)^{i_2+j_2} \cdots (p+1)!^{i_p + j_p}b_2(\delta)^{i_1+j_1}\cdots b_{p+1}(\delta)^{i_p+j_p}\right] \right|\\
&\leq&C_1|y'(s_0)| \leq C_1(4a_2(\delta)x_0 + 4Nx_0^{N-1}|h(0, 0)|)AC\\
&\leq&4ACC_1x_0a_2(\delta)\left[1+ \frac{3+c}{(N-1)(1-c)}\right].
\end{eqnarray*}
where the first inequality is because everything inside the brackets is bounded. 
If we choose $\epsilon$ small enough so that 
\[4ACC_1x_0\left[1+ \frac{3+c}{(N-1)(1-c)}\right] \leq \frac{1}{2^{n+1}}L,\]
then
\[\frac{n!b_n(\delta)}{\sqrt{x'(s_0)^2 + y'(s_0)^2}} \leq -\frac{1}{2^{n+1}}La_2(\delta) < 0,\]
as desired.
\end{proof} 

The next lemma not only determines the sign of $l_{\delta}(0, 0)$, but also gives an upper bound of $l_{\delta}(0, 0)$.
\begin{lemma}\label{lem6}
\[l_{\delta}(0, 0) \leq \frac{1-c}{2^{N-1}}h(0, 0) < 0.\]
\end{lemma}

\begin{proof}
By Lemma \ref{lem4},
\[p_{33}\frac{\partial^N k_\delta}{\partial u^N}(u, v) = \sum_{a+b = 1}^{N} \frac{\partial^{a+b} g_{\delta}}{\partial x^{a}\partial y^b} \sum_{I, J} c_{I, J} A^{a-|I|}B^{b-|J|} (\frac{\partial A}{\partial u})^{i_1}\cdots (\frac{\partial^p A}{\partial u^p})^{i_p}(\frac{\partial B}{\partial u})^{j_1}\cdots (\frac{\partial^p B}{\partial u^p})^{j_p}.\]
When $a+b =1$ or $p=N-1$ Corollary \ref{cor3} suggests that we have in the above sum
\[(g_{\delta})_x \frac{\partial^{N-1} A}{\partial u^{N-1}} + (g_{\delta})_y \frac{\partial^{N-1} B}{\partial u^{N-1}}.\]
Evaluating at $(u, v)=(0, 0)$ gives us
\begin{eqnarray*}
&&p_{33}N!l_{\delta}(0, 0) - (g_{\delta})_x(x_0, y_0)p_{13}N!l_{\delta}(0, 0) - (g_{\delta})_y(x_0, y_0)p_{13}N!l_{\delta}(0, 0) \\
&=& \sum_{a+b = 2}^{N} \frac{\partial^{a+b} g_{\delta}}{\partial x^{a}\partial y^b}(x_0, y_0) \sum_{I, J} c_{I, J} x'(s_0)^{a-|I|}y'(s_0)^{b-|J|} (p_{13})^{|I|}(p_{23})^{|J|}\\
&&  (2!)^{i_1+j_1}(3!)^{i_2+j_2} \cdots (p+1)!^{i_p + j_p}b_2(\delta)^{i_1+j_1}\cdots b_{p+1}(\delta)^{i_p+j_p}
\end{eqnarray*}
Since $p_{33} = Zx'(s_0)$, $p_{13} = -x'(s_0)^2z'(s_0)Z$, and $p_{23} = -x'(s_0)y'(s_0)z'(s_0) Z$, one has
\begin{eqnarray*}
&&\frac{N!l_{\delta}(0, 0)}{\sqrt{x'(s_0)^2 + y'(s_0)^2}}\\
&=&\sum_{a+b = 2}^{N} \frac{\partial^{a+b} g_{\delta}}{\partial x^{a}\partial y^b}(x_0, y_0) \sum_{I, J} c_{I, J} (-1)^{|I| + |J|}x'(s_0)^ay'(s_0)^b\left(\frac{z'(s_0)}{\sqrt{x'(s_0)^2 + y'(s_0)^2}}\right)^{|I|+|J|}\\
&&(2!)^{i_1+j_1}(3!)^{i_2+j_2} \cdots (p+1)!^{i_p + j_p}b_2(\delta)^{i_1+j_1}\cdots b_{p+1}(\delta)^{i_p+j_p}.
\end{eqnarray*}
Next let's use an analogous argument in Lemma \ref{lem3.2}. When $b=0$, $|J|=0$, then we have
\begin{eqnarray*}
&&\sum_{a=2}^N \frac{\partial^a g_{\delta}}{\partial x^a}(x_0, y_0) \sum_{I, J}c_{I, J}(-1)^{|I|}x'(s_0)^a\left(\frac{z'(s_0)}{\sqrt{x'(s_0)^2 + y'(s_0)^2}} \right)^{|I|} \\
&&2!^{i_1}3!^{i_2}\cdots(p+1)!^{i_p}b_2(\delta)^{i_1}b_3(\delta)^{i_2}\cdots b_{p+1}(\delta)^{i_p}\\
&\leq& \frac{\partial^N g_\delta}{\partial x^N}(x_0, y_0) x'(s_0)^N \leq \frac{1}{2^N}\frac{\partial^N g_\delta}{\partial x^N}(x_0, y_0),
\end{eqnarray*}
Since
\begin{eqnarray*}
\frac{\partial^N g_\delta}{\partial x^N}(x_0, y_0) &=& \sum_{q=0}^N \left( \begin{matrix} N \\ q \end{matrix} \right)N(N-1)\cdots(N-q+1)x_0^{N-q}\partial_x^{q}h_\delta(x_0, y_0) \\
&& + y_0 \left[ \sum_{q=0}^N \left( \begin{matrix} N \\ q \end{matrix} \right) \frac{\partial^{N-q} x}{\partial x^{N-q}}(x_0) \partial_x^q i_\delta(x_0, y_0) + y_0\frac{\partial^N j_\delta}{\partial x^N}(x_0, y_0) \right].
\end{eqnarray*}
On the one hand, if $\eta$ and $\epsilon$ are sufficiently small, 
\begin{eqnarray*}
\sum_{q=0}^N \left( \begin{matrix} N \\ q \end{matrix} \right)N(N-1)\cdots(N-q+1)x_0^{N-q}\partial_x^{q}h_\delta(x_0, y_0) \leq N!(1-\frac{c}{4})h(0, 0).
\end{eqnarray*}
On the other hand, 
\begin{eqnarray*}
&&\left|  y_0 \left[ \sum_{q=0}^N \left( \begin{matrix} N \\ q \end{matrix} \right) \frac{\partial^{N-q} x}{\partial x^{N-q}}(x_0) \partial_x^q i_\delta(x_0, y_0) + y_0\frac{\partial^N j_\delta}{\partial x^N}(x_0, y_0) \right] \right| \leq C_1|y_0| \\
&\leq& (18a_2(\delta)x_0^2 + 18N|h(0, 0)|x_0^N)ACC_1 \leq \frac{c}{4}N!|h(0, 0)|,
\end{eqnarray*}
for $\eta$ and $\epsilon$ sufficiently small. Thus
\[\frac{\partial^N g_\delta}{\partial x^N}(x_0, y_0) \leq N!(1-\frac{c}{2})h(0, 0).\]
When $b \neq 0$, there is always one copy of $y'(s_0)$ so we have 
\begin{eqnarray*}
&&\left|y'(s_0)\left[ \sum_{a+b = 2}^{N} \frac{\partial^{a+b} g_{\delta}}{\partial x^{a}\partial y^b}(x_0, y_0) \sum_{I, J} c_{I, J} (-1)^{|I| + |J|}x'(s_0)^ay'(s_0)^{b-1}\left(\frac{z'(s_0)}{\sqrt{x'(s_0)^2 + y'(s_0)^2}}\right)^{|I|+|J|} \right. \right.\\
&&\left. \left. (2!)^{i_1+j_1}(3!)^{i_2+j_2} \cdots (p+1)!^{i_p + j_p}b_2(\delta)^{i_1+j_1}\cdots b_{p+1}(\delta)^{i_p+j_p}\right] \right| \\
&\leq& C_1 |y'(s_0)| \leq C_1(4a_2(\delta)x_0 + 4Nx_0^{N-1}|h(0, 0)|)AC\\
&\leq&\frac{N!}{2^N}\frac{c}{2}|h(0, 0)|,
\end{eqnarray*}
for $\eta$ and $\epsilon$ sufficiently small. It follows that
\begin{eqnarray*}
&& \frac{N!l_{\delta}(0, 0)}{\sqrt{x'(s_0)^2 + y'(s_0)^2}} \leq \frac{N!}{2^N}(1-\frac{c}{2})h(0, 0) - \frac{N!}{2^N}\frac{c}{2}h(0, 0) = \frac{N!}{2^N}(1-c)h(0, 0) \\
&\Rightarrow& l_{\delta}(0, 0) \leq \frac{\sqrt{x'(s_0)^2 + y'(s_0)^2}}{2^N}(1-c)h(0, 0) \leq \frac{2}{2^N}(1-c)h(0, 0) < 0.
\end{eqnarray*}
\end{proof}

(2). Since $b_2(\delta) < 0$, we could show as in Case 2 that $\gamma$ is initially a straight line. The surface $S$ in the $(u, w)$-plane is the curve
\[w=k_{\delta}(u, 0) = b_2(\delta)u^2 + \cdots + b_{N-1}(\delta)u^{N-1} + u^Nl_{\delta}(u, 0).\]
If the line segment re-enters the surface at some switch point, then the curve can't be concave downward there. Otherwise the line lies above the surface.

\begin{eqnarray*}
w'(u) &=&2b_2(\delta)u + \cdots + (N-1)b_{N-1}(\delta)u^{N-2} + Nu^{N-1}l_{\delta}(u, 0),\\
w''(u) &=&2b_2(\delta) + 6b_3(\delta) + \cdots + (N-1)(N-2)b_{N-1}(\delta)u^{N-3} \\
&& + N(N-1)u^{N-2}l_{\delta}(u, 0) + 2Nu^{N-1}(l_{\delta})_u(u, 0) + u^N(l_{\delta})_{uu}(u, 0).
\end{eqnarray*}
On the one hand, 
\[2b_2(\delta) + 6b_3(\delta) + \cdots + (N-1)(N-2)b_{N-1}(\delta)u^{N-3} < 0.\]
On the other hand, we can choose $\epsilon$ small enough so that for all $u < \epsilon$,
\begin{eqnarray*}
&&l_{\delta}(u, 0) \leq (1-\frac{c}{2})l_{\delta}(0, 0), \frac{|2Nu(l_{\delta})_u(u, 0) + u^2(l_{\delta})_{uu}(u, 0)|}{N(N-1)} \leq \frac{c}{2}|l_{\delta}(0, 0)|\\
&\Rightarrow& N(N-1)u^{N-2}l_{\delta}(u, 0) + 2Nu^{N-1}(l_{\delta})_u(u, 0) + u^N(l_{\delta})_{uu}(u, 0) \\
&&\leq N(N-1)u^{N-2}(1-c)l_{\delta}(0, 0) < 0.
\end{eqnarray*}
So $w''(u) < 0$ and the graph is concave downward. Therefore $\gamma$ never re-enters the surface at a switch point. It follows that $\gamma$ is a straight line that either terminates at some point on the surface or exits the $\epsilon$-ball.

Let's summarize Case 4. In general, $\gamma$ is initially a boundary segment lying on the surface, then it leaves $S$ in a straight line that exits the $\epsilon$-ball. So there is at most one interval in this case. 

In the end, we are going to mention when the lowest degree in the Taylor expansion of $g(x, y)$ is greater than 2. Suppose the lowest degree is $k \geq 2$, the kth Taylor polynomial has the following form:
\[a_0x^k + a_1x^{k-1}y + \cdots + a_{k-1}xy^{k-1} + a_ky^k\]
Consider the line $y=mx$, substituing it into the above polynomial gives us
\[x^k(a_0 + a_1m + \cdots + a_km^k).\]
Setting it to zero gives as at most $k$ distinct solutions for $m$ if $a_k\neq0$, otherwise adding the vertical line $x=0$. Therefore the plane can be sliced into at most $2k$ distinct pies using these slopes, such that within each slice, the graph of $g$ along a ray is either concave upward or downward near the origin. The rest of the proof is analogous to the case when $k=2$.

\end{proof}

\end{section}

%% file: conclusion.tex
\begin{section}{Conclusion}
It seems naturally that our theorems could be generalized to higher-dimensional Euclidean space. However, the proof in Theorem 1 or 2 does not apply when $n > 3$, because the intersection of $M_1$ and $M_2$ becomes a surface instead of a curve. Furthermore, the proof in Theorem 3 involves dividing the plane into finitely many slices using the lowest degree Taylor polynomial, which does not make sense when we have more than two variables. Therefore we are looking for new methods and we conjecture that all Theorems 1, 2, and 3 do generalize when $n > 3$ despite of more free variables. Furthermore, the idea for two obstacles can be naturally generalized to finitely many obstacles. 
Moreover, this paper provides some good tools for the 3-dimensional stratification problem for finiteness property in a semi-algebraic set. 
\end{section}